\documentclass[reqno]{amsart}
\usepackage[utf8]{inputenc}

\usepackage{ amssymb, nccmath}

\usepackage{amsmath}

\usepackage{enumerate}
\usepackage{latexsym,amscd,amsfonts}
\usepackage{graphicx}
\usepackage{graphics}
\usepackage{tikz,tkz-tab}
\usepackage{float}

\def \Z {\mathbb Z}
\def \R {\mathbb R}

\newcommand{\ex}{\mathrm{e}}
\newcommand{\ii}{\mathrm{i}}
\newcommand{\dd}{\mathrm{d}}

\newtheorem{lem}{Lemma}

\newtheorem{theo}{Theorem}
\newtheorem{prop}{Proposition}
\newtheorem{cor}{Corollary}
\newenvironment{demo}{\noindent \textit{Proof:}}{\begin{flushright} $\Box$ \end{flushright}}

\title[The periodic Airy-Schr\"odinger operator]{The band spectrum of the periodic Airy-Schr\"odinger operator on the real line}

\author[H. Boumaza]{Hakim Boumaza}
\email{boumaza@math.univ-paris13.fr}
\address{Universit\'e Paris 13, Sorbonne Paris Cit\'e\\
LAGA\\
CNRS, UMR 7539\\
F-93430, Villetaneuse\\
France}

\author[O. Lafitte]{Olivier Lafitte}
\email{lafitte@math.univ-paris13.fr}
\address{Universit\'e Paris 13, Sorbonne Paris Cit\'e\\
LAGA\\
CNRS, UMR 7539\\
F-93430, Villetaneuse\\
France}
\thanks{CEA/DEN/DM2S, F-91191 Gif sur Yvette Cedex}


\dedicatory{A tribute to Louis Boutet de Monvel (1941-2014)}

\begin{document}


\begin{abstract}
We introduce the periodic Airy-Schr\"odinger operator and we study its band spectrum. This is an example of an explicitly solvable model with a periodic potential which is not differentiable at its minima and maxima. We define a semiclassical regime in which the results are stated for a fixed value of the semiclassical parameter and are thus estimates instead of asymptotic results. We prove that there exists a sequence of explicit constants, which are zeroes of classical functions, giving upper bounds of the semiclassical parameter for which the spectral bands are in the semiclassical regime. We completely determine the behaviour of the edges of the first spectral band with respect to the semiclassical parameter. Then, we investigate the spectral bands and gaps situated in the range of the potential. We prove precise estimates on the widths of these spectral bands and these spectral gaps and we determine an upper bound on the integrated spectral density in this range. Finally, in the semiclassical regime, we get  estimates of the edges of every spectral bands and thus of the widths of every spectral bands and spectral gaps.
\end{abstract}


\maketitle
\vspace{.5cm}


\section{The model and the semiclassical parameter}\label{sec_model_can_sol}
\vskip 2mm

\noindent For periodic Schr\"odinger operators, the semiclassical behaviour of the bottom of the spectrum and of the widths of the spectral bands and gaps is very well known when the potential is analytic or $C^{\infty}$. We introduce a periodic Schr\"odinger operator whose potential is not differentiable at its minima and maxima points and for which we still determine the semiclassical behaviour of the band edges and of the widths of the spectral bands and gaps. 
\vskip 2mm

\noindent The spectrum of the Schr\"odinger operator for a linear potential well is known. It is pure point spectrum and the eigenvalues are given by the absolute values of the zeroes of the first Airy function and its derivative. Our periodic Schr\"odinger operator with singular potential is a periodization of the truncated linear potential well. We expect to find a spectrum made of spectral bands, each containing one of the eigenvalues of the Schr\"odinger operator with a linear potential well.

\vskip 2mm

\subsection{The periodic Airy-Schr\"odinger operator} Let $2L_0\in \R_+^*$ be a characteristic length modelling the distance between two ions in a one dimensional periodic lattice of ions. The motion of electrons of mass $m$ in this lattice is studied through the following $2L_0$-periodic Schr\"odinger operator acting on the Sobolev space $\mathsf{H}^2(\mathbb{R})$, 
\begin{equation}\label{eq_def_H}
H= -\frac{\hbar^2}{2m} \frac{\mathrm{d}^2 }{\mathrm{d}x^2} + V, 
\end{equation}
where $\hbar$ is the reduced Planck constant and $V$ is the $2L_0$-periodic function on $\R$ defined by 
$$\forall x \in [-L_0,L_0],\ V(x) =V_0\left( \mfrac{|x|}{L_0} -1\right),$$
$V_0\in \R_+^*$ being a reference potential. The ions, in this model, are located at points $2nL_0$ for $n\in \Z$, this points corresponding to the minima of the potential $V$. 
\vskip 3mm

We call $H$ the \emph{periodic Airy-Schr\"odinger operator} on $\R$.

\vskip 3mm

\noindent Since $V$ is periodic and locally integrable, the theory of periodic operators (\cite{RS4}) asserts that the operator $H$ has purely absolutely continuous spectrum and that this spectrum is the union of spectral bands:  
$$\sigma(H)= \bigcup_{p\geq 0} \left[ E_{\mathrm{min}}^p, E_{\mathrm{max}}^p \right],$$
where $E_{\mathrm{min}}^p$ and $E_{\mathrm{max}}^p$ are the spectral band edges. For $p\geq 0$, we shall call $[ E_{\mathrm{min}}^p, E_{\mathrm{max}}^p ]$ the $p$-th spectral band and $(E_{\mathrm{max}}^p,E_{\mathrm{min}}^{p+1} )$ the $p$-th spectral gap. We will precise these notations and characterize these spectral band edges in Section \ref{sec_band_spec}.
\vskip 3mm

\noindent Traditional results describe the spectral bands near the minimum of a $C^2$ potential (\cite{H}). We generalize this analysis to the case of the potential $V$ which is not differentiable on its minima points. We are able to count the number of spectral bands in $[-V_0,0]$, the range of the potential $V$, for any value of a dimensionless parameter defined in (\ref{def_param}), and to describe precisely these spectral bands.
\vskip 3mm

\noindent This dimensionless parameter is called in the sequel, the semiclassical parameter.
\vskip 2mm

\subsection{The semiclassical parameter} In order to study the spectrum of $H$, one consider the equation in $\psi \in \mathsf{H}^2(\R)$,
\begin{equation}\label{eq_schrodinger}
 -\frac{\hbar^2}{2m}\psi''+V(x)\psi=E\psi\ .
\end{equation}
As $V$ is affine on every interval $[nL_0, (n+1)L_0]$ for every $n\in \Z$, one recognizes in (\ref{eq_schrodinger}), after rescaling and translating the variable, the Airy equation: $u''=xu$.
\vskip 3mm

\noindent On the interval $[0,L_0]$, the rescaling is done through a parameter $\theta$ given as the unique real number such that 
$$-\frac{\hbar^2}{2m}\theta^2\frac{V_0^2}{L_0^2}(\theta(V(x)-E))+V(x)-E=0,$$ 
namely
\begin{equation}\label{def_param_0}
\theta=\left(\frac{2mL_0^2}{\hbar^2V_0^2} \right)^{\frac13}.
\end{equation}
\vskip 2mm

\noindent \textbf{Notation.} For any real number $E$, we set $\mathbf{E}=\theta E$. In particular, we set, for any $p\geq 0$, $\mathbf{E}_{\mathrm{min}}^p=\theta E_{\mathrm{min}}^p$ and $\mathbf{E}_{\mathrm{max}}^p=\theta E_{\mathrm{max}}^p$. We also denote by $\mathbf{V}$ the function $x \mapsto \theta V(x)$.

\vskip 3mm

\noindent Note that $\mathbf{V}(x)$ is dimensionless for all values of $x\in \R$ and thus we introduce the semiclassical parameter $\mathsf{h}$ and the counting parameter $\mathsf{c}$: 
\begin{equation}\label{def_param}
\mathsf{h}:=\frac{\hbar}{L_0(2mV_0)^{\frac12}}\quad \mbox{ and }\quad \mathsf{c}:=\theta V_0=\left(\frac{2m V_0 L_0^2}{\hbar^2} \right)^{\frac13} =\mathsf{h}^{-\frac{2}{3}}.
\end{equation}

\vskip 3mm

\noindent We use the semiclassical parameter $\mathsf{h}$ to rewrite the periodic Airy-Schr\"odinger operator in a form analog to the operator studied in \cite{H}. Then, the equation (\ref{eq_schrodinger}) is equivalent to 
\begin{equation}\label{eq_H_Harrell_form}
-\mathsf{h}^2 \frac{\dd^2 }{\dd z^2} \phi + q(z) \phi =\mathsf{h}^{\frac23} \mathbf{E} \phi 
\end{equation}
with $q$ the $2$-periodic function equal to $q(z)=|z| - 1$ on the interval $[-1,1]$. Setting $z=\mathsf{h}^{\frac23}y$, the equation (\ref{eq_schrodinger}) is also equivalent to
\begin{equation}\label{eq_H_classic_form}
- \frac{\dd^2 }{\dd y^2} \phi + (q(y)+1) \phi =( \mathbf{E}+\mathsf{h}^{-\frac23}) \phi.  
\end{equation}
\vskip 3mm

\noindent The parameter $\mathsf{c}$ is the natural dimensionless parameter which allows to count the spectral bands of $H$ in the range of $V$. 
\vskip 2mm
\noindent The semiclassical parameter $\mathsf{h}$ allows to define a semiclassical limit and a notion of ``semiclassical regime'' in which we study the spectral properties of the periodic Airy-Schr\"odinger operator. 
\vskip 2mm

\noindent In the sequel, we call the \textit{semiclassical limit}, the limit when the semiclassical parameter $\mathsf{h}$ tends to $0$.
\vskip 2mm

\noindent A property $(\mathcal{P}_{\mathsf{h}})$ of $H$ depending on $\mathsf{h}$ is said to be true in the \textit{semiclassical regime}, if:
$$ \exists \mathsf{h}_0 >0,\ \forall \mathsf{h}\in (0,\mathsf{h}_0), \ (\mathcal{P}_{\mathsf{h}}) \ \mbox{ is true.} $$
The interval $(0,\mathsf{h}_0)$ is called an interval of validity of the semiclassical approximation for the property $(\mathcal{P}_{\mathsf{h}})$.
\vskip 3mm

\noindent This idea of a semiclassical regime follows the discussion of \cite{HL} about the domain of validity of the semiclassical analysis for their model. One of the purposes of semiclassical analysis is to prove estimates for asymptotic expansions, when the semiclassical parameter tends to $0$, of quantities depending on it. If these asymptotic expansions are actually estimates valid for any ``small'' value of the semiclassical parameter (\emph{i.e.} for $\mathsf{h}\in (0,\mathsf{h}_0)$ fixed), we want to say that the considered quantity is in the semiclassical regime.

\noindent In the literature, many results of semiclassical analysis are stated in the setting of the semiclassical regime, providing estimates uniform on $\mathsf{h}\in (0,\mathsf{h}_0)$. For example, the results for Schr\"odinger operators with multiple wells (\cite{HS1}) are written in this form in \cite{R}. The semiclassical theory of the Harper model (\cite{HS2,HS3}) also provides estimates written for the semiclassical regime. For periodic Schr\"odinger operators with regular potentials, the results in \cite{O} are also stated in the semiclassical regime. Usually, for these results, the constant $\mathsf{h}_0$ is not given explicitely, but sometimes it is (see \cite{HMR} for results with an explicit $\mathsf{h}_0$). For general references about semiclassical analysis we refer to the textbooks \cite{DS,martinez,Z}.
\vskip 2mm

\noindent In the present article, we are able to give explicit values of $\mathsf{h}_0$ for which the spectral band edges and the widths of the spectral bands and gaps are in the semiclassical regime. For example, we prove that the bottom of the spectrum of $H$ is in the semiclassical regime whenever $\mathsf{h}$ is smaller than a power of the absolute value of the first zero of the derivative of the Airy function. Thus, for the bottom of the spectrum, $\mathsf{h}$ ``small'' means $\mathsf{h}$ smaller than approximately $0.973$ hence is quite large.
\vskip 3mm

\noindent Our study of the periodic Airy-Schr\"odinger operator was first motivated by giving a rigourous mathematical treatment of the numerical results in \cite{Castro}, a paper which deals with the question of quark confinement. It was also motivated by previous results of one of the authors on the semiclassical analysis of the Rayleigh-Taylor instability \cite{CLa,L}. Note that the periodic Airy-Schr\"odinger model is also closely related with an infinite periodized sawtooth junction PN in photonic crystals, giving an explicitly solvable situation in a quantum setting (see \cite{morozov}).

\vskip 2mm
\noindent The periodic Airy-Schr\"odinger operator gives an example of operator for which one is able to give explicit values of the semiclassical parameter for which considered spectral quantities are in the semiclassical regime. These special values of the semiclassical parameter are given through the zeroes of the Airy function and its derivative and the zeroes of the derivative of the odd canonical solution of the Airy equation.
\vskip 2mm

\noindent Moreover, we identified two sequences $\{c_p\}_{p\geq 0}$ and $\{\tilde{c}_p\}_{p\geq 0}$ of values of the counting parameter $\mathsf{c}$ such that $c_k < \mathsf{c} < \tilde{c}_k$ or $\tilde{c}_k \leq \mathsf{c} < c_{k+1}$ if and only if there are exactly $k+1$ spectral bands in the range of the potential. This picture of the counting of energy levels in the range of the potential can be found in \cite{F}, in the simple case of the rectangular potential hole, which is not periodic. In \cite{F}, explicit values of the ``size parameter'' for which an eigenvalue enters the range of the potential are not given.

\vskip 6mm

\section{Notations and main results}\label{sec_main_results}

\subsection{The canonical solutions of the Airy equation} Let $u$ and $v$ be the canonical solutions of the Airy equation, satisfying
$$u(0)=1,\ u'(0)=0\quad \mbox{ and }\quad v(0)=0,\ v'(0)=1,$$
which Wronskian is $1$.
 The Airy function $Ai$ plays a very special role for the ordinary differential equation $u''=xu$. It generates the unique family of subdominant solutions of this Sturmian equation. It is also the unique solution of the Airy equation which is in $\mathcal{S}'(\R)$ and such that its Fourier transform satisfies $\widehat{Ai}(0)=1$. The function $Bi$ is the solution of the Airy equation satisfying the initial conditions $Bi(0)=\sqrt{3} Ai(0)$ and $Bi'(0)=-\sqrt{3} Ai'(0)$. One has the expression of $u$ and $v$ in terms of the classical Airy functions $Ai$ and $Bi$:
\begin{equation}\label{eq_expr_u_Ai}
\forall x \in \R,\ u(x)=\pi(Bi'(0)Ai(x)-Ai'(0)Bi(x))
\end{equation}
and
\begin{equation}\label{eq_expr_v_Ai}
\forall x \in \R,\   v(x)=\pi(Ai(0)Bi(x)-Bi(0)Ai(x)).
\end{equation}

\noindent Both $u$ and $v$ are analytic functions on $\R$. Moreover, $u$ and $v$ are strictly increasing and positive on $(0,+\infty)$. Thus, the zeroes of $u$, $v$ and their derivatives are all non-positive real numbers. 
\vskip 2mm

\noindent \textbf{Notation.} We denote by 
\begin{enumerate}
\item $\{ -\tilde{c}_{2j} \}_{j \geq 0}$ the set of the zeroes of $u$, 
\item $\{ -\tilde{c}_{2j+1} \}_{j \geq 0}\cup \{ 0\}$ the set of the zeroes of $u'$,
\item $\{ -c_{2j+1} \}_{j \geq 0} \cup \{ 0\} $  the set of the zeroes of $v$,
\item $\{ -c_{2j} \}_{j \geq 0}$ the set of the zeroes of $v'$.
\end{enumerate}
This definition is precised in Section \ref{sec_ordering_zeroes}.
\vskip 2mm
\noindent An important property of these zeroes, which is proven in Corollary \ref{prop_sep_ck}, is: $$\forall p\geq 0,\ -\tilde{c}_p < -c_p.$$

\vskip 3mm

\noindent Approximate values of the $c_p$ and $\tilde{c}_p$ can be given. For example, 

$$c_0 \simeq 1.515,\qquad \tilde{c}_0 \simeq 1.986, \qquad c_1\simeq 2.666, \qquad \tilde{c}_1 \simeq 2.948.$$
For more approximate values of the $c_p$'s, see Appendix \ref{sec_app_variations}.
\vskip 5mm

\subsection{The first spectral band}

\noindent Our first result gives the initialization of the counting of the spectral bands which are included in the range of $\mathbf{V}$, $[-\mathsf{c},0]$. 

\begin{theo}\label{thm_intro_transition}
For $\mathsf{c} \in (0,c_0]$, there is no rescaled spectral gap in $[-\mathsf{c},0]$. The first rescaled spectral gap intersects $[-\mathsf{c},0]$ as soon as $\mathsf{c} > c_0$.
\end{theo}

\vskip 3mm

\noindent \textbf{Remark.} Theorem \ref{thm_intro_transition} asserts that: for $\mathsf{c} \leq c_0$, the only gap in $(-\infty, 0]$ is the "ground state gap" $(-\infty, E_{\mathrm{min}}^0)$. The first non trivial gap intersects $[-V_0,0]$ as soon as $\mathsf{c} > c_0$.
\vskip 3mm

\noindent We get precise estimates of the rescaled ground state $\mathbf{E}_{\mathrm{min}}^0 $ in both semiclassical regime and when $\mathsf{h}$ tends to infinity. Before stating them, we need to introduce notations for the zeroes of the Airy function $Ai$ and its derivative. 

\vskip 3mm

\noindent \textbf{Notation.} We denote by $\{-a_j\}_{j\geq 1}$ the set of the zeroes of $Ai$ and by $\{-\tilde{a}_j\}_{j\geq 1}$ the set of the zeroes of $Ai'$ where the real numbers $-a_j$ and $-\tilde{a}_j$ are arranged in decreasing order. These sets are both subsets of $(-\infty,0]$. Moreover, for every $j\geq 1$, $-a_j \in (-\tilde{a}_{j+1},-\tilde{a}_j)$.

\vskip 3mm

\noindent We set $\alpha= -\frac{Ai(0)}{Ai'(0)} >0$. The number $\alpha$ is the inverse of the slope at $0$ of the Airy function $Ai$. An approximate value of $\alpha$ is: $\alpha \simeq 1,372$.
\vskip 3mm

\begin{theo}\label{thm_bottom}
We have the following estimates on $\mathbf{E}_{\mathrm{min}}^0$:
\begin{enumerate}
 \item For every $\mathsf{c}>0$
\begin{equation}\label{eq_thm_bottom1}
      -\mathsf{c} < \mathbf{E}_{\mathrm{min}}^0 < \min\left(-\mfrac{ \mathsf{c}}{2},\ -\mathsf{c} + \tilde{a}_1 \right).
       \end{equation}

\item Let $\tau >0$ be arbitrary small. For every $ \mathsf{h} \in (0,(\tilde{a}_1+\tau)^{-\frac32}],$
\begin{equation}\label{eq_thm_bottom2}
\mathbf{E}_{\mathrm{min}}^0 = -\mathsf{h}^{-\frac23}+\tilde{a}_1 -\alpha \sqrt{3} \mfrac{(u'(-\tilde{a}_1))^2  }{\tilde{a}_1 } \ex^{-\frac43 \mathsf{h}^{-1} + 2\tilde{a}_1\mathsf{h}^{-\frac13}} \left(1  +  \mathcal{O} \left(\mathsf{h}^{\frac13} \right) \right). 
\end{equation}
\end{enumerate}
\end{theo}

\vskip 2mm

\noindent \textbf{Remark.} The estimate (\ref{eq_thm_bottom2}) is in the semiclassical regime and an explicit interval of validity of the semiclassical approximation is given by $(0,\tilde{a}_1^{-\frac32} )$. Note that $\tilde{a}_1^{-\frac32}\simeq 0.973$ is the value mentioned in the introduction.
\vskip 2mm

\noindent After multiplication by $\mathsf{h}^{\frac23}$, the first term in the estimate (\ref{eq_thm_bottom2}) is equal to $-1$ and corresponds to the minimum of the potential $q$ defined in (\ref{eq_H_Harrell_form}).

\vskip 3mm

\noindent Then, we determine the first rescaled spectral band in the limit when $\mathsf{h}$ tends to infinity.

\begin{theo}\label{thm_intro_quantum_1st_band}

\begin{enumerate}
 \item One has,
$$\lim_{\mathsf{h} \to +\infty} \mathsf{h}^{\frac23} \mathbf{E}_{\mathrm{min}}^0 = -\frac12\quad \mbox{ and }\quad \lim_{\mathsf{h} \to +\infty} \mathbf{E}_{\mathrm{max}}^0 = +\infty.$$

\item More precisely, when $\mathsf{h}$ tends to infinity,
\begin{equation}\label{eq_thm_bottom3}
\mathsf{h}^{\frac23} \mathbf{E}_{\mathrm{min}}^0 = -\frac12- \frac{1}{120} \frac{1}{\mathsf{h}^2 } +\mathcal{O}\left( \frac{1}{\mathsf{h}^4 } \right).
\end{equation}
\end{enumerate}

\end{theo}

\vskip 2mm
\noindent \textbf{Remark}. We rewrite the first point in Theorem \ref{thm_intro_quantum_1st_band}: when $V_0$ is fixed and $L_0$ tends to $0$,
$$\lim_{L_0 \to 0} E_{\mathrm{min}}^0 = -\frac{V_0}{2}\quad \mbox{ and }\quad \lim_{L_0 \to 0} E_{\mathrm{max}}^0 = +\infty.$$
Thus, when $L_0$ tends to $0$, which implies at $V_0$ fixed that $\mathsf{h}$ tends to infinity, the first spectral band tends to cover all the interval $\left[-\frac{V_0}{2},  +\infty \right)$. This corresponds to a model with infinitely closed atoms, where the tight-binding approximation is no longer relevant. Thus, the limit when $\mathsf{h}$ tends to infinity is not a ``physical limit''. 
\vskip 5mm

\subsection{Counting and estimates of the spectral bands in the range of $V$}
 
\noindent We have estimates on the widths of the spectral bands and the spectral gaps which are located in the range of $\mathbf{V}$.

\vskip 3mm
\noindent Two $\mathsf{c}$-dependent integers are of interest in this paper:
\begin{enumerate}
\item  the unique integer $p_0$ such that
\begin{equation}\label{def_p_0}
\tilde{c}_{p_0+1}-\tilde{c}_{p_0} < \mathsf{c} \leq \tilde{c}_{p_0}-\tilde{c}_{p_0-1}\quad \mbox{ when }\quad  \mathsf{c} < c_0; 
\end{equation}
\item the unique integer $k_0$ such that 
\begin{equation}\label{def_k_0}
 c_{k_0} < \mathsf{c} < \tilde{c}_{k_0} \quad \mbox{ or }\quad  \tilde{c}_{k_0} \leq \mathsf{c} < c_{k_0 +1} \quad \mbox{ when }\quad  \mathsf{c}>c_0.
\end{equation}
\end{enumerate}
\noindent Denote the integer part of a  real number $x$ by $[x]$. Due to the intervals in which we localize the real numbers $c_p$ and $\tilde{c}_p$ in Proposition \ref{prop_loc_ck}, one has $k_0=\left[\frac{4}{3\pi}\mathsf{c}^{\frac32} \right]=\left[\frac{4}{3\pi}\frac{1}{\mathsf{h}} \right]$ or $k_0=\left[\frac{4}{3\pi}\mathsf{c}^{\frac32} \right]-1$.
\vskip 3mm

\noindent Let $p\geq 0$ an integer and denote by $\delta_p$ the width of the rescaled $p$-th spectral band and by $\gamma_p$ the width of the rescaled $p$-th spectral gap with $\delta_{p} =  \mathbf{E}_{\mathrm{max}}^{p}- \mathbf{E}_{\mathrm{min}}^{p}$ and $\gamma_{p}  = \mathbf{E}_{\mathrm{min}}^{p+1}- \mathbf{E}_{\mathrm{max}}^{p}$.
\vskip 3mm

\noindent Let $I$ be the strictly decreasing function defined on $[1,+\infty)$ by 
$$\forall y \geq 1,\ I(y)= \left(\mfrac32 \right)^{\frac13}\frac{y^{\frac32}+1}{y^2+y+1}.$$

\begin{theo}\label{thm_intro_pbands}
Let $\mathsf{c} > c_0$ and $k_0=\left[\frac{4}{3\pi}\mathsf{c}^{\frac32} \right]$ or $\left[\frac{4}{3\pi}\mathsf{c}^{\frac32} \right]-1$ introduced in (\ref{def_k_0}). 
\begin{enumerate}
\item The $k_0 +1$ first rescaled spectral bands are included in the range of $\mathbf{V}$, $[-\mathsf{c}, 0]$.
\item One has, for every $ p\in \{2,\ldots, k_0\}$, 
 \begin{equation}\label{eq_thm_intro_bandes}
0 <  \delta_{p} \leq  \left( \frac{\pi}{3}+\frac{7}{3\pi}\frac{p+\frac13}{p(p+\frac23)} \right)  \left( \frac{3}{\pi} \right)^{\frac{1}{3}}\frac{1}{p^{\frac13}}, 
 \end{equation}
and for every $p\in \{2,\ldots, k_0-1\}$,
  \begin{equation}\label{eq_thm_intro_gaps}
 0< I\left(\left( \mfrac76 \right)^{\frac23} \right) \frac{2^{\frac13}\pi^{\frac23}  }{9} \frac{1}{(p+1)^{\frac13}} <  \gamma_{p}\leq  \left( \pi+\frac{7}{3\pi}\frac{p}{p^2 -1}\right) \left( \frac{3}{\pi} \right)^{\frac{1}{3}} \frac{1}{(p-1)^{\frac13}}.
 \end{equation}
In particular, all the gaps in $\sigma(H)\cap [-V_0,0]$ are open.
\end{enumerate}
\end{theo}
\vskip 3mm

\noindent Recall that the spectral gaps are said to be open when they are not empty.
\vskip 3mm

\noindent The number of bands in the range of the potential is equal to $\left[\frac{4}{3\pi}\frac{1}{\mathsf{h}} \right]$ or to $\left[\frac{4}{3\pi}\frac{1}{\mathsf{h}} \right]+1$, for all $\mathsf{h} < c_0^{-\frac23}$.
\vskip 3mm

\noindent Note that we do not prove in this paper a lower bound of $\delta_p$. A still open conjecture is wether or not $\delta_p$ has an exponential lower bound.

\noindent The fact that all the gaps are open is also a consequence of general results which state that the potentials for which there is a finite number of gaps are analytic functions (see \cite{Shu} for references on the topic and a discussion of the results of Skriganov).

\noindent The inequality $ I\left(\left( \tfrac76 \right)^{\frac23} \right) \frac{2^{\frac13}\pi^{\frac23}  }{9} \frac{1}{(p+1)^{\frac13}} <  \gamma_{p}$ in  (\ref{eq_thm_intro_gaps}) shows that $\gamma_p$ cannot be smaller than any power of $p$. It was expected since, by results of Hochstadt, it would imply that $V$ is a smooth function (see \cite{Ho}). Moreover, an exponentially small upper bound of $\gamma_p$ is characteristic from the analyticity of $V$ (see \cite{T}). For general results on singular potentials for the Hill equation, we refer to \cite{DM}. 

\vskip 3mm

\noindent The inequality (\ref{eq_thm_intro_bandes}) implies an upper bound for the spectral density in the range of the rescaled potential $\mathbf{V}$ in the semiclassical limit, which is also characterized by ``$\mathsf{c}$ tends to infinity''.  Let $k_0$ be the integer defined in (\ref{def_k_0}). For any $\mathsf{c} > c_0$ we denote by $D(\mathsf{c})$ the sum of the lengths of the $k_0+1$ first rescaled spectral bands (which are all included in the range of $\mathbf{V}$) divided by the length of the range of $\mathbf{V}$:
$$\forall \mathsf{c} >c_0,\ D(\mathsf{c}) = \frac{1}{\mathsf{c}} \sum_{p=0}^{k_0}  \delta_p.$$

\begin{cor}\label{cor_spectral_density}
When $\mathsf{c}$ tends to infinity,  $D(\mathsf{c})$ admits a limit denoted by $D_V$. Moreover,
\begin{equation}\label{eq_cor_spectral_density}
0 < D_V \leq \left( \frac{2}{3} \right)^{\frac13} . 
\end{equation}

\end{cor}
\vskip 2mm

\noindent The limit $D_V$ can be interpreted as the integrated spectral density in the range of the potential $V$ in the semiclassical limit. 

\vskip 3mm

\noindent Note that the number of gaps intersecting $[-V_0,0]$ increases by one each time the semiclassical parameter is equal to one of the numbers $c_p$, $p\geq 0$. To complete the first point of Theorem \ref{thm_intro_pbands} we observe that the roots of the canonical solutions of the Airy equation and their derivatives characterize the values of the counting parameter $\mathsf{c}$ for which a spectral band either enters in the range of the potential $[-V_0,0]$ or completes its entrance:

\begin{prop}\label{thm_zero_edges}
There exists a unique spectral band for which either the upper or the lower edge is equal to $0$ if and only if $\mathsf{c} \in \{ c_p, \tilde{c}_p \}_{p\geq 0}$.
\end{prop}

\noindent  For any function $f$, let $Z(f)$ denotes the set of the zeroes of $f$. Then, Proposition \ref{thm_zero_edges} implies that, for every $p\geq 0$,
$$Z\left( \mathsf{c} \mapsto \mathbf{E}_{\mathrm{max}}^{p}(\mathsf{c})  \right) = Z(v) \cup Z(v')$$
and
$$Z\left( \mathsf{c} \mapsto \mathbf{E}_{\mathrm{min}}^{p+1}(\mathsf{c})  \right) = Z(u) \cup Z(u').$$

\vskip 5mm

\subsection{Spectral bands and spectral gaps in the semiclassical regime}

\noindent Thanks to the explicit form of the bands, other estimates are proven. 
\vskip 3mm

\noindent \textbf{Notation.} Let $j\geq 0$ and define the real numbers $\mathfrak{a}_{2j}= \tilde{a}_{j+1}$ and  $\mathfrak{a}_{2j+1}=a_{j+1}$.

\vskip 3mm

\begin{theo}\label{thm_intro_width}
\begin{enumerate}
 \item Let $p\geq 0$. The rescaled and shifted $p$-th spectral band, 
$$[\mathsf{h}^{-\frac23} +\mathbf{E}_{\mathrm{min}}^{p},  \mathsf{h}^{-\frac23} +\mathbf{E}_{\mathrm{max}}^{p} ]$$
tends to the singleton $\{ \mathfrak{a}_p \}$ when $\mathsf{h}$ tends to $0$.
\item  Its width is, for every $\mathsf{h} \in (0, c_{2j}^{-\frac32})$,
\begin{equation}\label{intro_thm_width_2j}
 \delta_{2j} = 2\alpha \sqrt{3} \mfrac{(u'(-\mathfrak{a}_{2j}))^2  }{\mathfrak{a}_{2j} } \ex^{-\frac43 \mathsf{h}^{-1} + 2\mathfrak{a}_{2j}\mathsf{h}^{-\frac13}} \left( 1 +  \mathcal{O} \left(\mathsf{h}^{\frac13} \right) \right)
\end{equation}
when $p=2j$ and, for every $\mathsf{h} \in (0, c_{2j+1}^{-\frac32})$,
\begin{equation}\label{intro_thm_width_2j1} 
\delta_{2j+1}  =  2 \alpha \sqrt{3} (u(-\mathfrak{a}_{2j+1}))^2 \ex^{-\frac43 \mathsf{h}^{-1} + 2\mathfrak{a}_{2j+1}\mathsf{h}^{-\frac13}} \left( 1 +  \mathcal{O} \left(\mathsf{h}^{\frac13} \right) \right)
\end{equation}
when $p=2j+1$.
\item  Let $\mathsf{h}\in (0,c_0^{-\frac32})$ be arbitrary small. For every $p\in \{0,\ldots, \left[\frac{4}{3\pi}\frac{1}{\mathsf{h}} \right]  \}$, (\ref{intro_thm_width_2j}) holds true if $p=2j$ is even and (\ref{intro_thm_width_2j1}) holds true if $p=2j+1$ is odd.
\end{enumerate}


\end{theo}
\vskip 5mm

\noindent We can interpret the first statement of Theorem \ref{thm_intro_width} as a convergence of the band spectrum of the periodic Airy-Schr\"odinger operator to the pure point spectrum of the Schr\"odinger operator for a linear potential well, since it is a classical result that this spectrum is exactly $\{ \mathfrak{a}_p \}_{p\geq 0}$. Note that the eigenspace associated with $\mathfrak{a}_{2j}$ is spanned by $x\mapsto Ai(|x|-\mathfrak{a}_{2j})$ and the eigenspace associated with $\mathfrak{a}_{2j+1}$ is spanned by $x\mapsto \mathrm{sign}(x)\cdot Ai'(|x|-\mathfrak{a}_{2j+1})$. 
\vskip 2mm

\noindent Also note that the shifting of the rescaled spectral bands corresponds to the shifting of the spectrum of the operator $-\frac{\dd^2}{\dd y^2} + (q+1)$ as defined in (\ref{eq_H_classic_form}). 
\vskip 2mm

\noindent The second statement can be interpreted as the fact that the width of the rescaled $p$-th spectral band is expressed in the semiclassical regime and that an explicit interval of validity of the semiclassical approximation is given by $(0,c_p^{-\frac32})$.

\vskip 2mm

\noindent It shows that the width of the spectral bands (multiplied by $\mathsf{h}^{\frac23}$ to take in account the scale in $\mathsf{h}$ of the right member of (\ref{eq_H_Harrell_form})) is proportional to 
\begin{equation}
\ex^{-\frac43 \mathsf{h}^{-1}+2\mathfrak{a}_p \mathsf{h}^{-\frac13} } \left(\mathsf{h}^{\frac23}+ \mathcal{O} \left( \mathsf{h}\right) \right).
\end{equation}
By comparison, the result of Theorem 1.1 of \cite{H} shows thinner spectral bands of widths proportional to an exponential term equal to 
\begin{equation}
\mathsf{h} \ex^{-\frac43 \mathsf{h}^{-1}(1-\mathfrak{a}_p\mathsf{h})^{\frac32} } \left(1+ \mathcal{O} \left( \mathsf{h}^{\frac14} \right) \right) = \ex^{-\frac43 \mathsf{h}^{-1}+2\mathfrak{a}_p } \left(\mathsf{h}+ \mathcal{O} \left( \mathsf{h}^{\frac54} \right) \right).
\end{equation} 
This result is expected since the size of the spectral bands depends strongly on the regularity of the potential. In \cite{H}, the potential is assumed to be at least two times differentiable at its minima and maxima points while in our case, it is not even differentiable at these points. 

\noindent The widening of the spectral bands is similar to a Gevrey $3$ effect, on the size of the bands, of the singularity of the potential at its minima. This Gevrey $3$ effect is characterized by the $\ex^{2\mathfrak{a}_p \mathsf{h}^{-\frac13} }$ term. 
The Gevrey $3$ effect appears also in the diffraction of a wave by a strictly convex obstacle (for a mixed boundary condition \cite{OL2} or a Neumann boundary condition \cite{OL3}). It is also the main result of \cite{BLR} for the scattering frequencies of the wave operator near a convex analytic obstacle. 
\vskip 2mm

\noindent The asymptotic result of Harrell on the width of the spectral bands concerns  the eigenvalues $E$ close to a fixed odd number $2n+1$, which corresponds to eigenvalues of $-\frac{\hbar^2}{2m}\frac{\dd^2}{\dd y^2}+q(y)$ close to $\frac{2n+1}{\sqrt{m} }\hbar$ (since $q''(0)=\frac12$). For results on the spectrum of the perturbed harmonic oscillator, see \cite{BBL}. In our results, it corresponds to $\mathsf{h}$ tends to  $0$ and $\mathbf{E}+ \mathsf{h}^{-\frac23}$ close to an eigenvalue of $-\frac{\dd^2}{\dd x^2}+\vert x\vert $, hence to the values in $\{\mathfrak{a}_p\}_{p\geq 0}$. The result of Harrell was later precised in \cite{HS1} and generalized for the first spectral band in the multidimensional case (see \cite{Ou,S}).
\vskip 2mm

\noindent The last statement of Theorem \ref{thm_intro_width} is a quantitative result in the \textbf{semiclassical limit}. It means that for every $p\geq 0$, there exists a real number $\mathsf{h}_p>0$ such that the estimate (\ref{intro_thm_width_2j}) or (\ref{intro_thm_width_2j1}) holds true for every $\mathsf{h}\in (0,\mathsf{h}_p)$ and the sequence $(\mathsf{h}_p)_{p\geq 0}$ tends to $0$ when $p$ tends to infinity. In particular, when the semiclassical parameter $\mathsf{h}$ tends to $0$, $\left[\frac{4}{3\pi}\frac{1}{\mathsf{h}} \right]$ tends to infinity and the estimates (\ref{intro_thm_width_2j}) and (\ref{intro_thm_width_2j1}) are valid for all the rescaled spectral bands in the range of $\mathbf{V}$.  

\vskip  5mm

\noindent The rescaled spectral gaps have constant widths, compared with $\mathsf{h}$, which are the differences of two consecutives eigenvalues of the Schr\"odinger operator for a linear potential well.

\begin{theo}\label{thm_intro_gap}
\begin{enumerate}
 \item Let $\tau >0$ be arbitrary small. For every $p\geq 0$ and every $h\in (0,(\mathfrak{a}_{p+1}+\tau)^{-\frac32}]$,
\begin{equation}\label{eq_thm_intro_gap_p}
 \gamma_{p}  = \mathfrak{a}_{p+1}-\mathfrak{a}_{p} + \mathcal{O}\left( \ex^{-\frac43 \mathsf{h}^{-1} + 2\mathfrak{a}_{p+1}\mathsf{h}^{-\frac13}}   \right).
\end{equation}
\item Let $h\in (0,c_1^{-\frac32})$ be arbitrary small. For every $p\in \{0,\ldots,  \left[\frac{4}{3\pi}\frac{1}{\mathsf{h}} \right] -1 \}$, (\ref{eq_thm_intro_gap_p}) holds true.
\end{enumerate}
\end{theo}

\noindent Thus, the width of the rescaled $p$-th spectral gap is in the semiclassical regime and an explicit interval of validity of the semiclassical approximation is given by $(0,\mathfrak{a}_{p+1}^{-\frac32})$.
\vskip 2mm

\noindent More precise estimates for the widths of the gaps are found in Proposition \ref{prop_asymp_gap}. 
\vskip 2mm

\noindent Theorems \ref{thm_intro_width} and \ref{thm_intro_gap} are stated in the case where $\mathsf{c}+\mathbf{E}$ remains bounded and close to one of the $\mathfrak{a}_p$ for $p$ fixed, and $\mathsf{c}>\mathfrak{a}_{p+1}$. It corresponds to looking at spectrals bands and gaps close to the bottom of the potential and thus to the bottom of the spectrum in the semiclassical limit. 

\noindent We can prove similar estimates in the case where $\mathbf{E}$ remains bounded and close to one of the $-\mathfrak{a}_p$  for $p$ fixed, and $\mathsf{c}+\mathbf{E}$ tends to $+\infty$. This second case corresponds to energy bands close to the maximum of the potential. In this second case, we number the bands starting with the band closest to $0$ in $[-\mathsf{c},0]$ labeled as the band $1$, while the lowest band is labeled by $k_0+1$ defined in (\ref{def_k_0}).

\noindent Another case which can be dealt with, is when one assume that $\mathbf{E}$ is close to $-\frac{\mathfrak{a}_p}{\mathsf{c}^{\frac12}}$ for some fixed $p$ and there exists constants $C_1>0$ and $C_2>0$ such that $\frac{\mathfrak{a}_p}{\mathsf{c}^{\frac12}}\in [C_1,C_2]$. In this case our bootstrap technique proof of Theorem \ref{thm_bottom} applies since $\ex^{-\frac43(\mathsf{c} - \mathfrak{a}_p)^{\frac32} } \cdot \ex^{\frac43\mathsf{c}^{\frac32} } \in [\ex^{C_1}, \ex^{C_2}]$. This proves a result similar to Theorem $8.1$ of \cite{marz} (it is the case $\mu \in [-Ch,0]$ in its notations) for the widths of the gaps and the bands, but with a different order of magnitude of these widths. These differences are due to the fact that in \cite[Theorem 8.1]{marz}, the potential is supposed to be analytic and for us it is only piecewise analytic.
\vskip 5mm

\subsection{Outline of the paper}

\noindent In Section \ref{sec_band_spec}, we recall the classical theory of periodic operators and their band spectra. It provides the equations to solve to determine the edges of the spectral bands. In Section \ref{sec_asymtptotics}, we study the asymptotics of the canonical solutions, their derivatives and the asymptotics of their zeroes. These very precise asymptotics imply an important result of separation and ordering of the zeroes of the canonical solutions and their derivatives (Corollary \ref{prop_sep_ck}). This result is the key result which allows to distinguish the upper edge of a spectral band from the lower edge of the next one among the solutions of the equations obtained in Section \ref{sec_band_spec}. This identification of the spectral edges is performed in Sections \ref{sec_band1} and \ref{sec_pbands}.

\noindent Section \ref{sec_prelim} is devoted to the study of families of strictly monotonous and continuous functions which allows to prove in Sections \ref{sec_band1} and \ref{sec_pbands} the existence of solutions to the equations which define the spectral edges. Section \ref{sec_prelim} describes also the graphical interpretation of these equations, in terms of the functions $\frac{v}{u}$ and $\frac{v'}{u'}$, which has guided our analysis throughout this paper. 

\noindent Section \ref{sec_bottom} is devoted to the proof of Theorem \ref{thm_bottom}. This proof contains most of the ideas used later in Section \ref{sec_asymtptotics} to get all the estimates of the width of the spectral bands and of the spectral gaps in the semiclassical regime. We also investigate in Section \ref{sec_E0max} the behaviour of the upper edge of the first spectral band in both semiclassical limit and when $\mathsf{h}$ tends to infinity. In the limit when $\mathsf{h}$ tends to infinity, the integer $p_0$ introduced in (\ref{def_p_0}) plays a crucial role. In Section \ref{sec_pbands} we characterize the spectral edges of all the spectral bands in the range of $V$, we count these bands and we prove Theorem \ref{thm_intro_pbands}. We also prove the result on the integrated spectral density in the range of $V$. 

\noindent In Appendix \ref{sec_zk}, the monotonicity of the functions $z_k$ introduced in Lemma \ref{lem_zeroes} is proven. This monotonicity result is particularly technical and requires a version of the Sturm Picone's lemma about interlacing of zeroes of solutions of ordinary differential equations adapted to our setting (\cite{CL}). Such result is proven in Appendix \ref{sec_sturm_picone}.

\section{The band structure of the spectrum of $H$}\label{sec_band_spec}

\noindent In this section we recall the equations characterizing the spectral edges of the spectrum of the operator $H$, using the general theory of periodic Schr\"odinger operators  (\cite{RS4}). Let $\omega \in [-L_0,L_0]$. We start by considering the restriction $H(\omega)$ of $H$ to $\mathsf{H}^2([-L_0,L_0])$, the Sobolev space of functions $\psi\in \mathsf{H}^2(\R)$ which satisfy 
\begin{equation}\label{eq_boundary_cond_0}
\forall x \in \R,\ \psi(x+2L_0)= \ex^{\ii (\frac{\pi}{L_0}\omega + \pi)} \psi(x).
\end{equation}
Note that, as $\mathsf{H}^2([-L_0,L_0])\subset C^1([-L_0,L_0])$, this condition is equivalent to the boundary conditions:  
\begin{equation}
\label{eq_boundary_cond} \psi(L_0)=\ex^{\ii (\frac{\pi}{L_0}\omega + \pi)} \psi(-L_0)\quad \mbox{and} \quad \psi'(L_0)=\ex^{\ii (\frac{\pi}{L_0}\omega + \pi)} \psi'(-L_0).
 \end{equation}
The operator $H(\omega)$ is self-adjoint. It is Hilbert-Schmidt and thus compact. Its spectrum is pure point and the eigenvalues of $H(\omega)$ are solutions of explicit equations. According to \cite{RS4}, $H$ is the direct integral of the operators $H(\omega)$: 
$$H=\int_{[-L_0,L_0]}^{\oplus} H(\omega) \dd \omega.$$
This decomposition in direct integral allows to recover the spectrum of $H$ from the spectra of the $H(\omega)$'s.
\vskip 2mm

\noindent From the canonical solutions $u$ and $v$ of the Airy equation, one defines the canonical pair of odd and even solutions of (\ref{eq_schrodinger}) on the interval $[-L_0,L_0]$.  These functions, denoted by $\mathsf{U}$ and $\mathsf{V}$, are defined, for every $x\in [-L_0,L_0]$, by 
\[
\mathsf{U}(x)=-v'(-\mathsf{c} - \mathbf{E})u(\mathbf{V}(x)-\mathbf{E}) + u'(-\mathsf{c} - \mathbf{E})v(\mathbf{V}(x)-\mathbf{E})
 \]
and
\[
 \mathsf{V}(x)=\mathrm{sign}(x)\left( -v(-\mathsf{c} - \mathbf{E})u(\mathbf{V}(x)-\mathbf{E}) + u(-\mathsf{c} - \mathbf{E})v(\mathbf{V}(x)-\mathbf{E})\right).
 \]
They form a basis of even and odd $C^1$ solutions of the equation (\ref{eq_schrodinger}) on the interval $[-L_0,L_0]$. Their Wronskian satisfies 
$$\forall x\in [L_0,L_0],\ (\mathsf{U}\mathsf{V}^{'} - \mathsf{U}^{'}\mathsf{V})(x)= \frac{\mathsf{c}}{L_0}.$$

\noindent Since any solution of (\ref{eq_schrodinger}) is a linear combination of $\mathsf{U}$ and $\mathsf{V}$, the boundary conditions (\ref{eq_boundary_cond}) rewrite:
\begin{align}
 A \mathsf{U}(L_0) + B \mathsf{V}(L_0)  &= -\ex^{\ii (\frac{\pi}{L_0}\omega)}\left( A \mathsf{U}(L_0) -B \mathsf{V}(L_0) \right)\\
A \mathsf{U}^{'}(L_0) + B \mathsf{V}^{'}(L_0)  &= -\ex^{\ii (\frac{\pi}{L_0}\omega)}\left( -A \mathsf{U}^{'}(L_0) +  B \mathsf{V}^{'}(L_0) \right), 
\end{align}
for $A,B \in \R$. Thus, $A$ and $B$ are solutions of the linear system: 
\begin{equation}\label{eq_sys_AB}
\left\lbrace \begin{array}{lcl}
             \mathsf{U}(L_0) \left(1+\ex^{\ii (\frac{\pi}{L_0}\omega)} \right) A + \mathsf{V}(L_0) \left(1-\ex^{\ii (\frac{\pi}{L_0}\omega)} \right)  B & = & 0 \\
	      \mathsf{U}^{'}(L_0) \left(1-\ex^{\ii (\frac{\pi}{L_0}\omega)} \right) A + \mathsf{V}^{'}(L_0) \left(1+\ex^{\ii (\frac{\pi}{L_0}\omega)} \right)  B & = & 0
             \end{array}\right. . \end{equation}

\noindent Considering the system (\ref{eq_sys_AB}), one gets that $E\in \R$ is an eigenvalue of $H(\omega)$ if and only if
\begin{equation} \label{eq_sys_vp_Ha}\left|\begin{matrix}
         \mathsf{U}(L_0) (1+\ex^{\ii (\frac{\pi}{L_0}\omega)} ) \ &\   \mathsf{V}(L_0) (1-\ex^{\ii (\frac{\pi}{L_0}\omega)} ) \\
	 \mathsf{U}^{'}(L_0) (1-\ex^{\ii (\frac{\pi}{L_0}\omega)} )\  &\   \mathsf{V}^{'}(L_0) (1+\ex^{\ii (\frac{\pi}{L_0}\omega)} )
        \end{matrix} \right|=0.
        \end{equation}
The determinant in (\ref{eq_sys_vp_Ha}) being analytic in $E$, equation (\ref{eq_sys_vp_Ha}) has only a discrete set of solutions in $E$ (since this determinant is not equal to $0$ for every $E$) which is consistent with the fact that $H(\omega)$ is a compact operator. For $\omega=-L_0$, we get that (\ref{eq_sys_vp_Ha}) is equivalent to 
\begin{equation}
\label{eq_vp_Ha_d} -4 \mathsf{U}^{'}(L_0)\cdot \mathsf{V}(L_0) =0. 
\end{equation}
For $\omega=0$, we get that (\ref{eq_sys_vp_Ha}) is equivalent to 
\begin{equation}\label{eq_vp_Ha_0} 
 4 \mathsf{U}(L_0)\cdot  \mathsf{V}^{'}(L_0)  =0. 
 \end{equation}
Since the expression in the left member of (\ref{eq_vp_Ha_d}) is analytic in $E$ we denote by $E_{\mathrm{min}}^0\leq E_{\mathrm{max}}^1 \leq E_{\mathrm{min}}^2 \leq E_{\mathrm{max}}^3 \leq \ldots$ the elements of the spectrum of $H(-L_0)$,
$$\sigma(H(-L_0)):=  \{ E_{\mathrm{min}}^0, E_{\mathrm{max}}^1, E_{\mathrm{min}}^2, E_{\mathrm{max}}^3,\ldots \} = \left\{E\in \R,\ \mathsf{U}^{'}(L_0) \mathsf{V}(L_0)=0 \right\} $$

\noindent and we denote by $E_{\mathrm{max}}^0\leq E_{\mathrm{min}}^1 \leq E_{\mathrm{max}}^2 \leq E_{\mathrm{min}}^3\leq \ldots$   the elements of the spectrum of $H(0)$
$$\sigma(H(0)):= \{ E_{\mathrm{max}}^0, E_{\mathrm{min}}^1, E_{\mathrm{max}}^2, E_{\mathrm{min}}^3,\ldots \}  = \left\{E\in \R,\ \mathsf{U}(L_0)  \mathsf{V}^{'}(L_0)=0 \right\}.$$

\noindent Then, using \cite[Theorem XIII.90]{RS4}, the spectrum of $H$ is the band spectrum:
$$\sigma(H)= \bigcup_{p\geq 0} \left[ E_{\mathrm{min}}^p, E_{\mathrm{max}}^p \right].$$
Moreover $\sigma(H)$ is purely absolutely continuous and $H$ has no eigenvalues. To compute the spectrum of $H$, it remains to determine the edges of the spectral bands $E_{\mathrm{min}}^p$ and $E_{\mathrm{max}}^p$ for $p\geq 0$.
\vskip 3mm

\noindent \textbf{Remark.} Recall that the result here is true for any periodic potential (at least locally integrable) and thus our computations and the description of the spectrum as a band spectrum are valid for any symmetric potential.

\section{Asymptotics of the canonical solutions and of their zeroes}\label{sec_asymtptotics}

\noindent The aim of this Section is to obtain precise estimates of the roots of $u$, $v$, $u'$ and $v'$. This will be a consequence of the asymptotic expansions of the canonical solutions of the Airy equation on $\R_{-}$.

\subsection{Asymptotics of the canonical solutions}\label{sec_asymtptotics_canonical}

\noindent In order to obtain asymptotic expansions of the canonical solutions of the Airy equation, we start with asymptotic expansions of the Airy functions $Ai$, $Bi$ and their derivatives which are deduced from Bessel functions for which asymptotic expansions are well known.

\noindent For this purpose, one defines functions $P(\nu,\cdot)$ and $Q(\nu,\cdot)$ for any real number $\nu$ through the Bessel functions $J_{\nu}$ and $Y_{\nu}$ (see \cite{O})
$$P(\nu,\xi) =  \sqrt{\mfrac{\pi \xi}{2}}  \left( J_{\nu}(\xi)\cos\left(\xi-\tfrac12 \nu \pi - \tfrac14 \pi \right) + Y_{\nu}(\xi) \sin\left(\xi-\tfrac12 \nu \pi - \tfrac14 \pi\right)\right) $$
and 
$$ Q(\nu,\xi) =  \sqrt{\mfrac{\pi \xi}{2}}  \left(  Y_{\nu}(\xi) \cos\left(\xi-\tfrac12 \nu \pi - \tfrac14 \pi\right)-J_{\nu}(\xi)\sin\left(\xi-\tfrac12 \nu \pi - \tfrac14 \pi\right)\right).$$
The functions  $P(\nu,\cdot)$ and $Q(\nu,\cdot)$ have known expansions which are used to get asymptotic expansions of the canonical solutions and their derivatives. Note that, when $\xi$ tends to $+\infty$, $P(\nu,\xi) \sim 1$ and $Q(\nu,\xi)\sim \frac{4\nu^2 -1}{8\xi}$.

\begin{prop}\label{prop_asym_uv}
For every $x>0$, we set $\xi=\frac23 x^{\frac32}$. We have:
\begin{align}
u(-x) & =     2 \pi^{\frac12} x^{-\frac14} Ai'(0) \left(  \sin\left(\xi - \tfrac{7\pi}{12}\right)P\left(\tfrac{1}{3},\xi\right) + \cos\left(\xi - \tfrac{7\pi}{12}\right)Q\left(\tfrac{1}{3},\xi\right)  \right), \label{eq_u_asym} \\
u'(-x) & =  - 2 \pi^{\frac12} x^{\frac14} Ai'(0) \left(  \cos\left(\xi - \tfrac{7\pi}{12}\right)P\left(\tfrac{2}{3},\xi\right) - \sin\left(\xi - \tfrac{7\pi}{12}\right)Q\left(\tfrac{2}{3},\xi\right)  \right), \label{eq_up_asym} \\
v(-x) & =  - 2 \pi^{\frac12} x^{-\frac14} Ai(0) \left(  \sin\left(\xi + \tfrac{\pi}{12}\right)P\left(\tfrac{1}{3},\xi\right) + \cos\left(\xi + \tfrac{\pi}{12}\right)Q\left(\tfrac{1}{3},\xi\right)  \right), \label{eq_v_asym} \\
v'(-x) & =  -2 \pi^{\frac12} x^{\frac14} Ai(0) \left( - \cos\left(\xi + \tfrac{\pi}{12}\right)P\left(\tfrac{2}{3},\xi\right) + \sin\left(\xi + \tfrac{\pi}{12}\right)Q\left(\tfrac{2}{3},\xi\right)  \right).\label{eq_vp_asym}
\end{align}
\end{prop}

\vskip 3mm

\noindent The proof of Proposition \ref{prop_asym_uv} is given in Appendix \ref{sec_proofs}.
\vskip 3mm

\noindent \textbf{Remark.} Note that $x\mapsto \pi^{-\frac12}x^{-\frac14} \cos(\xi-\frac{\pi}{4})$ and  $x\mapsto \pi^{-\frac12}x^{-\frac14} \sin(\xi-\frac{\pi}{4})$ are solutions of the equation $y''=\left(-x + \frac{5}{16x^2}\right)y$ and thus are approximate solutions of $y''=-xy$. Hence the existence of $P$ and $Q$ can be seen as an application of the Duhamel principle.
\vskip 3mm

\subsection{Asymptotics and ordering of the zeroes of the canonical solutions}\label{sec_ordering_zeroes}  

Before having precise intervals in which we localize the zeroes of $u$, $u'$, $v$ and $v'$, we localize them between zeroes of the classical Airy function $Ai$ and its derivative. Since the zeroes of $Ai$ and $Ai'$ are known, it will guide us to choose the good intervals in which we will verify that $u$, $u'$, $v$ and $v'$ do vanish and outside of which they do not. 
\vskip 2mm

\noindent We start by looking at the variations of the functions $\frac{Bi}{Ai}$ and $\frac{Bi'}{Ai'}$. The functions $\frac{Bi}{Ai}$ and $\frac{Bi'}{Ai'}$ have the following behaviours:
\begin{enumerate}
\item  On every interval $(-a_{j+1},-a_j)$, the function $\frac{Bi}{Ai}$ is continuous, increasing and is a bijection from $(-a_{j+1},-a_j)$ to $\R$. Moreover,  $\frac{Bi}{Ai}$ is continuous, increasing and is a bijection from $(-a_1,+\infty)$ to $\R$.
\item  On every interval $(-\tilde{a}_{j+1},-\tilde{a}_j)$, the function $\frac{Bi'}{Ai'}$ is continuous, increasing and is a bijection from $(-\tilde{a}_{j+1},-\tilde{a}_j)$ to $\R$. Moreover, $\frac{Bi'}{Ai'}$ is continuous and increasing on $(-\tilde{a}_1,0]$ from $-\infty$ to $0$. It is also continuous, decreasing and a bijection from $(-\infty,0]$ to $[0,+\infty)$.
\end{enumerate}

\vskip 4mm
\vskip 5mm

\noindent Note that the roots of $u$, $v$, $u'$ and $v'$ are exactly the solutions of the  equations $\frac{Bi}{Ai}(x)=\frac{Bi'}{Ai'}(0)$, $\frac{Bi}{Ai}(x)=\frac{Bi}{Ai}(0)$,    $\frac{Bi'}{Ai'}(x)=\frac{Bi'}{Ai'}(0)$ and $\frac{Bi'}{Ai'}(x)=\frac{Bi}{Ai}(0)$, respectively. We identify the solutions of these equations and the sequences of zeroes introduced in Section \ref{sec_model_can_sol}.

\vskip 2mm

\noindent The equation 
$$\frac{Bi}{Ai}(x)=\frac{Bi'}{Ai'}(0)=-\sqrt{3}$$
has a countable number of solutions which are negative and no positive solution. Let us denote by $( -\tilde{c}_{2j} )_{j \geq 0}$ the sequence of all the solutions arranged in decreasing order with, for every $j\geq 0$, $-\tilde{c}_{2j} \in (-a_{j+1},-\tilde{a}_{j+1})$. Since the sets of the zeroes of $Ai$ and of $u$ are disjoints, \textbf{the set of the zeroes of $u$ is exactly $\{ -\tilde{c}_{2j} \}_{j \geq 0}$}.
\vskip 2mm

\noindent Similarly, the equation 
$$\frac{Bi'}{Ai'}(x)=\frac{Bi'}{Ai'}(0)=-\sqrt{3}$$
has a countable number of solutions which are negative and no positive solution except $0$. Let us denote by $(-\tilde{c}_{2j+1} )_{j \geq 0}$ the sequence of all the solutions arranged in decreasing order with, for every $j\geq 0$, $-\tilde{c}_{2j+1} \in (-\tilde{a}_{j+2},-a_{j+1})$. Since the sets of the zeroes of $Ai'$ and of $u'$ are disjoints, \textbf{the set of the zeroes of $u'$ is exactly $\{ -\tilde{c}_{2j+1} \}_{j \geq 0}\cup \{ 0\}$}.
\vskip 2mm

\noindent The equation
$$\frac{Bi}{Ai}(x)=\frac{Bi}{Ai}(0)=\sqrt{3}$$
has a countable number of solutions which are negative and no positive solution except $0$. Let us denote by $( -c_{2j+1} )_{j \geq 0}$ the sequence of all the solutions arranged in decreasing order with, for every $j\geq 0$, $-c_{2j+1} \in (-\tilde{a}_{j+2},-a_{j+1})$. Since the sets of the zeroes of $Ai$ and of $v$ are disjoints, \textbf{the set of the zeroes of $v$ is exactly $\{ -c_{2j+1} \}_{j \geq 0} \cup \{ 0\} $} . 
\vskip 2mm

\noindent Finally, the equation
$$\frac{Bi'}{Ai'}(x)=\frac{Bi}{Ai}(0)=\sqrt{3}$$
also has a countable number of solutions which are negative and no positive solution. Let us denote by $( -c_{2j} )_{j \geq 0}$ the sequence of all the solutions arranged in decreasing order with, for every $j\geq 0$, $-c_{2j} \in (-a_{j+1},-\tilde{a}_{j+1})$. Since the sets of the zeroes of $Ai'$ and of $v'$ are disjoints, \textbf{the set of the zeroes of $v'$ is exactly $\{ -c_{2j} \}_{j \geq 0}$}. 
\vskip 3mm

\noindent From the asymptotic expansions (\ref{eq_u_asym}), (\ref{eq_up_asym}), (\ref{eq_v_asym}) and (\ref{eq_vp_asym}) and the distribution of the sequences of zeroes of $Ai$ and $Ai'$, we get intervals in which we localize the constants $-c_p$ and $-\tilde{c}_p$. We also obtain the variations of $u$ and $v$.

\begin{prop}\label{prop_loc_ck}
\begin{enumerate}
 \item For every $j\geq 0$, the function $v'$ has a unique zero in the interval $(-\left(\frac32 (j\pi + \frac{\pi}{2}) \right)^{\frac23},-\left(\frac32 (j\pi + \frac{\pi}{3}) \right)^{\frac23} )$ and does not vanish outside of these intervals. Thus, 
\begin{equation}\label{eq_interval_c2j} 
-c_{2j} \in \left(-\left(\tfrac32 (j\pi + \tfrac{\pi}{2}) \right)^{\frac23},-\left(\tfrac32 (j\pi + \tfrac{\pi}{3}) \right)^{\frac23} \right).
\end{equation}

\item For every $j\geq 0$, the function $u$ has a unique zero in the interval

\noindent $(-\left(\frac32 (j\pi + \frac{2\pi}{3}) \right)^{\frac23},-\left(\frac32 (j\pi + \frac{\pi}{2}) \right)^{\frac23} )$ and does not vanish outside of these intervals. Thus, 
\begin{equation}\label{eq_interval_tc2j} 
-\tilde{c}_{2j} \in \left(-\left(\tfrac32 (j\pi + \tfrac{2\pi}{3}) \right)^{\frac23},-\left(\tfrac32 (j\pi + \tfrac{\pi}{2}) \right)^{\frac23} \right).
\end{equation}
\item For every $j\geq 0$, the function $v$ has a unique zero in the interval

\noindent  $(-\left(\frac32 (j\pi + \pi) \right)^{\frac23},-\left(\frac32 (j\pi + \frac{5\pi}{6}) \right)^{\frac23} )$ and does not vanish outside of these intervals. Thus, 
\begin{equation}\label{eq_interval_c2j1} 
-c_{2j+1} \in \left(-\left(\tfrac32 (j\pi + \pi) \right)^{\frac23},-\left(\tfrac32 (j\pi + \tfrac{5\pi}{6}) \right)^{\frac23} \right).
\end{equation}

\item For every $j\geq 0$, the function $u'$ has a unique zero in the interval

\noindent  $(-\left(\frac32 (j\pi + \frac{7\pi}{6}) \right)^{\frac23},-\left(\frac32 (j\pi + \pi) \right)^{\frac23} )$ and does not vanish outside of these intervals. Thus, 
\begin{equation}\label{eq_interval_tc2j1} 
-\tilde{c}_{2j+1} \in \left(-\left(\tfrac32 (j\pi + \tfrac{7\pi}{6}) \right)^{\frac23},-\left(\tfrac32 (j\pi + \pi) \right)^{\frac23} \right).
\end{equation}
\end{enumerate}
\end{prop}

\begin{prop}\label{prop_variations_uv}
The variations of $u$ and $v$ and their signs between two consecutive zeroes are:
\begin{enumerate}
\item $u$ is positive on $(-\tilde{c}_0,+\infty)$ and for every $j\geq 0$, $(-1)^j u$ is negative on $[-\tilde{c}_{2j+2}, -\tilde{c}_{2j} ]$. It is strictly increasing on $(-\tilde{c}_1,+\infty)$,  and for every $j\geq 0$, $(-1)^j u$ is strictly decreasing on $[-\tilde{c}_{2j+3}, -\tilde{c}_{2j+1} ]$.
\item $v$ is positive on $[0,+\infty)$, negative on $[-c_1,0]$ and for every $j\geq 0$, $(-1)^j v$ is positive on $[-c_{2j+3}, -c_{2j+1} ]$. It is strictly increasing on $(-c_0,+\infty)$,  and for every $j\geq 0$, $(-1)^j v$ is strictly decreasing on $[-c_{2j+2}, -c_{2j} ]$.
\end{enumerate}
\end{prop}

\noindent The respective behaviour of $u$ and $v$ on respectively the intervals $(-\tilde{c}_0,+\infty)$ and $(-c_1,+\infty)$ are different than their respective behaviour on respectively the intervals $(-\infty,-\tilde{c}_0]$ and $(-\infty,-c_1]$.

\vskip 4mm

\noindent The proofs of Propositions \ref{prop_loc_ck} and \ref{prop_variations_uv} are given in Appendix \ref{sec_proofs}. The signs and variations of $u$, $u'$, $v$ and $v'$ are summarized in Appendix \ref{sec_app_variations}. From Proposition \ref{prop_loc_ck} we deduce immediately the ordering of the zeroes of the canonical solutions and their derivatives.

\begin{cor}\label{prop_sep_ck}
For every $p\geq 0$, $-\tilde{c}_p < -c_p$.
\end{cor}

\noindent From Proposition \ref{prop_loc_ck}, we also deduce asymptotics of the sequences $(c_p)_{p\geq 0}$ and  $(\tilde{c}_p)_{p\geq 0}$.

\begin{cor}\label{cor_asymp_ck}
One has
\begin{equation}\label{eq_approx_ck}
c_p =  \left(\frac{3p\pi}{4}\right)^{\frac23} + \mathcal{O}\left( \frac{1}{p^{\frac13} }\right) \ \mbox{ and }\ \tilde{c}_p-c_p = \left( \frac{\pi}{9\sqrt{2}} \right)^{\frac23} \frac{1}{p^{\frac13} } + \mathcal{O}\left( \frac{1}{p^{\frac43} } \right).
\end{equation}
\end{cor}

\begin{demo}
The first estimate in (\ref{eq_approx_ck}) follows directly from the asymptotics proven in Lemma \ref{lem_asymp_ck}. For the difference between $c_p$ and $\tilde{c}_p$, one uses (\ref{eq_approx_c2j})-(\ref{eq_approx_tc2j}) and (\ref{eq_approx_c2j1})-(\ref{eq_approx_tc2j1}). Indeed, for every $j\geq 0$,
$$\tilde{c}_{2j} -  c_{2j} =  \left(\frac{3j\pi}{2}\right)^{\frac23} \cdot\left(  \frac{7}{18j}- \frac{5}{18j} + \mathcal{O}\left( \frac{1}{j^{2} }\right) \right)  =   \left( \frac{\pi}{18} \right)^{\frac23} \frac{1}{j^{\frac13} } + \mathcal{O}\left( \frac{1}{j^{\frac43} } \right)$$
and similarly,
$$\tilde{c}_{2j+1} -  c_{2j+1}  =\left( \frac{\pi}{18} \right)^{\frac23} \frac{1}{j^{\frac13} } + \mathcal{O}\left( \frac{1}{j^{\frac43} } \right), $$
which proves (\ref{eq_approx_ck}). 
\end{demo}

\noindent We introduce, for every $p\geq 0$,
$$\xi_p = \frac{2}{3} c_p^{\frac{3}{2}}\quad \mbox{ and } \quad \tilde{\xi}_p = \frac{2}{3} \tilde{c}_p^{\frac{3}{2}}.$$

\begin{lem}\label{lem_asymp_ck}
Let $j\geq 0$. One has 
 \begin{align}
\xi_{2j}  \in & \left[ \frac{5\pi}{12}+j\pi - \frac{7}{12(j\pi+\frac{\pi}{3})},\frac{5\pi}{12}+j\pi + \frac{7}{12(j\pi+\frac{\pi}{3})} \right] \label{eq_approx_xi2j} \\
  \mbox{and}\quad  & c_{2j}  =  \left(\frac{3j\pi}{2}\right)^{\frac23} \cdot\left( 1  + \frac{5}{18j} + \mathcal{O}\left( \frac{1}{j^{2} }\right) \right), \label{eq_approx_c2j} 
\end{align}

\begin{align}
\xi_{2j+1}  \in & \left[ \frac{11\pi}{12}+j\pi - \frac{5}{36(j\pi+\frac{5\pi}{6})},\frac{11\pi}{12}+j\pi + \frac{5}{36(j\pi+\frac{5\pi}{6})} \right] \label{eq_approx_xi2j1} \\
  \mbox{and} \quad & c_{2j+1} =  \left( \frac{3j\pi}{2} \right)^{\frac23}\cdot \left(1 + \frac{11}{18j} + \mathcal{O}\left( \frac{1}{j^{2}} \right) \right),  \label{eq_approx_c2j1} 
\end{align} 

 \begin{align}
 \tilde{\xi}_{2j}  \in & \left[ \frac{7\pi}{12}+j\pi - \frac{5}{36(j\pi+\frac{\pi}{2})},\frac{7\pi}{12}+j\pi + \frac{5}{36(j\pi+\frac{\pi}{2})} \right] \label{eq_approx_txi2j} \\  
 \mbox{and}\quad  &     \tilde{c}_{2j} =  \left(\frac{3j\pi}{2} \right)^{\frac23}\cdot \left(1 + \frac{7}{18j}+ \mathcal{O}\left( \frac{1}{j^{2}} \right)\right), \label{eq_approx_tc2j}
\end{align}

\begin{align}
 \tilde{\xi}_{2j+1}  \in & \left[ \frac{13\pi}{12}+j\pi - \frac{7}{12(j+1)\pi},\frac{13\pi}{12}+j\pi + \frac{7}{12(j+1)\pi}  \right] \label{eq_approx_txi2j1} \\    
 \mbox{and}\quad &    \tilde{c}_{2j+1}  =  \left( \frac{3j\pi}{2}  \right)^{\frac23}\cdot \left(1 +  \frac{13}{18j} +\mathcal{O}\left( \frac{1}{j^{2}} \right) \right). \label{eq_approx_tc2j1}
 \end{align}
\end{lem}

\vskip 2mm

\noindent The proof of Lemma \ref{lem_asymp_ck} is given in Appendix \ref{sec_proofs}.
\vskip 3mm

\noindent \textbf{Remark.} From the proof of the asymptotic expansions of $c_p$ and $\tilde{c}_p$, one could obtain asymptotic expansion of these sequences at any order, using the developments of the functions $P$ and $Q$ (\cite[9.2.9 and 9.2.10]{AS}). One would then get similar formula as those for the zeroes of the functions $Ai$, $Ai'$, $Bi$ and $Bi'$ (\cite[10.4.94 and below]{AS}).

\section{Preliminaries to the computation of the band edges}\label{sec_prelim}

\subsection{Characterization of the spectral band edges}

\noindent The band edges are characterized by the functions $\mathsf{U}$, $\mathsf{V}$ and their derivatives,  through the equations (\ref{eq_vp_Ha_d}) and (\ref{eq_vp_Ha_0}).

\vskip 3mm

\noindent To find $\mathbf{E}_{\mathrm{min}}^p$ and $\mathbf{E}_{\mathrm{max}}^p$ and thus the band edges $E_{\mathrm{min}}^p$ and $E_{\mathrm{max}}^p$ for any $p\geq 0$, we have to solve the four equations:
\begin{alignat}{2}
\mathsf{U}^{'}(L_0) & =  u'(-\mathsf{c}-\mathbf{E})v'(-\mathbf{E})-v'(-\mathsf{c}-\mathbf{E})u'(-\mathbf{E}) && =0, \label{eq_band_lin1} \\
\mathsf{V}(L_0)& =  u(-\mathbf{E})v(-\mathsf{c}-\mathbf{E})-v(-\mathbf{E})u(-\mathsf{c}-\mathbf{E})&& =0, \label{eq_band_lin2} \\
\mathsf{U}(L_0) & =  v(-\mathbf{E})u'(-\mathsf{c}-\mathbf{E})-u(-\mathbf{E})v'(-\mathsf{c}-\mathbf{E})&& =0, \label{eq_band_lin3} \\
\mathsf{V}^{'}(L_0)& =  v'(-\mathbf{E})u(-\mathsf{c}-\mathbf{E})-u'(-\mathbf{E})v(-\mathsf{c}-\mathbf{E})&& =0. \label{eq_band_lin4} 
\end{alignat}

\vskip 5mm

\noindent We have the four equivalences:
\begin{enumerate}
 \item for $ \mathbf{E}\notin \{\tilde{c}_{2j+1} -\mathsf{c} \}_{j\geq 0} \cup \{  \tilde{c}_{2j+1} \}_{j\geq 0}$,
$$u'(-\mathsf{c}-\mathbf{E})v'(-\mathbf{E})-v'(-\mathsf{c}-\mathbf{E})u'(-\mathbf{E}) =0 \  \Leftrightarrow \  \mfrac{v'(-\mathsf{c}-\mathbf{E})}{u'(-\mathsf{c}-\mathbf{E})}=\mfrac{v'(-\mathbf{E})}{u'(-\mathbf{E})},$$ 
\item for $\mathbf{E} \notin \{ \tilde{c}_{2j} -\mathsf{c} \}_{j\geq 0}\cup  \{  \tilde{c}_{2j} \}_{j\geq 0}$,
$$ u(-\mathbf{E})v(-\mathsf{c}-\mathbf{E})-v(-\mathbf{E})u(-\mathsf{c}-\mathbf{E}) =0\  \Leftrightarrow \  \mfrac{v(-\mathsf{c}-\mathbf{E})}{u(-\mathsf{c}-\mathbf{E})}=\mfrac{v(-\mathbf{E})}{u(-\mathbf{E})},$$
\item for $ \mathbf{E}\notin \{ \tilde{c}_{2j+1} -\mathsf{c} \}_{j\geq 0}\cup \{  \tilde{c}_{2j} \}_{j\geq 0}$,
$$u(-\mathbf{E})v'(-\mathsf{c}-\mathbf{E})-v(-\mathbf{E})u'(-\mathsf{c}-\mathbf{E}) =0\  \Leftrightarrow \ \mfrac{v'(-\mathsf{c}-\mathbf{E})}{u'(-\mathsf{c}-\mathbf{E})}=\mfrac{v(-\mathbf{E})}{u(-\mathbf{E})},$$
\item for $\mathbf{E}\notin  \{  \tilde{c}_{2j} -\mathsf{c} \}_{j\geq 0}\cup  \{ \tilde{c}_{2j+1} \}_{j\geq 0}$,
$$ u'(-\mathbf{E})v(-\mathsf{c}-\mathbf{E})-v'(-\mathbf{E})u(-\mathsf{c}-\mathbf{E}) =0\  \Leftrightarrow \ \mfrac{v(-\mathsf{c}-\mathbf{E})}{u(-\mathsf{c}-\mathbf{E})}=\mfrac{v'(-\mathbf{E})}{u'(-\mathbf{E})}.$$
\end{enumerate}

\vskip 5mm

\noindent We look at the cases where the conditions on $-\mathsf{c}-\mathbf{E}$ and $-\mathbf{E}$ are not satisfied. Assume that:
$$\mathsf{c}\in \{ {\tilde c}_q-{\tilde c}_r,\ |\ q>r\geq 0\}:=\mathsf{Z}.$$
Then, for $\mathsf{c} ={\tilde c}_q-{\tilde c}_r  \in \mathsf{Z}$, in the set $\{ c_p,\tilde{c}_p \}_{p\geq 0}$, $\mathbf{E}=\tilde{c}_r$ is the unique  solution of the equation: (\ref{eq_band_lin1}) if  $q$ and $r$ are odd, (\ref{eq_band_lin2}) if  $q$ and $r$ are even, (\ref{eq_band_lin3}) if  $q$ is odd and $r$ is even and (\ref{eq_band_lin4}) if $q$ is even and $r$ is odd.

\noindent Conversely, if  $\mathsf{c}\notin \mathsf{Z}$, then none of the $\tilde{c}_p$ is a solution in $-\mathbf{E}$ of any of the equations (\ref{eq_band_lin1}), (\ref{eq_band_lin2}), (\ref{eq_band_lin3})  and (\ref{eq_band_lin4}). 

\vskip 3mm

\noindent \textbf{Assumption.} From now on, we assume that $\mathsf{c}\notin \mathsf{Z}.$ 
\vskip 3mm

\noindent Note that, however, all our results hold true when $\mathsf{c}\in \mathsf{Z}$, this assumption is only made for convenience's sake. 
 
\vskip 3mm

\noindent Since we have the four equivalences above, the band edges of the rescaled spectral bands of $H$ are solutions of the four following equations:
\begin{alignat}{2}
\frac{v'}{u'}(-\mathsf{c}-\mathbf{E})&& =\frac{v'}{u'}(-\mathbf{E}), & \ \mbox{ for }\ \mathbf{E}\notin \{  \tilde{c}_{2j+1} -\mathsf{c} \}_{j\geq 0}\cup \{  \tilde{c}_{2j+1} \}_{j\geq 0},\label{eq_vp_Ha_d_a} \\
\frac{v}{u}(-\mathsf{c}-\mathbf{E})&& =\frac{v}{u}(-\mathbf{E}), & \ \mbox{ for }\  \mathbf{E} \notin  \{ \tilde{c}_{2j}-\mathsf{c} \}_{j\geq 0}\cup  \{ \tilde{c}_{2j} \}_{j\geq 0},\label{eq_vp_Ha_d_b} \\
\frac{v'}{u'}(-\mathsf{c}-\mathbf{E})&& =\frac{v}{u}(-\mathbf{E}), & \  \mbox{ for }\  \mathbf{E}\notin \{ \tilde{c}_{2j+1} -\mathsf{c} \}_{j\geq 0}\cup \{ \tilde{c}_{2j} \}_{j\geq 0},\label{eq_vp_Ha_0_b} \\
\frac{v}{u}(-\mathsf{c}-\mathbf{E})&& =\frac{v'}{u'}(-\mathbf{E}),& \  \mbox{ for }\   \mathbf{E}\notin \{ \tilde{c}_{2j}-\mathsf{c} \}_{j\geq 0}\cup \{ \tilde{c}_{2j+1} \}_{j\geq 0}.\label{eq_vp_Ha_0_a} 
\end{alignat}

\subsection{Variations of $\frac{v}{u}$ and $\frac{v'}{u'}$}
\vskip 5mm

\begin{figure}
\centering
\resizebox{0.9\textwidth}{!}{%
  \includegraphics{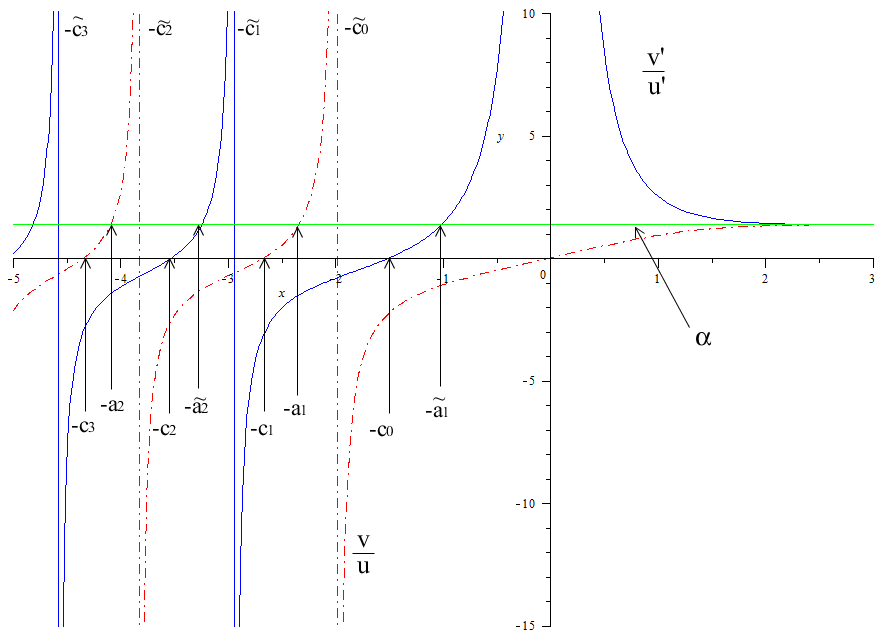}
}
\caption{Graphs of $\frac{v}{u}$ and $\frac{v'}{u'}$.}
\label{fig_usurv}       
\end{figure}

\vskip 5mm

\noindent Using the value of the Wronskian of $u$ and $v$ one has:
$$\forall x\in [0,+\infty),\ \left( \frac{v}{u} \right)'(x)=\frac{v'(x)u(x)-v(x)u'(x)}{u^2(x)} =\frac{1}{u^2(x)}>0$$
and
$$\forall x\in (0,+\infty),\ \left( \frac{v'}{u'} \right)'(x)=\frac{xv(x)u'(x)-xv'(x)u(x)}{(u'(x))^2} = -\frac{x}{(u'(x))^2}<0.$$

\noindent Thus, the functions $\frac{v}{u}$ and $\frac{v'}{u'}$ have the following behaviour. Let $j\geq 0$.
\begin{enumerate}
\item On every interval $(-\tilde{c}_{2j+2},-\tilde{c}_{2j})$, the function $\frac{v}{u}$ is continuous, strictly increasing and is a bijection from $(-\tilde{c}_{2j+2},-\tilde{c}_{2j})$ to $\R$. The function $\frac{v}{u}$ is continuous, strictly increasing and is a bijection from $(-\tilde{c}_0,+\infty)$ to $(-\infty, \alpha)$.
\item On every interval $(-\tilde{c}_{2j+3},-\tilde{c}_{2j+1})$, $\frac{v'}{u'}$ is continuous, strictly increasing and is a bijection from $(-\tilde{c}_{2j+3}, -\tilde{c}_{2j+1})$ to $\R$. The function $\frac{v'}{u'}$ is continuous and strictly increasing on $(-\tilde{c}_1,0)$ from $-\infty$ to $+\infty$. It is also continuous, strictly decreasing and a bijection from $(0,+\infty)$ to $(\alpha,+\infty)$.
\end{enumerate}

\noindent We remark that $\alpha$ is the common limit at infinity of the two functions $\frac{v}{u}$ and $\frac{v'}{u'}$, thanks to the limits $\frac{Bi(x)}{Ai(x)} \xrightarrow[x\to +\infty]{} +\infty$ and  $\frac{Bi'(x)}{Ai'(x)} \xrightarrow[x\to +\infty]{} -\infty$. 

\vskip 4mm

\noindent The limit $\alpha$ also provides lower bounds for the derivatives of $\frac{v}{u}$ and $\frac{v'}{u'}$ on the negative half-line.

\begin{lem}\label{lem_minor_uv_alpha}
For every $x<0$,
\begin{equation}\label{eq_lem_uv_alpha}
 \frac{1}{u^2(x)} > \frac{\alpha}{2}\qquad \mbox{ and } \qquad -\frac{x}{(u'(x))^2} > \alpha.
\end{equation}
\end{lem}

\begin{demo}
Using (\ref{eq_expr_u_Ai}), the equality $Bi'(0)=\sqrt{3}Ai'(0)$ and the expression of $\alpha$, one has:
\begin{equation}\label{eq_demo_lem_uv_alpha1}
\forall x \leq 0,\   \frac{1}{u^2(x)} = \frac{\alpha}{2}\cdot \frac{2\alpha}{\pi^2 (Ai(0))^2}\cdot \frac{1}{(\sqrt{3} Ai(x)+Bi(x))^2}.
\end{equation}
Let $M=\sqrt{Ai^2+Bi^2}$ be the Airy modulus. The function $M$ is strictly increasing on $(-\infty,0]$ (see \cite{O}) and we have:
$$\forall x\leq 0,\ (\sqrt{3} Ai(x)+Bi(x))^2 \leq 4(M(x))^2 \leq 4(M(0))^2.$$
With (\ref{eq_demo_lem_uv_alpha1}), it leads to 
$$\forall x \leq 0,\ \frac{1}{u^2(x)} \geq \frac{\alpha}{2}\cdot \frac{2\alpha}{\pi^2 (Ai(0))^2}\cdot \frac{1}{4((Ai(0))^2+(Bi(0))^2)}.$$
Since $\frac{2\alpha}{\pi^2 (Ai(0))^2}\cdot \frac{1}{4((Ai(0))^2+(Bi(0))^2)} \simeq 1.09 $, one deduces the first lower bound in (\ref{eq_lem_uv_alpha}).

\vskip 3mm

\noindent Derivating (\ref{eq_expr_u_Ai}) and using  the equality $Bi'(0)=\sqrt{3}Ai'(0)$ and the expression of $\alpha$, one has:
\begin{equation}\label{eq_demo_lem_uv_alpha2}
\forall x < 0,\   -\frac{x}{(u'(x))^2} = \alpha\cdot \frac{\alpha}{\pi^2 (Ai(0))^2}\cdot \frac{-x}{(\sqrt{3} Ai'(x)+Bi'(x))^2}.
\end{equation}
Let $N=\sqrt{(Ai')^2+(Bi')^2}$. Then, $x\mapsto |x|^{-\frac14} N(x)$ is strictly increasing on $(-\infty,0)$ (see \cite{O}) and we have :
$$\forall x< 0,\ \frac{-x}{(\sqrt{3} Ai'(x)+Bi'(x))^2} \geq  \frac{-x\cdot |x|^{-\frac12}}{4(N(x))^2 |x|^{-\frac12} } .$$
Since $x\mapsto \frac{|x|^{\frac12}}{ (N(x))^2 |x|^{-\frac12}}$ is strictly decreasing on $(-\infty,0)$, one has:
$$\forall x < -\tilde{a}_1,\  \frac{-x}{(\sqrt{3} Ai'(x)+Bi'(x))^2} \geq \frac{\tilde{a}_1}{4(N(-\tilde{a}_1 ))^2 }. $$
As $(N(-\tilde{a}_1 ))^2= |Bi'(-\tilde{a}_1)|$, with (\ref{eq_demo_lem_uv_alpha2}) one has:
$$\forall x <  -\tilde{a}_1, \ -\frac{x}{(u'(x))^2} \geq \alpha\cdot \frac{\alpha}{\pi^2 (Ai(0))^2}\cdot \frac{\tilde{a}_1}{4|Bi'(-\tilde{a}_1)|}.$$
With $\frac{\alpha}{\pi^2 (Ai(0))^2}\cdot \frac{\tilde{a}_1}{4|Bi'(-\tilde{a}_1)|} \simeq 2.74$, we deduce that 
$$\forall x <  -\tilde{a}_1, \ -\frac{x}{(u'(x))^2} > \alpha.$$
Moreover, the function $x\mapsto -\frac{x}{(u'(x))^2}$ is convex on the interval $(-\tilde{c}_1,0)$ and its minimum on this interval has an approximate value equal to $2.03> \alpha$. Since $-\tilde{c}_1 < -\tilde{a}_1$, we deduce the second lower bound in (\ref{eq_lem_uv_alpha}).
\end{demo}

\subsection{Some auxilliary functions}\label{sec_2analytic}

\noindent One sets, for $x\geq 0$ and $z\in \R$, 
{\small \begin{equation}\label{eq_def_fx}
  f_x(z)=v'(x-z)u(x)-u'(x-z)v(x) = \pi \left( Bi'(x-z)Ai(x)-Ai'(x-z)Bi(x)  \right)
\end{equation}}
and 
{\small \begin{equation}\label{eq_def_gx}
 g_x(z)= v(x-z)u(x)-u(x-z)v(x)= \pi \left( Bi(x-z)Ai(x)-Ai(x-z)Bi(x) \right).
\end{equation}}
\noindent The expressions  in terms of the Airy functions allow us to use classical properties of the $Ai$ and $Bi$ functions instead of the properties of $u$ and $v$ when it makes proofs easier.
\vskip 2mm

\noindent The functions $f_x$ and $g_x$ are non-zero solutions of differential equations which satisfy the assumptions of Sturm's theorem, thus their zeroes are isolated on the real line. We denote by 
$$z_0(x) < z_2(x) < \cdots < z_{2j}(x) < \ldots$$
the zeroes of $f_x$ arranged in increasing order.
Then, since $0$ is the first zero of $g_x$ for every $x$, we denote by 
$$0 < z_1(x) < z_3(x) < \cdots < z_{2j+1}(x) < \ldots$$
the zeroes of $g_x$ arranged in increasing order.

\noindent We characterize these zeroes and prove that none of them is negative. 
\vskip 2mm

\noindent Let $j\geq 0$ an integer. Let $x\geq 0$ and denote by $\psi_{2j}(x)$ the unique solution of the equation
\begin{equation}\label{eq_psij1}
\frac{v'}{u'}(z)=\frac{v}{u}(x),\ z\in [-c_{2j},-\tilde{a}_{j+1}). 
\end{equation}
We also denote by $\psi_{2j+1}(x)$ the unique solution of the equation
\begin{equation}\label{eq_psij2}
\frac{v}{u}(z)=\frac{v}{u}(x),\ z\in [-c_{2j+1},-a_{j+1}). 
\end{equation}

\begin{lem}\label{lem_psi_j}
For every $k\geq 0$, the function $\psi_k$ is well defined, continuous and strictly increasing. 
\end{lem}

\begin{demo}
Let $j\geq 0$. Recall that $\alpha$ is the common limit of $\frac{v}{u}$ and $\frac{v'}{u'}$ at infinity. The function $\frac{v'}{u'}$ is a bijection from $[-c_{2j},-\tilde{a}_{j+1})$ to $[0,\alpha)$ and we denote by 
$$\left( \frac{v'}{u'} \right)_{2j}^{-1} \ : \ [0,\alpha) \to [-c_{2j},-\tilde{a}_{j+1})$$
its reciprocal function. Since for every $x\geq 0$, $\frac{v}{u}(x)\in [0,\alpha)$, we have:
$$\forall x \geq 0,\ \psi_{2j}(x) = \left( \frac{v'}{u'} \right)_{2j}^{-1} \left( \frac{v}{u}(x) \right).$$
Thus, the function $\psi_{2j}$ is well defined and it is continuous by continuity of $\frac{v}{u}$ on $[0,+\infty)$ and of the inverse of $\frac{v'}{u'}$ on $ [0,\alpha)$. Since $\frac{v'}{u'}$ is strictly increasing on $[-c_{2j},-\tilde{a}_{j+1})$, its reciprocal function is strictly increasing on $[0,\alpha)$ and since $\frac{v}{u}$ is strictly increasing on $[0,+\infty)$, we deduce that  $\psi_{2j}$ is strictly increasing on $[0,+\infty)$. 

\noindent The function $\frac{v}{u}$ is a bijection from $(-c_{2j+1},-a_{j+1}]$ to $[0,\alpha)$ and we denote by 
$$\left( \frac{v}{u} \right)_{2j+1}^{-1} \ : \ [0,\alpha) \to [-c_{2j+1},-a_{j+1})$$
its reciprocal function. With the same arguments as before, we have that 
$$\forall x \geq 0,\ \psi_{2j+1}(x) = \left( \frac{v}{u} \right)_{2j+1}^{-1} \left( \frac{v}{u}(x) \right)$$
and thus $\psi_{2j+1}$ is well defined, continuous and strictly increasing.
\end{demo}

\begin{lem}\label{lem_zeroes}
Let $k\geq 0$. Then, for every $x\geq 0$,
$$z_k(x)\geq 0\quad \mbox{ and }\quad z_k(x)=x-\psi_k(x).$$
Therefore, $z_k$ is continuous on $[0,+\infty)$. Moreover, for every $k\geq 0$, the function $z_k$ is strictly increasing from $[0,+\infty)$ to $[c_k,+\infty)$.
\end{lem}

\noindent This result is proven in Appendix \ref{sec_zk}.

\vskip 3mm

\noindent Let $j\geq 0$. For $x\geq 0$, denote by $\psi^{2j}(x) \in (-a_{j+1}, -\tilde{c}_{2j}]$ the unique solution of the equation
\begin{equation}\label{eq_Eomax1_even}
\frac{v}{u}(z)=\frac{v'}{u'}(x),\ z\in (-a_{j+1},-\tilde{c}_{2j}]. 
\end{equation}
\noindent We also denote by $\psi^{2j+1}(x) \in (-\tilde{a}_{j+2},-\tilde{c}_{2j+1})$ the unique solution of the equation
\begin{equation}\label{eq_Eomax1_odd}
\frac{v'}{u'}(z)=\frac{v'}{u'}(x),\ z\in (-\tilde{a}_{j+2},-\tilde{c}_{2j+1}). 
\end{equation}

\begin{lem}\label{lem_Eomax1}
For every $k\geq 0$, the function $\psi^k$ is well defined, continuous and strictly decreasing on $[0,+\infty)$.  
\end{lem}

\begin{demo} We assume that $k$ is even, $k=2j$ for one $j\geq 0$. The function $\frac{v}{u}$ is a continuous bijection from $(-a_{j+1},-\tilde{c}_{2j}]$ to $[\alpha,+\infty)$ and we denote its reciprocal function by 
$$\left( \frac{v}{u} \right)_{2j}^{-1}\ : \ [\alpha,+\infty) \to (-a_{j+1},-\tilde{c}_{2j}].$$
Since for every $x\geq 0$, $\frac{v'}{u'}(x)\in [\alpha,+\infty)$, we have:
$$\forall x\geq 0,\ \psi^{2j}(x)=\left(\frac{v}{u}\right)_{2j}^{-1} \left( \frac{v'}{u'}(x)\right).$$
Then $\psi^{2j}$ is well defined, continuous and it is the unique solution of (\ref{eq_Eomax1_even}). Moreover, $\left(\frac{v}{u}\right)_{2j}^{-1}$ is a strictly increasing function from $ [\alpha,+\infty)$ to  $(-a_{j+1},-\tilde{c}_{2j}]$ and $\frac{v'}{u'}$ is strictly decreasing on $[0,+\infty)$. Thus, $\psi^{2j}$ is strictly decreasing on $[0,+\infty)$.
\vskip 2mm

\noindent  We assume that $k$ is odd, $k=2j+1$ for one $j\geq 0$. The function $\frac{v'}{u'}$ is a continuous bijection from $(-\tilde{a}_{j+2},-\tilde{c}_{2j+1})$ to $[\alpha,+\infty)$ and we denote its reciprocal function by 
$$\left( \frac{v'}{u'} \right)_{2j+1}^{-1}\ : \ [\alpha,+\infty) \to (-\tilde{a}_{j+2},-\tilde{c}_{2j+1}).$$
Since for every $x\geq 0$, $\frac{v'}{u'}(x)\in [\alpha,+\infty)$, we write:
$$\forall x\geq 0,\ \psi^{2j+1}(x)=\left(\frac{v'}{u'}\right)_{2j+1}^{-1} \left( \frac{v'}{u'}(x)\right).$$
Then $\psi^{2j+1}$ is well defined, continuous and it is the unique solution of (\ref{eq_Eomax1_odd}). Moreover, $\left(\frac{v'}{u'}\right)_{2j+1}^{-1}$ is a strictly increasing function from $[\alpha,+\infty)$ to $(-\tilde{a}_{j+2}, -\tilde{c}_{2j+1})$ and $\frac{v'}{u'}$ is strictly decreasing on $[0,+\infty)$. Thus, $\psi^{2j+1}$ is strictly decreasing on $[0,+\infty)$.
\end{demo}

\vskip 3mm

\section{The first spectral band}\label{sec_band1}

\subsection{Lower bound of the continuous spectrum}

 For values of $E$ such that $\mathbf{E}<-\mathsf{c}$, we show that there is no solutions to equations (\ref{eq_vp_Ha_d_a}), (\ref{eq_vp_Ha_d_b}), (\ref{eq_vp_Ha_0_a}) and (\ref{eq_vp_Ha_0_b}). This is natural, from a physical point of view, since in this case, the energy is smaller than the minimum of the potential and it is as if the potential was not ``seen'' at this energy. This writes:

\begin{prop}\label{prop_spec}
For every $\mathsf{c}\geq 0$, $\sigma(H) \subset [-V_0,+\infty)$. 
\end{prop}

\begin{demo} We remark that on the interval $[0,+\infty)$, $\frac{v}{u}< \alpha$ and $\frac{v'}{u'}>\alpha$. Thus, these two functions cannot have a common value on this interval and equations  (\ref{eq_vp_Ha_0_a}) and (\ref{eq_vp_Ha_0_b}) do not have any solution with $0<-\mathsf{c}-\mathbf{E}<-\mathbf{E}$.

\noindent Moreover, we already know that $\frac{v}{u}$ is strictly increasing on $[0,+\infty)$ and $\frac{v'}{u'}$ is strictly decreasing on $(0,+\infty)$. Since $-\mathsf{c}-\mathbf{E} \neq -\mathbf{E}$ the equations (\ref{eq_vp_Ha_d_a}), (\ref{eq_vp_Ha_d_b}) do not have any solution when $0<-\mathsf{c}-\mathbf{E}<-\mathbf{E}$.
\end{demo}
\vskip 3mm

\noindent \textbf{Remark.} This result holds true for every $\mathsf{c}$ strictly positive. In particular we do not need to assume the semiclassical parameter $\mathsf{h}$ to be small.

\subsection{The bottom of the spectrum}\label{sec_bottom}

\noindent At first, we determine the bottom of the spectrum. It has to be a solution of either the equation (\ref{eq_vp_Ha_d_a}) or the equation (\ref{eq_vp_Ha_d_b}) with $-\mathbf{E}\geq 0$ and thus $-\mathsf{c}-\mathbf{E} <0$.
\vskip 2mm

\noindent We start by proving that for every $\mathsf{c} >0$, the equation (\ref{eq_vp_Ha_d_a}) has a unique solution with $-\mathbf{E}>0$ and $-\mathsf{c}-\mathbf{E} \in [-\tilde{a}_1,0)$. The function $\frac{v'}{u'}$ is an increasing continuous bijection from $(-\tilde{a}_1,0)$ to $[\alpha,+\infty)$ and thus, since for every $x>0$, $\frac{v'}{u'}(x)\in [\alpha,+\infty)$, one define:
\begin{equation}\label{def_psi}
\forall x>0,\ \psi(x)=\left(\frac{v'}{u'}\Big|_{(-\tilde{a}_1,0)} \right)^{-1} \left( \frac{v'}{u'}(x) \right).
\end{equation}
The function $\psi$ does not belong to the family of the functions $\psi_{2k+1}$, due to the difference of behaviour of $\frac{v'}{u'}$ on $(-\tilde{a}_1,+\infty)$ compared to $(-\infty,-\tilde{a}_1)$.

\noindent The function $\psi$ is a continuous decreasing bijection from $(0,+\infty)$ to $(-\tilde{a}_1,0)$ and $x\mapsto x-\psi(x)$ is a continuous increasing bijection from $(0,+\infty)$ to $(0,+\infty)$. Thus,
$$\forall \mathsf{c} >0,\ \exists! \tilde{x}>0,\ \tilde{x}-\psi(\tilde{x})=\mathsf{c}.$$
One sets $\tilde{\mathbf{E}}=-\tilde{x}$ and since $\tilde{x}=\mathsf{c} +\psi(\tilde{x}) \in (-\tilde{a}_1 + \mathsf{c}, \mathsf{c})$ one has
\begin{equation}\label{eq_bottom_estim1}
-\mathsf{c} < \tilde{\mathbf{E}} < -\mathsf{c} + \tilde{a}_1. 
\end{equation}
We remark that $\tilde{\mathbf{E}}$ is the smallest solution of (\ref{eq_vp_Ha_d_a}).

\vskip 2mm

\noindent Now we turn to the smallest solution of (\ref{eq_vp_Ha_d_b}). Since the function $z_1$ is an increasing continuous bijection from $[0,+\infty)$ to $[c_1,+\infty)$,
$$ \forall \mathsf{c} \geq c_1, \exists! x_1\geq 0,\ z_1(x_1)=x_1-\psi_1(x_1) = \mathsf{c} .$$
One sets $\breve{\mathbf{E}}_1=-x_1$ and since $x_1=\mathsf{c} +\psi_1(x_1) \in (-c_1 + \mathsf{c}, -a_1 + \mathsf{c})$ one has
$$-\mathsf{c}+ a_1 < \breve{\mathbf{E}}_1 < -\mathsf{c} + c_1.$$
\vskip 2mm

\noindent Since $\tilde{a}_1 < a_1$ we have $\tilde{\mathbf{E}}< \breve{\mathbf{E}}_1$ and thus  $\mathbf{E}_{\mathrm{min}}^0= \tilde{\mathbf{E}}$. In particular, $\mathbf{E}_{\mathrm{min}}^0$ is the smallest solution of (\ref{eq_vp_Ha_d_a}). Moreover, it is the smallest solution among those of (\ref{eq_vp_Ha_d_a}) and (\ref{eq_vp_Ha_d_b}).
\vskip 2mm

\noindent Now that we have identified the bottom of the spectrum of $H$, we prove the estimates stated in Theorem \ref{thm_bottom}.

\vskip 3mm
\begin{demo} (of Theorem \ref{thm_bottom}).

\noindent \textbf{(1)} Since $\mathbf{E}_{\mathrm{min}}^0= \tilde{\mathbf{E}}$, by (\ref{eq_bottom_estim1}) we have, for every $\mathsf{c}>0$, $-\mathsf{c} < \mathbf{E}_{\mathrm{min}}^0 <  -\mathsf{c} + \tilde{a}_1$.
\vskip 2mm
\noindent Let $\mathsf{c} >0$. By the previous inequality, $-\mathbf{E}_{\mathrm{min}}^0 \in (-\tilde{a}_1 + \mathsf{c}, \mathsf{c})$. We recall that, since $u$ and $v$ are linear combinations of $Ai$ and $Bi$ and since $\mathbf{E}_{\mathrm{min}}^0$ is a solution of (\ref{eq_vp_Ha_d_a}), it also satisfies:
$$\mfrac{Bi'}{Ai'}(-\mathbf{E}_{\mathrm{min}}^0 - \mathsf{c}) =\mfrac{Bi'}{Ai'}(-\mathbf{E}_{\mathrm{min}}^0),\quad -\mathbf{E}_{\mathrm{min}}^0 - \mathsf{c}\in (-\tilde{a}_1,0). $$
We introduce the functions
$$F_{\pm}(\cdot, \mathsf{c})\ :\ \begin{array}{ccl}
                                 I_{\pm} & \to & \R \\
x  & \mapsto & \frac{Bi'}{Ai'}\left(x-\frac{\mathsf{c}}{2} \right)-\frac{Bi'}{Ai'}\left(x+\frac{\mathsf{c}}{2} \right)
                                 \end{array}
$$
with $I_{-}=\left(-\frac{\mathsf{c}}{2},\frac{\mathsf{c}}{2} \right)$ and  $I_{+}=\left(-\tilde{a}_1 +\frac{\mathsf{c}}{2},\frac{\mathsf{c}}{2} \right)$. We remark that, for every $\mathsf{c} >0$, the unique zero of $F_{\pm}(\cdot, \mathsf{c})$ on $I_{\pm}$ is  $-\mathbf{E}_{\mathrm{min}}^0-\frac{\mathsf{c}}{2}$.

\noindent Thanks to the fact that $0$ is the point of maximum of $\frac{Bi'}{Ai'}$ on $(-\tilde{a}_1,+\infty)$, one has:
$$F_{\pm}\left(\mfrac{\mathsf{c}}{2}, \mathsf{c}\right) =\mfrac{Bi'}{Ai'}(0)- \mfrac{Bi'}{Ai'}(\mathsf{c}) >0.$$
\vskip 2mm

\noindent \emph{First case:} $\mathsf{c} \geq  2\tilde{a}_1$. In this case, $-\tilde{a}_1 +  \frac{\mathsf{c}}{2} >0$ hence the unique zero of $F_{+}(\cdot, \mathsf{c})$ on $I_+$ is strictly positive, hence $-\mathbf{E}_{\mathrm{min}}^0-\frac{\mathsf{c}}{2} > 0$ which implies that $\mathbf{E}_{\mathrm{min}}^0 < -\frac{\mathsf{c}}{2}$. Thus, point $(1)$ is proved for every $\mathsf{c} \geq  2\tilde{a}_1$.
\vskip 2mm

\noindent \emph{Second case:} $0< \mathsf{c} < 2 \tilde{a}_1$. We prove that the unique zero of  $F_{\pm}(\cdot, \mathsf{c})$ on $I_{\pm}$ is strictly positive. Since $F_{\pm}(\frac{\mathsf{c}}{2}, \mathsf{c})>0$, it is sufficient to prove that $F_{\pm}(0, \mathsf{c})<0$. Indeed, if $0< \mathsf{c} < \tilde{a}_1$, then $0\in I_{-}$ and we study the unique zero of $F_{-}$ in $I_{-}$. If $\tilde{a}_1 \leq \mathsf{c} < 2\tilde{a}_1$, then $-\tilde{a}_1 + \frac{\mathsf{c}}{2} <0$ and $0\in I_+$. In this case, we study the unique zero of $F_+$ in $I_+$. 
\vskip 2mm

\noindent We have:
$$F_{\pm}(0, \mathsf{c}) = \mfrac{Bi'}{Ai'}\left(-\mfrac{\mathsf{c}}{2} \right)-\mfrac{Bi'}{Ai'}\left(\mfrac{\mathsf{c}}{2} \right).$$
Let $y=\frac{\mathsf{c}}{2}$ so that $y\in (0,\tilde{a}_1)$ and set:
$$\forall y\in (0,\tilde{a}_1),\ G(y)= \mfrac{Bi'}{Ai'}(-y)-\mfrac{Bi'}{Ai'}(y). $$
One has, for every $y\in (0,\tilde{a}_1)$,
$$G'(y)  = \frac{y}{\pi (Ai'(-y))^2 (Ai'(y))^2  } \left( (Ai'(-y) - Ai'(y))(Ai'(-y) + Ai'(y) ) \right). $$
On $(0,\tilde{a}_1)$, $\frac{y}{\pi (Ai'(-y))^2 (Ai'(y))^2  } >0$ and as $Ai'$ is negative on $(-\tilde{a}_1,+\infty)$, $Ai'(-y) + Ai'(y) <0$ for every $y\in (0,\tilde{a}_1)$. 

\noindent Let 
$$K\ :\ \begin{array}{ccl}
          (0,\tilde{a}_1) & \to & \R \\
y  & \mapsto & Ai'(-y)-Ai'(y)
\end{array}.$$
Then $K(0)=0$ and
$$\forall y\in (0,\tilde{a}_1),\ K'(y)=-Ai''(-y)-Ai''(y) = y(Ai(-y)-Ai(y)).$$
But, the Airy function $Ai$ is decreasing on $(-\tilde{a}_1,+\infty)$ hence, for $y\in (0,\tilde{a}_1)$, $Ai(-y)-Ai(y) >0$ and $K'(y) >0$. Thus, $K$ is strictly increasing on $(0,\tilde{a}_1)$ and
$$\forall y\in (0,\tilde{a}_1),\ K(y) > K(0)=0.$$
Thus, for every $ y\in (0,\tilde{a}_1)$, $G'(y) <0$. Since $G(0)=0$, for every $ y\in (0,\tilde{a}_1)$, $G(y) <0$, which rewrites
$$\forall \mathsf{c} \in (0,2\tilde{a}_1),\ F_{\pm}(0, \mathsf{c}) <0.$$
Thus, the unique zero of $F_{\pm}(\cdot, \mathsf{c})$ in $I_{\pm}$ is strictly positive. Thus, $-\mathbf{E}_{\mathrm{min}}^0-\frac{\mathsf{c}}{2} \in \left(0,\frac{\mathsf{c}}{2}  \right)$ and $\mathbf{E}_{\mathrm{min}}^0 < -\frac{\mathsf{c}}{2}$. 
\vskip 2mm

\noindent Taking in account the result in both cases, point $(1)$ is proved.
\vskip 5mm

\noindent \textbf{(2)} Since $\mathbf{E}_{\mathrm{min}}^0$ is the smallest solution of (\ref{eq_vp_Ha_d_a}), if one sets $X=-\mathbf{E}_{\mathrm{min}}^0 - \mathsf{c} + \tilde{a}_1$, then $X$ satisfies
\begin{equation}\label{eq_demo_bottom_1}
\frac{v'}{u'}(X-\tilde{a}_1) =  \frac{v'}{u'}(X + \mathsf{c}-\tilde{a}_1), \quad X\in [0,\tilde{a}_1).
\end{equation}

\noindent Since $x\mapsto x-\psi(x)$ is increasing on $(0,+\infty)$, when $\mathsf{c}$ tends to infinity, $X-\tilde{a}_1$ tends to $-\tilde{a}_1$ and $X+\mathsf{c}-\tilde{a}_1$ tends to $+\infty$. Therefore, both members of (\ref{eq_demo_bottom_1}) tends to $\alpha$ when $\mathsf{c}$ tends to $+\infty$. Let $\tau >0$ be arbitrary small. We assume that $\mathsf{c} \geq \tilde{a}_1 + \tau$.

\noindent Using \cite[10.4.61, 10.4.66]{AS} in the equality (\ref{eq_demo_bottom_1}), one gets

\begin{equation}\label{eq_demo_bottom_4}
 \frac{v'}{u'}(X-\tilde{a}_1)= \alpha + \alpha \sqrt{3} \ex^{-\frac43 (X + \mathsf{c} - \tilde{a}_1 )^{\frac32}} \left(1+ \mathcal{O}\left(  \left(X + \mathsf{c} - \tilde{a}_1 \right)^{-\frac32} \right) \right).
\end{equation}

\noindent Let $\epsilon = \ex^{-\frac43 (\mathsf{c}-\tilde{a}_1)^{\frac32}}$, which amounts to $\mathsf{c}-\tilde{a}_1 = \left( -\frac34 \ln(\epsilon)\right)^{\frac23}$. In particular, $\mathsf{c}$ tends to $+\infty$ if and only if $\epsilon$ tends to $0^+$.

\noindent  We use the following identity valid for strictly positive real numbers $a$ and $b$:
\begin{equation}\label{eq_demo_bottom_ab}
 a^{\frac32}-b^{\frac32} - \tfrac32 b^{\frac12} (a-b) = (a-b)^2 \frac{2\sqrt{\frac{a}{b} } +1 }{2\left( \sqrt{\frac{a}{b} } +1  \right)^2  } b^{-\frac12}
\end{equation}
with $a=X+\mathsf{c}-\tilde{a}_1 = \epsilon Y +\mathsf{c}-\tilde{a}_1$ and $b=\mathsf{c}-\tilde{a}_1 =  \left( -\frac34 \ln(\epsilon)\right)^{\frac23}$.
Then, 
\begin{equation}
\begin{split}
 -\tfrac43 \left(X+\mathsf{c}-\tilde{a}_1\right)^{\frac32} +  \tfrac43 \left( \left( -\tfrac34 \ln(\epsilon)\right)^{\frac23} \right)^{\frac32} & + 2  \left( -\tfrac34 \ln(\epsilon)\right)^{\frac13} \epsilon Y  \nonumber \\
& = \left( -\tfrac34 \ln(\epsilon)\right)^{-\frac13} (\epsilon Y)^2 Q(\epsilon, Y)  
\end{split}
\end{equation}
where 
$$Q(\epsilon, Y) =\frac43 \frac{2 \frac{(\epsilon Y + \mathsf{c} - \tilde{a}_1)^{\frac12}}{\left( -\frac34 \ln(\epsilon)\right)^{\frac13}  } +1 }{2 \left(   \frac{(\epsilon Y + \mathsf{c} - \tilde{a}_1)^{\frac12}}{\left( -\frac34 \ln(\epsilon)\right)^{\frac13} }+1 \right)^2  }. $$
The condition $\mathsf{c} \geq \tilde{a}_1 + \tau$ implies $0 < \epsilon \leq \ex^{-\frac43 \tau^{\frac32}}$. Note that in the sequel, all the proven estimates will be valid for every $\epsilon \in (0,\ex^{-\frac43 \tau^{\frac32}}]$, which shows that they are estimates and not asymptotic expansions.

\noindent Since $ \epsilon Y +\mathsf{c}-\tilde{a}_1 = -\mathbf{E}_{\mathrm{min}}^0 >0$,
$$ \frac{(\epsilon Y + \mathsf{c} - \tilde{a}_1)^{\frac12}}{\left( -\frac34 \ln(\epsilon)\right)^{\frac13}} \geq 0.$$
If $\varphi : \R_+ \to \R$ is defined, for every $x\in \R_+$, by $\varphi(x)=\frac{2x+1}{2(x+1)^2}$, then $\varphi$ is decreasing on $\R_+$, $\varphi(0)=\frac12$ and $\varphi$ tends to $0$ when $x$ tends to $+\infty$. Thus, for every $\epsilon\in (0,\ex^{-\frac43 \tau^{\frac32}}]$ and every $Y$,  
$$0 <  Q(\epsilon, Y) \leq \frac23.$$
We have,
\begin{equation}\label{eq_demo_bottom_8}
 \ex^{-\frac43 \left(X+\mathsf{c}-\tilde{a}_1\right)^{\frac32} } = \epsilon \cdot \ex^{-  2  \left( -\frac34 \ln(\epsilon)\right)^{\frac13} \epsilon Y +  \left( -\frac34 \ln(\epsilon)\right)^{-\frac13} (\epsilon Y)^2 Q(\epsilon, Y)}.
\end{equation}

\noindent The condition $X\in [0,\tilde{a}_1)$ implies that 
$$ \left( -\tfrac34 \ln(\epsilon)\right)^{\frac23} = -\tilde{a}_1+\mathsf{c}    \leq X+\mathsf{c} -\tilde{a}_1 \leq \tilde{a}_1 +  \left( -\tfrac34 \ln(\epsilon)\right)^{\frac23}.$$
Moreover,
\begin{equation}\label{eq_demo_bottom_13}
\left( \left( \mfrac34 \right)^{\frac23}  + \mfrac{\tilde{a}_1}{|\ln (\epsilon)|^{\frac23}  }\right)^{-\frac32} \frac{1}{|\ln (\epsilon)|} \leq \left(X + \mathsf{c} - \tilde{a}_1 \right)^{-\frac32} \leq  \mfrac43  \frac{1}{|\ln (\epsilon)|}.  
\end{equation}
Using (\ref{eq_demo_bottom_4}) and the identity $\alpha =  \frac{v'}{u'}(-\tilde{a}_1)$, which comes directly from the expressions of $u$ and $v$ in terms of $Ai$ and $Bi$, one gets, for every $\epsilon \in (0,\ex^{-\frac43 \tau^{\frac32}}]$,
\begin{equation}
\begin{split}
& \frac{1}{\epsilon} \left(  \frac{v'}{u'}(\epsilon Y-\tilde{a}_1)-\frac{v'}{u'}(-\tilde{a}_1)\right) =\label{eq_demo_bottom_5} \\
&  \alpha \sqrt{3} \ex^{-  2  \left( -\frac34 \ln(\epsilon)\right)^{\frac13} \epsilon Y } \ex^{  \left( -\frac34 \ln(\epsilon)\right)^{-\frac13} (\epsilon Y)^2 Q(\epsilon, Y)} \left( 1+  \mathcal{O}\left( \frac{1}{|\ln (\epsilon)|} \right) \right).
\end{split}
\end{equation}
To estimate the left member of (\ref{eq_demo_bottom_5}) one uses:
{\small $$\frac{v'}{u'}(\epsilon Y-\tilde{a}_1)-\frac{v'}{u'}(-\tilde{a}_1) = \int_0^1 \left(\frac{v'}{u'} \right)' (-\tilde{a}_1 +t \epsilon Y) \epsilon Y \dd t =  \int_0^1 \frac{\tilde{a}_1 -t \epsilon Y}{(u'(-\tilde{a}_1 +t \epsilon Y))^2} \epsilon Y \dd t.$$}
\noindent Using the lower bound for $-\frac{x}{(u'(x))^2}$ given in (\ref{eq_lem_uv_alpha}), 
\begin{equation}\label{eq_demo_bottom_9}
\forall \epsilon \in (0,\ex^{-\frac43 \tau^{\frac32}}],\  \left| \frac{v'}{u'}(\epsilon Y-\tilde{a}_1)-\frac{v'}{u'}(-\tilde{a}_1) \right| \geq \alpha \epsilon Y.
\end{equation}
Let $A>0$ a constant given by the term $\mathcal{O}\left(  \left(X + \mathsf{c} - \tilde{a}_1 \right)^{-\frac32} \right) $ in (\ref{eq_demo_bottom_4}). Using (\ref{eq_demo_bottom_5}), (\ref{eq_demo_bottom_13}) and (\ref{eq_demo_bottom_9}), there exists $B=\frac43 A >0$ such that,
{\small $$\forall \epsilon \in (0,\ex^{-\frac43 \tau^{\frac32}}],\ Y\leq \sqrt{3} \ex^{-  2  \left( -\frac34 \ln(\epsilon)\right)^{\frac13} \epsilon Y } \ex^{ \frac23  \left( -\frac34 \ln(\epsilon)\right)^{-\frac13} (\epsilon Y)^2  }\left( 1+  B \left( \frac{1}{|\ln (\epsilon)|} \right) \right), $$}
which rewrites
{\small $$\forall \epsilon \in (0,\ex^{-\frac43 \tau^{\frac32}}],\ Y \ex^{2  \left( -\frac34 \ln(\epsilon)\right)^{\frac13} \epsilon Y } \leq \sqrt{3} \ex^{\frac23  \left( -\frac34 \ln(\epsilon) \right)^{-\frac13} (\epsilon Y)^2} \left( 1+  B \left( \frac{1}{|\ln (\epsilon)|} \right) \right). $$}
But, $\epsilon Y=X\in [0,\tilde{a}_1)$ is bounded and thus the right member of the previous inequality is bounded by a constant independent of $\epsilon $ and $Y$. Thus, there exists $C>0$ such that
\begin{equation}\label{eq_demo_bottom_6}
\forall \epsilon \in (0,\ex^{-\frac43 \tau^{\frac32}}],\ Y \ex^{2  \left( -\frac34 \ln(\epsilon)\right)^{\frac13} \epsilon Y } \leq  \sqrt{3} \ex^{\frac23 \tau^{-\frac12} \tilde{a}_1^2 }\left(1+\frac34 B \tau^{-\frac32} \right) =:C
\end{equation}
and
$$\forall \epsilon \in (0,\ex^{-\frac43 \tau^{\frac32}}],\ \epsilon \left( -\tfrac34 \ln(\epsilon)\right)^{\frac13}  Y \ex^{2  \left( -\frac34 \ln(\epsilon)\right)^{\frac13} \epsilon Y } \leq  \epsilon \left( -\tfrac34 \ln(\epsilon)\right)^{\frac13} C.$$
The function $x\mapsto x\ex^{2x}$ is $C^1$ and strictly increasing on $\R$, let us denote by $h$ its reciprocal function which is also $C^1$. From (\ref{eq_demo_bottom_6}) one gets
$$\forall \epsilon \in (0,\ex^{-\frac43 \tau^{\frac32}}],\ \epsilon \left( -\tfrac34 \ln(\epsilon)\right)^{\frac13} Y \leq h\left(\epsilon \left( -\tfrac34 \ln(\epsilon)\right)^{\frac13} C \right).$$
Since $h$ is of class $C^1$, $h(0)=0$, $h'(0)=1$ and $\epsilon \left( -\frac34 \ln(\epsilon)\right)^{\frac13}\xrightarrow[\epsilon \to 0]{} 0$, there exists $D>0$ such that
$$\forall \epsilon \in (0,\ex^{-\frac43 \tau^{\frac32}}],\ h\left(\epsilon \left( -\tfrac34 \ln(\epsilon)\right)^{\frac13} C \right) \leq \epsilon \left( -\tfrac34 \ln(\epsilon)\right)^{\frac13} C + D \epsilon^2 \left( -\tfrac34 \ln(\epsilon)\right)^{\frac23}$$ 
and 
\begin{equation}\label{eq_demo_bottom_7}
\forall \epsilon \in (0,\ex^{-\frac43 \tau^{\frac32}}],\ Y \leq C + D \epsilon \left( -\tfrac34 \ln(\epsilon)\right)^{\frac13},
\end{equation}
from which we deduce that $Y$ is bounded. Thus,  $X=\mathcal{O}(\epsilon)$, namely
$$X= \mathcal{O} \left(  \ex^{-\frac43 (\mathsf{c}-\tilde{a}_1)^{\frac32}} \right)$$
and since $X=-\mathbf{E}_{\mathrm{min}}^0 - \mathsf{c} + \tilde{a}_1$, we already have
$$\mathbf{E}_{\mathrm{min}}^0 = -\mathsf{c}+\tilde{a}_1+ \mathcal{O}\left( e^{-\frac43(-\tilde{a}_1+\mathsf{c})^{\frac32}}\right).$$
We refine the estimate. Since $Y$ is bounded, $\epsilon Y$ tends to $0$ when $\epsilon$ tends to $0$ and one has,
$$\forall \epsilon \in (0,\ex^{-\frac43 \tau^{\frac32}}],\ \frac{v'}{u'}(\epsilon Y-\tilde{a}_1)-\frac{v'}{u'}(-\tilde{a}_1) = \left(\frac{v'}{u'} \right)' (-\tilde{a}_1 ) (\epsilon Y) + \mathcal{O}(\epsilon^2).$$
Using (\ref{eq_demo_bottom_7}) to prove that 
$$\forall \epsilon \in (0,\ex^{-\frac43 \tau^{\frac32}}],\ \ex^{-  2  \left( -\frac34 \ln(\epsilon)\right)^{\frac13} \epsilon Y } \ex^{ \frac{1}{2}  \left( -\frac34 \ln(\epsilon)\right)^{-\frac13} (\epsilon Y)^2 } = 1+ \mathcal{O}\left( \epsilon |\ln(\epsilon)|^{\frac13} \right),$$
one deduces from (\ref{eq_demo_bottom_5}) that, for every $\epsilon \in (0,\ex^{-\frac43 \tau^{\frac32}}]$,
$$\left(\frac{v'}{u'} \right)'(-\tilde{a}_1 )\cdot Y + \mathcal{O}(\epsilon) =  \alpha \sqrt{3}  \left(1+ \mathcal{O}\left( \frac{1}{|\ln (\epsilon)|} \right) +  \mathcal{O}\left( \epsilon |\ln(\epsilon)|^{\frac13} \right) \right)$$
hence
$$\forall \epsilon \in (0,\ex^{-\frac43 \tau^{\frac32}}],\ \left(\frac{v'}{u'} \right)'(-\tilde{a}_1 )\cdot Y  =  \alpha \sqrt{3}  \left( 1+ \mathcal{O}\left( \frac{1}{|\ln (\epsilon)|} \right)  \right). $$
Since $\left(\frac{v'}{u'} \right)' (-\tilde{a}_1 ) = \frac{\tilde{a}_1}{(u'(-\tilde{a}_1))^2}$, for every $\mathsf{c}\geq \tilde{a}_1+\tau$,
\begin{equation}\label{eq_demo_bottom_10}
 Y=\alpha \sqrt{3} \mfrac{(u'(-\tilde{a}_1))^2  }{\tilde{a}_1} + \mathcal{O} \left( (\mathsf{c} - \tilde{a}_1)^{-\frac32} \right). 
\end{equation}
Thus, for every $\mathsf{c}\geq\tilde{a}_1+\tau$,
$$X=\epsilon Y = \alpha \sqrt{3} \frac{(u'(-\tilde{a}_1))^2  }{ \tilde{a}_1 } \ex^{-\frac43 (\mathsf{c}-\tilde{a}_1)^{\frac32}} +  \mathcal{O} \left( (\mathsf{c} - \tilde{a}_1)^{-\frac32} \ex^{-\frac43 (\mathsf{c}-\tilde{a}_1)^{\frac32}} \right)  $$
and finally, using $\mathbf{E}_{\mathrm{min}}^0=-X - \mathsf{c} + \tilde{a}_1$, for every $\mathsf{c}\geq\tilde{a}_1 + \tau$,
$$\mathbf{E}_{\mathrm{min}}^0 = -\mathsf{c}+\tilde{a}_1 - \alpha \sqrt{3}\mfrac{(u'(-\tilde{a}_1))^2  }{ \tilde{a}_1 } \ex^{-\frac43 (\mathsf{c}-\tilde{a}_1)^{\frac32}} +  \mathcal{O} \left( (\mathsf{c} - \tilde{a}_1)^{-\frac32} \ex^{-\frac43 (\mathsf{c}-\tilde{a}_1)^{\frac32}} \right).$$

\noindent Recall that
$$\ex^{-\frac43 (\mathsf{c}-\tilde{a}_1)^{\frac32}}\left(1+ \mathcal{O}\left((\mathsf{c}-\tilde{a}_1)^{-\frac32} \right) \right) =\ex^{-\frac43 (\mathsf{h}^{-\frac23}-\tilde{a}_1)^{\frac32}}\left(1+ \mathcal{O}\left(\mathsf{h}\right) \right). $$
As 
$$(1-\tilde{a}_1 \mathsf{h}^{\frac23})^{\frac32} = 1-\frac32 \tilde{a}_1  \mathsf{h}^{\frac23} + \frac38  \tilde{a}_1^2  \mathsf{h}^{\frac43} + \mathcal{O}\left( \mathsf{h}^2 \right),$$
one has
$$-\frac43 (\mathsf{h}^{-\frac23}-\tilde{a}_1)^{\frac32} = -\frac{4}{3\mathsf{h}} + 2\tilde{a}_1 \mathsf{h}^{-\frac13}+ \mathcal{O}\left(\mathsf{h}^{\frac13}\right),$$
hence
$$\ex^{-\frac43 (\mathsf{h}^{-\frac23}-\tilde{a}_1)^{\frac32}} = \ex^{ -\frac43 \mathsf{h}^{-1} + 2\tilde{a}_1 \mathsf{h}^{-\frac13}  } \left(1+\mathcal{O}\left( \mathsf{h}^{\frac13} \right) \right),$$
hence
$$\ex^{-\frac43 (\mathsf{c}-\tilde{a}_1)^{\frac32}}\left(1+ \mathcal{O}\left((\mathsf{c}-\tilde{a}_1)^{-\frac32} \right) \right) = \ex^{ -\frac43 \mathsf{h}^{-1} + 2\tilde{a}_1 \mathsf{h}^{-\frac13}  } \left(1+\mathcal{O}\left( \mathsf{h}^{\frac13} \right) \right),$$
which proves the second point.
\end{demo}
\vskip 3mm

\subsection{The upper edge of the first spectral band}\label{sec_E0max}

The upper edge of the first spectral band is the smallest value of $\mathbf{E}$ among the solutions of equations (\ref{eq_vp_Ha_0_b}) and (\ref{eq_vp_Ha_0_a}).

\vskip 3mm

\noindent We start by assuming that $\mathsf{c}\in [c_0,\tilde{c}_0)$. In this case (\ref{eq_vp_Ha_0_a}) has no solution with $-\mathbf{E} >0$ and we prove that (\ref{eq_vp_Ha_0_b}) has a unique solution such that  $-\mathbf{E} >0$ and $-\mathsf{c}-\mathbf{E}\in [-c_0,-\tilde{a}_1)$. Indeed, the function $z_0\ :\ [0,+\infty)\to [c_0,+\infty)$ is a continuous bijection and
$$\forall \mathsf{c} \geq c_0,\  \exists! x_0\geq 0,\ z_0(x_0)=x_0-\psi_0(x_0) = \mathsf{c} .$$
One sets $\breve{\mathbf{E}}_0=-x_0$ and since $x_0=\mathsf{c} +\psi_0(x_0) \in [-c_0 + \mathsf{c}, -\tilde{a}_1 + \mathsf{c})$ one has
\begin{equation}\label{eq_estim_Emin}
-\mathsf{c}+ \tilde{a}_1 < \breve{\mathbf{E}}_0 \leq -\mathsf{c} + c_0. 
\end{equation}
Thus, for $\mathsf{c}\in [c_0,\tilde{c}_0)$, $\mathbf{E}_{\mathrm{max}}^0=\breve{\mathbf{E}}_0$. 

\vskip 3mm
 
\noindent We then assume that $\mathsf{c} \geq \tilde{c}_0$. In this case, (\ref{eq_vp_Ha_0_b}) has still a unique solution such that  $-\mathbf{E} >0$ and $-\mathsf{c}-\mathbf{E}\in [-c_0,-\tilde{a}_1)$, namely  $\breve{\mathbf{E}}_0$, but one can also find a solution of (\ref{eq_vp_Ha_0_a}) with $-\mathbf{E} >0$. Indeed, the function from $[0,+\infty)\to [\tilde{c}_0,+\infty)$ which maps $x\geq 0$ to $x-\psi^0(x)$ is a continuous strictly increasing bijection. Thus, 
$$\forall \mathsf{c} \geq \tilde{c}_0,\  \exists! \tilde{x}_0\geq 0,\ \tilde{x}_0-\psi^0(\tilde{x}_0) = \mathsf{c} .$$
One sets $\tilde{\mathbf{E}}_0=-\tilde{x}_0$ and since $\tilde{x}_0=\mathsf{c} +\psi^0(\tilde{x}_0) \in (-a_1 + \mathsf{c}, -\tilde{c}_0 + \mathsf{c}]$ one has
$$-\mathsf{c}+ \tilde{c}_0 \leq \tilde{\mathbf{E}}_0 < -\mathsf{c} + a_1.$$   
Since $c_0 < \tilde{c}_0$, one has $\breve{\mathbf{E}}_0 <  \tilde{\mathbf{E}}_0 $ which implies  $\mathbf{E}_{\mathrm{max}}^0=\breve{\mathbf{E}}_0$ and $\mathbf{E}_{\mathrm{min}}^1=\tilde{\mathbf{E}}_0$. 

\vskip 3mm

\noindent As we identified $\mathbf{E}_{\mathrm{max}}^0$ among all the solutions of (\ref{eq_vp_Ha_0_b}) and (\ref{eq_vp_Ha_0_a}), we give more precise estimates. 

\begin{prop}\label{prop_upper_edge0}
We have the following estimates on $\mathbf{E}_{\mathrm{max}}^0$, in the semiclassical regime: 
\begin{enumerate}
\item For every $\mathsf{c} > c_0$, 
$$-\mathsf{c}+ \tilde{a}_1 <\mathbf{E}_{\mathrm{max}}^0 < -\mathsf{c} + c_0.$$
\item  For every $\mathsf{h}\in (0,c_0^{-\frac32})$,  
\begin{equation}\label{eq_prop_upper_edge0}
\mathbf{E}_{\mathrm{max}}^0 = -\mathsf{h}^{-\frac23}+\tilde{a}_1 +\alpha \sqrt{3} \mfrac{(u'(-\tilde{a}_1))^2  }{\tilde{a}_1 }\ex^{-\frac43 \mathsf{h}^{-1} + 2\tilde{a}_{1}\mathsf{h}^{-\frac13}} \left( 1 +  \mathcal{O} \left(\mathsf{h}^{\frac13} \right) \right).
\end{equation}
\end{enumerate}

\end{prop}

\begin{demo}
\noindent \textbf{(1)} The first point is (\ref{eq_estim_Emin}) and the inequality is strict since $z_0$ is strictly increasing.
\vskip 2mm

\noindent \textbf{(2)} For the second point, we follow the proof of point $(2)$ of Theorem \ref{thm_bottom}. We assume that $\mathsf{c} > c_0$ (for the identification of $\mathbf{E}_{\mathrm{max}}^0$) and thus $\mathsf{c} \geq \tilde{a}_1 + \tau$ with $\tau=c_0-\tilde{a}_1 >0$. One sets $X=-\mathbf{E}_{\mathrm{max}}^0 - \mathsf{c} + \tilde{a}_1$. Then $X$ satisfies
\begin{equation}\label{eq_demo_upper_1}
\frac{v'}{u'}(X-\tilde{a}_1) =  \frac{v}{u}(X + \mathsf{c}-\tilde{a}_1), \quad X\in [-c_0+\tilde{a}_1,0).
\end{equation}

\noindent Using \cite[10.4.59, 10.4.63]{AS} in the equality (\ref{eq_demo_bottom_1}), one gets

\begin{equation}\label{eq_demo_upper_4}
 \frac{v'}{u'}(X-\tilde{a}_1)= \alpha - \alpha \sqrt{3} \ex^{-\frac43 (X + \mathsf{c} - \tilde{a}_1 )^{\frac32}} \left(1+ \mathcal{O}\left(  \left(X + \mathsf{c} - \tilde{a}_1 \right)^{-\frac32} \right) \right).
\end{equation}

\noindent Again, we set $\epsilon = \ex^{-\frac43 (\mathsf{c}-\tilde{a}_1)^{\frac32}} \in (0,\ex^{-\frac43 (c_0-\tilde{a}_1)^{\frac32}}]$. We also define $Y$ as in the proof of point $(2)$ of Theorem \ref{thm_bottom}. Then, equality (\ref{eq_demo_bottom_8}) is still valid and the condition $X\in  [-c_0+\tilde{a}_1,0)$ implies that 
$$-c_0 + \tilde{a}_1 + \left( -\tfrac34 \ln(\epsilon)\right)^{\frac23} = -c_0 +\mathsf{c}   \leq X+\mathsf{c} -\tilde{a}_1 \leq \left( -\tfrac34 \ln(\epsilon)\right)^{\frac23}.$$
Moreover, for every $ \epsilon \in (0,\ex^{-\frac43 (c_0-\tilde{a}_1)^{\frac32}}]$,
$$\mfrac43  \frac{1}{|\ln (\epsilon)|}   \leq \left(X + \mathsf{c} - \tilde{a}_1 \right)^{-\frac32} \leq \left( \left( \mfrac34 \right)^{\frac23}  + \mfrac{\tilde{a}_1 - c_0}{|\ln (\epsilon)|^{\frac23}  }\right)^{-\frac32} \frac{1}{|\ln (\epsilon)|}.  $$
Using (\ref{eq_demo_upper_4}) and the relation $\alpha =  \frac{v'}{u'}(-\tilde{a}_1)$, one gets:
\begin{equation}
\begin{split}
\forall \epsilon & \in (0,\ex^{-\frac43 (c_0-\tilde{a}_1)^{\frac32}}],\  \frac{1}{\epsilon} \left(  \frac{v'}{u'}(\epsilon Y-\tilde{a}_1)-\frac{v'}{u'}(-\tilde{a}_1)\right) = \label{eq_demo_upper_5}\\
 &- \alpha \sqrt{3} \ex^{-  2  \left( -\frac34 \ln(\epsilon)\right)^{\frac13} \epsilon Y } \ex^{  \left( -\frac34 \ln(\epsilon)\right)^{-\frac13} (\epsilon Y)^2 Q(\epsilon Y)} \left( 1+  \mathcal{O}\left( \frac{1}{|\ln (\epsilon)|} \right) \right).
\end{split}
\end{equation}
Inequality (\ref{eq_demo_bottom_9}) is still valid and, using the fact that $\epsilon Y = X \in [-c_0-\tilde{a}_1,0)$, one gets that $Y$ is bounded and
$$X= \mathcal{O} \left(  \ex^{-\frac43 (\mathsf{c}-\tilde{a}_1)^{\frac32}} \right).$$
Then, the bootstrap argument gives the limit of $Y$ and 
$$\forall c>c_0,\ Y=-\alpha \sqrt{3} \mfrac{(u'(-\tilde{a}_1))^2  }{\tilde{a}_1} + \mathcal{O} \left( (\mathsf{c} - \tilde{a}_1)^{-\frac32} \right) $$
and for every $c> c_0$,
$$\mathbf{E}_{\mathrm{max}}^0 = -\mathsf{c}+\tilde{a}_1+ \alpha \sqrt{3}\mfrac{(u'(-\tilde{a}_1))^2  }{ \tilde{a}_1 } \ex^{-\frac43 (\mathsf{c}-\tilde{a}_1)^{\frac32}} +  \mathcal{O} \left( (\mathsf{c} - \tilde{a}_1)^{-\frac32} \ex^{-\frac43 (\mathsf{c}-\tilde{a}_1)^{\frac32}} \right),$$
which proves the second point with $\mathsf{c}=\mathsf{h}^{-\frac23}$ and the same computation as in the proof of the second point of Theorem \ref{thm_bottom}.
\end{demo}

\noindent We deduce from the estimates in $\mathsf{h}$ of $\mathbf{E}_{\mathrm{min}}^0$ and  $\mathbf{E}_{\mathrm{max}}^0$ the width of the first rescaled spectral band in the semiclassical regime.

\begin{prop}\label{prop_asym_width1}
For every $\mathsf{h}\in (0,c_0^{-\frac32})$, 
\begin{equation}\label{eq_prop_asym_width1}
 \mathbf{E}_{\mathrm{max}}^0 -\mathbf{E}_{\mathrm{min}}^0 =2 \alpha \sqrt{3}\mfrac{(u'(-\tilde{a}_1))^2  }{ \tilde{a}_1 } \ex^{-\frac43 \mathsf{h}^{-1} + 2\tilde{a}_{1}\mathsf{h}^{-\frac13}} \left( 1 +  \mathcal{O} \left(\mathsf{h}^{\frac13} \right) \right).
\end{equation}
\end{prop}

\vskip 5mm

\noindent For $\mathsf{c} \leq c_0$, the situation changes. The first rescaled spectral band recovers completely the range of the periodic potential $\mathbf{V}$ and is even larger. 
\vskip 2mm

\begin{prop}\label{prop_sb1_classic}
If $\mathsf{c}\in (0, c_0]$, $\mathbf{E}_{\mathrm{max}}^0\geq 0$ and we have
$$\left[\mathrm{min}\left(-\mfrac{\mathsf{c}}{2},-\mathsf{c} + \tilde{a}_1\right) ,0\right] \subset \left[\mathbf{E}_{\mathrm{min}}^0,\mathbf{E}_{\mathrm{max}}^0\right].$$
\end{prop}

\begin{demo}
\noindent If $\mathsf{c}\in (0,c_0]$, there is no longer a solution of (\ref{eq_vp_Ha_0_b}) or (\ref{eq_vp_Ha_0_a}) satisfying $-\mathbf{E} >0$. Thus, $\mathbf{E}_{\mathrm{max}}^0 \geq 0$. Using the upper bound on $\mathbf{E}_{\mathrm{min}}^0$ given in Theorem \ref{thm_bottom}, we have $-\mathsf{c} + \tilde{a}_1 \in [\mathbf{E}_{\mathrm{min}}^0,\mathbf{E}_{\mathrm{max}}^0]$. Using point $(1)$ of Theorem \ref{thm_bottom}, we also have $-\frac{\mathsf{c}}{2} \in [\mathbf{E}_{\mathrm{min}}^0,\mathbf{E}_{\mathrm{max}}^0]$, which proves the proposition.
\end{demo}

\vskip 3mm

\noindent Proposition \ref{prop_sb1_classic} along with Proposition \ref{prop_upper_edge0} imply Theorem \ref{thm_intro_transition}.  

\vskip 3mm 
\noindent The following proposition precise the behaviour of $\mathbf{E}_{\mathrm{max}}^0$.
\vskip 3mm

\begin{prop}\label{prop_E0max_0}
Let $\mathsf{c}\in (0, c_0]$. 
\begin{enumerate}
\item If $\mathsf{c} \in (\tilde{c}_1-\tilde{c}_0,c_0]$, then  $\mathbf{E}_{\mathrm{max}}^0 \in (-\mathsf{c} + \tilde{c}_{0}, 0]$.
\item If $\mathsf{c} <\tilde{c}_1-\tilde{c}_0$, let $p_0$ defined in (\ref{def_p_0}). Then, $\mathbf{E}_{\mathrm{max}}^0 \in [-\mathsf{c} + \tilde{c}_{p_0}, \tilde{c}_{p_0+1}]$ or $\mathbf{E}_{\mathrm{max}}^0 \in [-\mathsf{c} + \tilde{c}_{p_0-1}, \tilde{c}_{p_0}]$.
\item \begin{equation}\label{eq_prop_E0_max_0}
 \lim_{\mathsf{c} \to 0} \mathbf{E}_{\mathrm{max}}^0 = +\infty.
\end{equation}
\end{enumerate}
\end{prop}

\begin{demo}
\textbf{(2)} Since $(\tilde{c}_{p+1}-\tilde{c}_{p})_{p\geq 0}$ is strictly decreasing and converges to $0$, for any $\mathsf{c} \in [0,\tilde{c}_1-\tilde{c}_0)$, the integer $p_0\geq 1$ defined in (\ref{def_p_0}) is well defined and unique.

\noindent $\mathbf{E}_{\mathrm{max}}^0$ is a solution of either (\ref{eq_vp_Ha_0_a}) or (\ref{eq_vp_Ha_0_b}). Let $k\geq 1$. The restriction of the function $\frac{v}{u}$ to $(-\tilde{c}_{2k},-\tilde{c}_{2k-2} )$ is a strictly increasing and continuous bijection from $(-\tilde{c}_{2k},-\tilde{c}_{2k-2} )$ to $\R$ denoted by $\left(\frac{v}{u} \right)_{2k-2}$, and  $\frac{v'}{u'}$ induce a strictly increasing and continuous bijection from $(-\tilde{c}_{2k+1},-\tilde{c}_{2k-1})$ to $\R$ denoted by $\left(\frac{v'}{u'} \right)_{2k-1}$. Then, studying the sign of $f_x$ for $x\in (-\tilde{c}_{2k},-\tilde{c}_{2k-2} )$ and using the Sturm-Picone's lemma in a way similar as in the proof of Lemma \ref{lem_Eomax2}, we prove that $x\mapsto x-\left(\frac{v'}{u'} \right)_{2k-1}^{-1}\left(\frac{v}{u} \right)_{2k-2}(x)$ is strictly increasing and continuous from $(-\tilde{c}_{2k+1},-\tilde{c}_{2k-1})$ to $(\tilde{c}_{2k+1}-\tilde{c}_{2k}, \tilde{c}_{2k-1}-\tilde{c}_{2k-2})$.
\vskip 1mm

\noindent Thus, (\ref{eq_vp_Ha_0_b}) admits a unique solution $\underline{\mathbf{E}}_{k}$ with $-\underline{\mathbf{E}}_{k}\in (-\tilde{c}_{2k}, - \tilde{c}_{2k-2})$ and $-\mathsf{c}- \underline{\mathbf{E}}_{k}\in (-\tilde{c}_{2k+1}, - \tilde{c}_{2k-1})$. 
\vskip 2mm

\noindent To study the sign of $f_x$, we need to know the signs of $u$, $v$, $u'$ and $v'$ on the interval $(-\tilde{c}_{2k+1}, -\tilde{c}_{2k-2})$, since  $x\in (-\tilde{c}_{2k},-\tilde{c}_{2k-2} )$ and $x-z\in (-\tilde{c}_{2k+1},-\tilde{c}_{2k-1} )$. For example, we have on $(-\tilde{c}_{2k+1},-c_{2k+1})$,
$$(-1)^k u <0,\quad (-1)^k u' >0,\quad (-1)^k v >0,\quad (-1)^k v' <0$$ 
and on $(-c_{2k+1},-\tilde{c}_{2k})$,
$$(-1)^k u <0,\quad (-1)^k u' >0,\quad (-1)^k v <0,\quad (-1)^k v' <0$$ 
and the signs alternate on the successives intervals $(-\tilde{c}_{2k},-c_{2k} )$, $(-c_{2k}, -\tilde{c}_{2k-1})$,  $(-\tilde{c}_{2k-1}, -c_{2k-1})$ and $(-c_{2k-1}, -\tilde{c}_{2k-2})$.

\noindent The restriction of the function $\frac{v}{u}$  to $(-\tilde{c}_{2k+2},-\tilde{c}_{2k})$ is a strictly increasing and continuous bijection from $(-\tilde{c}_{2k+2},-\tilde{c}_{2k})$ to $\R$ denoted by $\left(\frac{v}{u} \right)_{2k}$. We set for every $x\in (-\tilde{c}_{2k+1},-\tilde{c}_{2k-1})$ and every $z\in (-\tilde{c}_{2k+1}+\tilde{c}_{2k+2},-\tilde{c}_{2k-1}+\tilde{c}_{2k} )$,  $\bar{g}_x(z)=v(x-z)u'(x)-u(x-z)v'(x)$. Then $\bar{g}_x$ satisfies the Airy equation and using the signs of  $u$, $v$, $u'$ and $v'$ given above and a Sturm-Picone's argument as in Lemma \ref{lem_Eomax2}, we prove that $x\mapsto x-\left(\frac{v}{u} \right)_{2k}^{-1}\left(\frac{v'}{u'} \right)_{2k-1}(x)$ is strictly increasing and continuous from  $(-\tilde{c}_{2k+1},-\tilde{c}_{2k-1})$ to $(-\tilde{c}_{2k+1}+\tilde{c}_{2k+2},-\tilde{c}_{2k-1}+\tilde{c}_{2k} )$. 
\vskip 1mm

\noindent Thus, (\ref{eq_vp_Ha_0_a}) admits a unique solution $\bar{\mathbf{E}}_{k}$ with $-\bar{\mathbf{E}}_{k}\in (-\tilde{c}_{2k+1}, - \tilde{c}_{2k-1})$ and $-\mathsf{c}- \bar{\mathbf{E}}_{k}\in (-\tilde{c}_{2k+2}, - \tilde{c}_{2k})$. 
\vskip 3mm

\noindent Since $\tilde{c}_{p_0+1} - \tilde{c}_{p_0} < \mathsf{c} < \tilde{c}_{p_0} - \tilde{c}_{p_0-1}$, we have either $\mathbf{E}_{\mathrm{max}}^0 = \underline{\mathbf{E}}_{k}$ or $\mathbf{E}_{\mathrm{max}}^0 = \bar{\mathbf{E}}_{k}$ for $k$ equal to the integer part of $\frac{p_0}{2}$.
\vskip 1mm

\noindent We deduce that $-\mathbf{E}_{\mathrm{max}}^0\in [-\tilde{c}_{p_0+1},-\tilde{c}_{p_0-1}]$ and $-\mathsf{c}-\mathbf{E}_{\mathrm{max}}^0\in [-\tilde{c}_{p_0+2},-\tilde{c}_{p_0}]$, or  $-\mathbf{E}_{\mathrm{max}}^0\in [-\tilde{c}_{p_0},-\tilde{c}_{p_0-2}]$ and $-\mathsf{c}-\mathbf{E}_{\mathrm{max}}^0\in [-\tilde{c}_{p_0+1},-\tilde{c}_{p_0-1}]$. This proves the second point. 
\vskip 3mm

\noindent \textbf{(1)} If $\mathsf{c} \in (\tilde{c}_1-\tilde{c}_0,c_0]$, then  $\mathsf{c} \in (\tilde{c}_2-\tilde{c}_1,\tilde{c}_0)$. The function  $\frac{v'}{u'}$ induces a strictly increasing and continuous bijection from $(-\tilde{c}_{1},0)$ to $\R$ denoted by $\left(\frac{v'}{u'} \right)_{0}$. Then, with the previous notations, the function  $x\mapsto x-\left(\frac{v}{u} \right)_{0}^{-1}\left(\frac{v'}{u'} \right)_{0}(x)$ is strictly increasing and continuous from  $(-\tilde{c}_{1},0)$ to $(\tilde{c}_{2}-\tilde{c}_{1},\tilde{c}_{0})$, using again Sturm-Picone's Lemma with $\bar{g}_x$ for $x\in (-\tilde{c}_{1},0)$. Thus,  (\ref{eq_vp_Ha_0_a}) admits a unique solution $\bar{\mathbf{E}}_{0}$ with $- \bar{\mathbf{E}}_{0}\in (-\tilde{c}_{1},0)$ and $-\mathsf{c}- \bar{\mathbf{E}}_{0}\in (-\tilde{c}_{2}, - \tilde{c}_{0})$. Since $\mathbf{E}_{\mathrm{max}}^0 = \bar{\mathbf{E}}_{0}$, we proved the first point.  

\vskip 3mm

\noindent \textbf{(3)} The integer $p_0$ tends to $+\infty$ when $\mathsf{c}$ tends to $0$. Indeed, using (\ref{eq_approx_tc2j}) and (\ref{eq_approx_tc2j1}), one has
$$p_0= \mathcal{O}\left( \mathsf{c}^{-3} \right).$$
Since $p_0\xrightarrow[\mathsf{c} \to 0]{}+\infty$ and $c_k\xrightarrow[k \to +\infty]{}+\infty$, we get (\ref{eq_prop_E0_max_0}) and prove the third point.
\end{demo}

\begin{demo}(of Theorem \ref{thm_intro_quantum_1st_band})

\noindent \textbf{(2)} We start with the proof of the second point since it will imply the first limit in the first point of Theorem \ref{thm_intro_quantum_1st_band}. We look at the behaviour of $\mathbf{E}_{\mathrm{min}}^0$ when $\mathsf{c}$ tends to $0$, since it is easier than considering the limit when $\mathsf{h}$ tends to infinity. Since $\mathbf{E}_{\mathrm{min}}^0$ is a solution of (\ref{eq_vp_Ha_d_a}) with $-\mathbf{E}_{\mathrm{min}}^0>0$ and $-\mathsf{c}-\mathbf{E}_{\mathrm{min}}^0 \in (-\tilde{a}_1,0)$, when $\mathsf{c}$ tends to $0$,  both $-\mathbf{E}_{\mathrm{min}}^0$ and $-\mathsf{c} - \mathbf{E}_{\mathrm{min}}^0$ tends to $0$. In order to avoid the technical difficulty induced by the fact that $\frac{v'}{u'}$ tends to $+\infty$ at $0$, we use the fact that $\mathbf{E}_{\mathrm{min}}^0$ is also the unique solution in $(-\mathsf{c},0)$ of the equation 
\begin{equation}\label{eq_demo_bottom_0_1}
\frac{u'}{v'}(-\mathsf{c}-\mathbf{E})= \frac{u'}{v'}(-\mathbf{E}).
\end{equation}

\noindent Note that $u$ and $v$ stand for $f$ and $g$ in \cite[10.4.3]{AS}. Thus, for $x$ in a neighborhood of $0$ where $v'$ does not vanish, 
$$\left( \frac{u'}{v'} \right)'(x) = \frac{x}{(v'(x))^2}\quad \mbox{ and } \quad v(x)=x+\frac{x^4}{12} + \mathcal{O}(x^7).$$
One deduces
$$\left( \frac{u'}{v'} \right)'(x) = x-\frac{2}{3} x^4 +  \mathcal{O}(x^7)$$
hence
\begin{equation}\label{eq_demo_bottom_0_2}
\left( \frac{u'}{v'} \right)'(x) =\left( \frac{u'}{v'} \right)'(0)+ \frac12 x^2-\frac{2}{15} x^5 +  \mathcal{O}(x^8).
\end{equation}
Let $y=-\mathbf{E}_{\mathrm{min}}^0 -\frac{\mathsf{c}}{2}$. Then, (\ref{eq_demo_bottom_0_1}) rewrites
\begin{equation}\label{eq_demo_bottom_0_3}
\frac{u'}{v'}\left(y-\mfrac{\mathsf{c}}{2}\right)= \frac{u'}{v'}\left(y+\mfrac{\mathsf{c}}{2}\right).
\end{equation}
Since $-\mathsf{c} < \mathbf{E}_{\mathrm{min}}^0 < 0$, 
\begin{equation}\label{eq_demo_bottom_0_5}
 |y|\leq \frac{\mathsf{c}}{2}.
\end{equation}
By (\ref{eq_demo_bottom_0_5}), $y= \mathcal{O}(\mathsf{c})$. Thus equation (\ref{eq_demo_bottom_0_3}) and equality (\ref{eq_demo_bottom_0_2}) imply
$$\frac12 \left( y^2 + \left( \mfrac{\mathsf{c}}{2} \right)^2 - \mathsf{c} y \right) - \frac{2}{15} \left( y -\mfrac{\mathsf{c}}{2}  \right)^5 = \hspace{5cm}$$
$$\hspace{2cm}  \frac12 \left( y^2 + \left( \mfrac{\mathsf{c}}{2} \right)^2 + \mathsf{c} y \right) - \frac{2}{15} \left( y +\mfrac{\mathsf{c}}{2}  \right)^5 + \mathcal{O}(\mathsf{c}^8)$$
that is
$$- \mathsf{c} y - \frac{2}{15} \left( -5\mathsf{c} y^4 - 20 \left( \mfrac{\mathsf{c}}{2} \right)^3 y^2 - 2 \left( \mfrac{\mathsf{c}}{2} \right)^5  \right) + \mathcal{O}(\mathsf{c}^8)= 0,$$
which implies 
\begin{equation}\label{eq_demo_bottom_0_4}
y -  \frac{2}{15}y^4 - \frac43  \left( \mfrac{\mathsf{c}}{2} \right)^3 y^2 = \frac{2}{15}\left( \mfrac{\mathsf{c}}{2} \right)^4 + \mathcal{O}(\mathsf{c}^7).
\end{equation}
Then, (\ref{eq_demo_bottom_0_5}) and (\ref{eq_demo_bottom_0_4}) give $y=\mathcal{O}(\mathsf{c}^4)$.  We write $y=\mathsf{c}^4 \rho$ where $\rho$ is a bounded function of $\mathsf{c}$. Then, (\ref{eq_demo_bottom_0_4})  rewrites
$$\rho- \frac{2}{15} \mathsf{c}^{12} \rho^4 = \frac{1}{120} + \mathcal{O}(\mathsf{c}^3),$$
hence
$$\rho= \frac{1}{120} + \mathcal{O}(\mathsf{c}^3) $$
and then, 
$$y= \frac{1}{120}\mathsf{c}^4 + \mathcal{O}(\mathsf{c}^7), $$
which proves (\ref{eq_thm_bottom3}) after using $\mathsf{c}=\mathsf{h}^{-\frac23}$ and multiplying by $\mathsf{h}^{\frac23}$. 
\vskip 3mm

\noindent \textbf{(1)} The first limit is a direct consequence of point $(2)$. The second limit follows from point $(3)$ of Proposition \ref{prop_E0max_0} and $\mathsf{c}=\mathsf{h}^{-\frac23}$.
\end{demo}

\bigskip

\subsection{The first spectral gap}

In the discussion before Proposition \ref{prop_upper_edge0}, we  identified both $\mathbf{E}_{\mathrm{max}}^0$ and $\mathbf{E}_{\mathrm{min}}^1$ in the case $\mathsf{c} \geq \tilde{c}_0$. We had obtained, for every $\mathsf{c} \geq \tilde{c}_0$, 
$$-\mathsf{c}+ \tilde{a}_1 < \mathbf{E}_{\mathrm{max}}^0  \leq -\mathsf{c} + c_0$$
and 
$$-\mathsf{c}+ \tilde{c}_0 \leq \mathbf{E}_{\mathrm{min}}^1 < -\mathsf{c} +a_1.$$
This yields a first estimate of the first spectral gap:
$$0< \tilde{c}_0-c_0 \leq \mathbf{E}_{\mathrm{min}}^1 -  \mathbf{E}_{\mathrm{max}}^0 < a_1-\tilde{a}_1.$$
In particular, the first gap is always open. 
\vskip 3mm

\noindent Similarly to the estimates for the edges of the first spectral band, we prove the following estimate for $\mathbf{E}_{\mathrm{min}}^1$ in the semiclassical regime.

\begin{prop}\label{prop_asym_E1min}
Let $\tau >0$ be arbitrary small. For every $\mathsf{h}\in (0,(a_1+\tau)^{-\frac32}]$,  
\begin{equation}\label{eq_prop_asym_E1min}
\mathbf{E}_{\mathrm{min}}^1 = -\mathsf{h}^{-\frac23}+a_1 -\alpha \sqrt{3} (u(-a_1))^2  \ex^{-\frac43 \mathsf{h}^{-1} + 2a_{1}\mathsf{h}^{-\frac13}} \left( 1 +  \mathcal{O} \left(\mathsf{h}^{\frac13} \right) \right).
\end{equation}
 \end{prop}

\begin{demo}
We follow the proof of point $(2)$ of Theorem \ref{thm_bottom}. Let $\tau >0$ be arbitrary small. We assume that $\mathsf{c} \geq a_1+\tau$. One sets $X=-\mathbf{E}_{\mathrm{min}}^1 - \mathsf{c} + a_1$ which satisfies
\begin{equation}\label{eq_demo_Emin1_1}
\frac{v}{u}(X-a_1) =  \frac{v'}{u'}(X + \mathsf{c}-a_1), \quad X\in [0,-\tilde{c}_0 + a_1).
\end{equation}
Using \cite[10.4.61, 10.4.66]{AS} in the equality (\ref{eq_demo_bottom_1}), one gets
\begin{equation}\label{eq_demo_Emin1_2}
 \frac{v}{u}(X-a_1)= \alpha + \alpha \sqrt{3} \ex^{-\frac43 (X + \mathsf{c} - a_1 )^{\frac32}} \left(1+ \mathcal{O}\left(  \left(X + \mathsf{c} - a_1 \right)^{-\frac32} \right) \right).
\end{equation}

\noindent We set $\epsilon = \ex^{-\frac43 (\mathsf{c}-a_1)^{\frac32}}\in (0,\ex^{-\frac43 \tau^{\frac32}}]$ and $Y=\frac{1}{\epsilon} X$. Since $\alpha = \frac{v}{u}(-a_1) $ and using the fact that $\left(\frac{v}{u} \right)' = \frac{1}{u^2}$ is bounded from below by $\frac{\alpha}{2}$ by (\ref{eq_lem_uv_alpha}), one shows with a similar proof as (\ref{eq_demo_bottom_9}) that: 
$$\forall \epsilon \in (0,\ex^{-\frac43 \tau^{\frac32}}],\ \left| \frac{v}{u}(X-a_1) - \frac{v}{u}(-a_1) \right| \geq \frac{\alpha}{2} \epsilon Y.$$
Then, following the proof of point $(2)$ of Theorem \ref{thm_bottom} and using that $\epsilon Y = X \in [0, -\tilde{c}_0+a_1]$ is bounded, one gets that $Y$ is bounded and 
$$X= \mathcal{O} \left(  \ex^{-\frac43 (\mathsf{c}-a_1)^{\frac32}} \right).$$
Since $Y$ is bounded one has:
$$\forall \epsilon \in (0,\ex^{-\frac43 \tau^{\frac32}}],\ \frac{v}{u}(\epsilon Y-a_1) - \frac{v}{u}(-a_1) = \left(\frac{v}{u} \right)' (-a_1) (\epsilon Y) + \mathcal{O}(\epsilon^2).$$
But, $ \left(\frac{v}{u} \right)' (-a_1) = \frac{1}{(u(-a_1))^2}$ and we get, similarly to (\ref{eq_demo_bottom_10}), 
$$\forall c\geq a_1+\tau,\ Y=\alpha \sqrt{3} (u(-a_1))^2 + \mathcal{O} \left( (\mathsf{c} - a_1)^{-\frac32} \right)$$
from which we obtain (\ref{eq_prop_asym_E1min}) by taking $\mathsf{c}=\mathsf{h}^{-\frac23}$ as before.
\end{demo}

\noindent Combining the estimates of $\mathbf{E}_{\mathrm{max}}^0$ and $\mathbf{E}_{\mathrm{min}}^1$ we deduce an estimate of the rescaled first gap in the semiclassical regime.

\begin{prop}\label{prop_asymp_gap1}
Let $\tau >0$ be arbitrary small. For every $\mathsf{h}\in (0,(a_1+\tau)^{-\frac32}]$,
\begin{equation} \label{eq_prop_asym_gap1}
\mathbf{E}_{\mathrm{min}}^1- \mathbf{E}_{\mathrm{max}}^0  =  a_1-\tilde{a}_1 - \alpha \sqrt{3}  (u(-a_1))^2 \ex^{-\frac43 \mathsf{h}^{-1} + 2a_{1}\mathsf{h}^{-\frac13}} \left( 1 +  \mathcal{O} \left(\mathsf{h}^{\frac13} \right) \right).
\end{equation}
\end{prop}

\begin{demo}
We combine (\ref{eq_prop_asym_E1min}) and point $(2)$ of Proposition \ref{prop_upper_edge0} and use $a_1 > \tilde{a}_1$.
\end{demo}

\section{Counting the spectral bands in the range of $V$}\label{sec_pbands}

In this Section, we prove Theorem \ref{thm_intro_pbands} by determining the band edges which are contained in the interval $[-\mathsf{c},0]$ for a fixed $\mathsf{c}$.

\vskip 3mm
\begin{prop}\label{prop_Eomax1}
Let $p\geq 0$ and assume that $\mathsf{c} \geq \tilde{c}_p$. Then, for every $k\in \{ 0,\ldots, p\}$, 
\begin{enumerate}
 \item If $k=2j$ is even, (\ref{eq_vp_Ha_0_a}) has a unique solution $\hat{\mathbf{E}}_{2j}$ with $-\hat{\mathbf{E}}_{2j} \in [0,+\infty)$, $-\mathsf{c}-\hat{\mathbf{E}}_{2j}\in (-a_{j+1},-\tilde{c}_{2j}]$ and satisfying:
\begin{equation}\label{eq_prop_Eomax1_even}
-\mathsf{c} +\tilde{c}_{2j} \leq  \hat{\mathbf{E}}_{2j} <  -\mathsf{c} + a_{j+1}.
\end{equation}
\item  If $k=2j+1$ is odd, (\ref{eq_vp_Ha_d_a}) has a unique solution $\hat{\mathbf{E}}_{2j+1}$ with $-\hat{\mathbf{E}}_{2j+1} \in [0,+\infty)$, $-\mathsf{c}-\hat{\mathbf{E}}_{2j+1}\in (-\tilde{a}_{j+2},-\tilde{c}_{2j+1}]$ and satisfying:
\begin{equation}\label{eq_prop_Eomax1_odd}
 -\mathsf{c} + \tilde{c}_{2j+1} \leq  \hat{\mathbf{E}}_{2j+1} <  -\mathsf{c}+ \tilde{a}_{j+2}.
\end{equation}    
\end{enumerate}
\end{prop}

\begin{demo} By Lemma \ref{lem_Eomax1}, for every $k\geq 0$, the function $x\mapsto x-\psi^k(x)$ is a strictly increasing and continuous bijection from $[0,+\infty)$ to $[\tilde{c}_k,+\infty)$. Thus, if $\mathsf{c} \geq \tilde{c}_p\geq \tilde{c}_k$, there exists a unique $x^k\geq 0$ such that $\mathsf{c} = x^k - \psi^k(x^k)$. Let $\hat{\mathbf{E}}_k$ be such that $\hat{\mathbf{E}}_k =- x^k$. Then, if $k=2j$,  $\hat{\mathbf{E}}_{2j}$ is the unique solution of (\ref{eq_vp_Ha_0_a}) with $-\hat{\mathbf{E}}_{2j} \in [0,+\infty)$ and $-\mathsf{c}-\hat{\mathbf{E}}_{2j} \in (-a_{j+1},-\tilde{c}_{2j}]$. Moreover,
$$-a_{j+1} < -\hat{\mathbf{E}}_k - \mathsf{c} \leq -\tilde{c}_k  <0 \leq - \hat{\mathbf{E}}_k, $$
and we get (\ref{eq_prop_Eomax1_even}). If $k=2j+1$,  $\hat{\mathbf{E}}_{2j+1}$ is the unique solution of (\ref{eq_vp_Ha_d_a}) with $-\hat{\mathbf{E}}_{2j+1}\in [0,+\infty)$ and $-\mathsf{c}-\hat{\mathbf{E}}_{2j+1}\in (-\tilde{a}_{j+2},-\tilde{c}_{2j+1}]$. Moreover,
$$-\tilde{a}_{j+2} < - \hat{\mathbf{E}}_k - \mathsf{c} \leq -\tilde{c}_k  <0 \leq -  \hat{\mathbf{E}}_k, $$
and we get (\ref{eq_prop_Eomax1_odd}).
\end{demo}
\vskip 3mm

\vskip 5mm

\begin{prop}\label{prop_Eomax2}
Assume that $\mathsf{c} \geq c_p$. Then, for every $k\in \{0,\ldots ,p\}$,
\begin{enumerate}
 \item If $k=2j$ is even, (\ref{eq_vp_Ha_0_b}) has a unique solution $\breve{\mathbf{E}}_{2j}$ with $-\breve{\mathbf{E}}_{2j} \in [0,+\infty)$, $-\mathsf{c}-\breve{\mathbf{E}}_{2j}\in [c_{2j},-\tilde{a}_{j+1})$ and satisfying:
\begin{equation}\label{eq_prop_Eomax2_even}
-\mathsf{c} + \tilde{a}_{j+1} < \breve{\mathbf{E}}_{2j} \leq  -\mathsf{c} + c_{2j}.
\end{equation}
\item  If $k=2j+1$ is odd, (\ref{eq_vp_Ha_d_b}) has a unique solution $\breve{\mathbf{E}}_{2j+1}$ with $-\breve{\mathbf{E}}_{2j+1} \in [0,+\infty)$, $-\mathsf{c}-\breve{\mathbf{E}}_{2j+1}\in [c_{2j+1},-a_{j+1})$ and satisfying:
\begin{equation}\label{eq_prop_Eomax2_odd}
 -\mathsf{c} + a_{j+1} < \breve{\mathbf{E}}_{2j+1} \leq -\mathsf{c} + c_{2j+1}.
\end{equation}    
\end{enumerate}
\end{prop}

\begin{demo} Let $k\in \{0,\ldots ,p\}$. Since $\mathsf{c} \geq c_p$, we have $\mathsf{c} \in [c_k,+\infty)$. Thanks to Lemma \ref{lem_Eomax2}, $z_k$ is continuous and strictly increasing and there exists a unique real number  $x_k\geq 0$ such that $\mathsf{c}=z_k(x_k)$. Let $\breve{\mathbf{E}}_k$ be such that $\breve{\mathbf{E}}_k = -x_k$. Then, $\breve{\mathbf{E}}_{2j}$ is the unique solution of (\ref{eq_vp_Ha_0_b}) such that $-\breve{\mathbf{E}}_{2j} \in [0,+\infty)$ and $-\mathsf{c}-\breve{\mathbf{E}}_{2j} \in  [c_{2j},-\tilde{a}_{j+1})$. Moreover,
$$-c_k \leq - \breve{\mathbf{E}}_k - \mathsf{c} < -\tilde{a}_{j+1}$$
and we get (\ref{eq_prop_Eomax2_even}).

\noindent Similarly,  $\breve{\mathbf{E}}_{2j+1}$ is the unique solution of (\ref{eq_vp_Ha_d_b}) such that $-\breve{\mathbf{E}}_{2j+1} \in [0,+\infty)$ and $-\mathsf{c}-\breve{\mathbf{E}}_{2j+1} \in [c_{2j+1},-a_{j+1})$. Moreover,
$$-c_k \leq - \breve{\mathbf{E}}_k - \mathsf{c} < -a_{j+1}$$
and we get (\ref{eq_prop_Eomax2_odd}).
\end{demo}
\vskip 2mm

\noindent We deduce from Proposition \ref{prop_Eomax1} and Proposition \ref{prop_Eomax2} the following proposition on the $p$ first spectral bands and the $p-1$ first spectral gaps of the operator $H$.

\begin{prop}\label{prop_Ekmin_Ekmax}
Let $p\geq 0$. Assume that $\mathsf{c} \geq \tilde{c}_p$. 
\begin{enumerate}
 \item For every $k\in \{0,\ldots ,p \}$, $\mathbf{E}_{\mathrm{min}}^{k+1} = \hat{\mathbf{E}}_k$ and  $\mathbf{E}_{\mathrm{max}}^k = \breve{\mathbf{E}}_k$.
 \item We have the estimates on the spectral gaps:
$$\forall j \geq 1,\ 0 < \tilde{c}_{2j-1}-c_{2j-1} \leq \mathbf{E}_{\mathrm{min}}^{2j}- \mathbf{E}_{\mathrm{max}}^{2j-1} \leq \tilde{a}_{j+1}-a_j $$
and 
$$\forall j \geq 0,\ 0 < \tilde{c}_{2j}-c_{2j} \leq \mathbf{E}_{\mathrm{min}}^{2j+1}- \mathbf{E}_{\mathrm{max}}^{2j} \leq a_{j+1}-\tilde{a}_{j+1}.$$
In particular, all the spectral gaps are open.
\end{enumerate}
\end{prop}

\begin{demo}
\noindent \textbf{(1)} Using the estimates obtained on $\hat{\mathbf{E}}_k$ and  $\breve{\mathbf{E}}_k$ and using the fact that $c_k < \tilde{c}_k$, we have $\breve{\mathbf{E}}_k < \hat{\mathbf{E}}_k$. Since 
$$-\mathsf{c} < \mathbf{E}_{\mathrm{min}}^0 < -\mathsf{c} + \tilde{a}_1 < \breve{\mathbf{E}}_0 <  \hat{\mathbf{E}}_0, $$
we have $\mathbf{E}_{\mathrm{max}}^0 =  \breve{\mathbf{E}}_0 $ and $\mathbf{E}_{\mathrm{min}}^{1} = \hat{\mathbf{E}}_0$. Then using $\breve{\mathbf{E}}_k < \hat{\mathbf{E}}_k$, we deduce the first point.

\noindent \textbf{(2)} These two estimates are deduced directly from the estimates proven on $\hat{\mathbf{E}}_k$ and  $\breve{\mathbf{E}}_k$ in Proposition \ref{prop_Eomax1} and Proposition \ref{prop_Eomax2}. We just have to be careful with the fact that $\mathbf{E}_{\mathrm{min}}^{2j}=\hat{\mathbf{E}}_{2j-1}$ and $\mathbf{E}_{\mathrm{min}}^{2j+1}=\hat{\mathbf{E}}_{2j}$ and to use the right estimate in Proposition \ref{prop_Eomax1} depending on the parity of $k$. 
\end{demo}

\noindent Propositions \ref{prop_Eomax1}, \ref{prop_Eomax2} and \ref{prop_Ekmin_Ekmax} imply the proof of Proposition \ref{thm_zero_edges}.
\vskip 3mm

\begin{demo} (of Proposition \ref{thm_zero_edges}).
For every $p\geq 0$, $-\mathbf{E}_{\mathrm{max}}^p (\mathsf{c}) = (z_p)^{-1}(\mathsf{c})$ and since $z_p$ is strictly increasing and continuous on $[0,+\infty)$, $\mathsf{c} \mapsto \mathbf{E}_{\mathrm{max}}^p (\mathsf{c}) $ is strictly decreasing and continuous on $[0,+\infty)$. Since for every $p\geq 0$, $\mathbf{E}_{\mathrm{max}}^p (c_p) = 0$, $c_p$ is the unique zero in $[0,+\infty)$ of the function $\mathsf{c} \mapsto \mathbf{E}_{\mathrm{max}}^p (\mathsf{c}) $. 
\vskip 2mm

\noindent Since for every $p\geq 0$, $-\mathbf{E}_{\mathrm{min}}^{p+1} (\mathsf{c}) = (z^p)^{-1}(\mathsf{c})$ (where $z^p\ :\ x \mapsto x-\psi^p(x)$ is strictly increasing and continuous on $[0,+\infty)$),   $\mathsf{c} \mapsto \mathbf{E}_{\mathrm{min}}^{p+1} (\mathsf{c}) $ is also strictly decreasing and continuous on $[0,+\infty)$, and since for every $p\geq 0$, $\mathbf{E}_{\mathrm{min}}^{p+1} (\tilde{c}_p) = 0$, $\tilde{c}_p$ is the unique zero in $[0,+\infty)$ of the function $\mathsf{c} \mapsto \mathbf{E}_{\mathrm{min}}^{p+1} (\mathsf{c}) $.
\end{demo}

\vskip 3mm

\noindent The estimates in Propositions \ref{prop_Eomax1}, \ref{prop_Eomax2} and \ref{prop_Ekmin_Ekmax} combined with the intervals given in Lemma \ref{lem_asymp_ck} lead to the proof of Theorem \ref{thm_intro_pbands}. Before that, we prove a technical lemma.

\begin{lem}\label{lem_a_b}
For every $y\in \R_+^{*}$, let $I(y)= (\frac32)^{\frac13}\frac{y^{\frac32}+1}{y^2+y+1}$. Then, for every $\eta >0$ and every real numbers $0< b<a$ such that  $\frac{a-b}{b}\in [0, \eta]$,
\begin{equation}\label{eq_lem_a_b}
(a-b)b^{-\frac13}I((1+\eta)^{\frac23})\leq  \left(\tfrac32 a \right)^{\frac23}-\left(\tfrac32 b\right)^{\frac23}\leq (a-b)b^{-\frac13}I(1). 
\end{equation}
\end{lem}

\begin{demo}
One checks that $I(1)=(\frac32)^{\frac13}$, $I'(1)=-\frac16$ and 
$$\forall y>0,\ I'(y)=-\frac{\frac32+(t-1)(\frac12+2t+\frac12(1+t)t^3)}{(1+t^2+t^4)^2},\quad \mbox{where } t^2=y.$$
Hence $I'(y)<0$ for $y>1$ and $I$ is strictly decreasing on $[1, 1+\eta]$ for all $\eta>0$.

\noindent In particular, for every real numbers $0< b<a$ such that  $\frac{a-b}{b}\in [0, \eta]$, $\left(\frac{a}{b} \right)^{\frac23} \in [1,(1+\eta)^{\frac23}] \subset [1,1+\eta]$ and  $$I((1+\eta)^{\frac23}) \leq I\left( \left(\tfrac{a}{b} \right)^{\frac23} \right) \leq I(1).$$
Since 
$$I\left( \left(\tfrac{a}{b} \right)^{\frac23} \right)=(a-b)^{-1}\left( \left(\tfrac32 a \right)^{\frac23}-\left(\tfrac32 b\right)^{\frac23} \right)b^{\frac13}, $$
we get (\ref{eq_lem_a_b}). 

\end{demo}

\vskip 3mm

\begin{demo}
(of Theorem \ref{thm_intro_pbands}). Let $\mathsf{c} > c_0$. The first point in Theorem \ref{thm_intro_pbands} is a direct consequence of point $(1)$ of Proposition \ref{prop_Ekmin_Ekmax} and Propositions \ref{prop_Eomax1} and \ref{prop_Eomax2} which ensure that for every $k\in \{0,\ldots , k_0\}$, $\hat{\mathbf{E}}_k$ and  $\breve{\mathbf{E}}_k$ are in $[-\mathsf{c},0]$.

\noindent For the second point, using Propositions \ref{prop_Eomax1}, \ref{prop_Eomax2} and \ref{prop_Ekmin_Ekmax} one deduces that 
$$\forall p\in \{2,\ldots, k_0\},\  \tilde{c}_p-c_p \leq \mathbf{E}_{\mathrm{min}}^{p+1}- \mathbf{E}_{\mathrm{max}}^{p} \leq \tilde{c}_p-\tilde{c}_{p-2}$$
and
$$\forall p\in \{2,\ldots, k_0\},\ 0 < \mathbf{E}_{\mathrm{max}}^{p}- \mathbf{E}_{\mathrm{min}}^{p} \leq c_p-\tilde{c}_{p-1}.$$
Let $p\in \{2,\ldots, k_0\}$. Assume that $p$ is even, that is $p=2j$ for $j\geq 1$. Then,
$$ \tilde{c}_{2j}-\tilde{c}_{2j-2} = \left(\tfrac32 \tilde{\xi}_{2j} \right)^{\frac23} -  \left(\tfrac32 \tilde{\xi}_{2j-2} \right)^{\frac23}.$$
Using (\ref{eq_approx_txi2j}) in Lemma \ref{lem_asymp_ck}, we have $\tilde{\xi}_{2j}-\tilde{\xi}_{2j-2}\in \left[ \pi-\frac{5}{9\pi}\frac{2j}{(2j)^2 -1}, \pi+\frac{5}{9\pi}\frac{2j}{(2j)^2 -1}  \right] $ and 
$$\tilde{\xi}_{2j-2}^{-\frac13} \leq \left(-\frac{5\pi}{12} + j\pi - \frac{5}{18\pi}\right)^{-\frac13} \leq  \left( \frac{(2j-1)\pi}{2}\right)^{-\frac13}.  $$
Thus, by Lemma \ref{lem_a_b},
$$ \tilde{c}_{2j}-\tilde{c}_{2j-2}\leq  \left(\pi+\frac{5}{9\pi}\frac{2j}{(2j)^2 -1} \right)  \left(\frac32\right)^{\frac13} \left( \frac{(2j-1)\pi}{2}\right)^{-\frac13}.$$
If $p$ is odd, that is $p=2j+1$ for $j\geq 1$, then, using (\ref{eq_approx_txi2j1}) in Lemma \ref{lem_asymp_ck}, we have $\tilde{\xi}_{2j+1}-\tilde{\xi}_{2j-1}\in \left[ \pi-\frac{7}{3\pi}\frac{2j+1}{(2j+1)^2 -1}, \pi+\frac{7}{3\pi}\frac{2j+1}{(2j+1)^2 -1}  \right]$ and $\tilde{\xi}_{2j-1}^{-\frac13} \leq   \left( \frac{2j\pi}{2}\right)^{-\frac13}$. Thus, by Lemma \ref{lem_a_b},
$$ \tilde{c}_{2j+1}-\tilde{c}_{2j-1}\leq  \left(\pi+\frac{7}{3\pi}\frac{2j+1}{(2j+1)^2 -1} \right)  \left(\frac32\right)^{\frac13} \left( \frac{2j\pi}{2}\right)^{-\frac13}.$$
Since $\frac{7}{3\pi} > \frac{5}{9\pi}$, we have
$$\forall  p\in \{2,\ldots, k_0\},\  \tilde{c}_p-\tilde{c}_{p-2} \leq  \left(\pi+\frac{7}{3\pi}\frac{p}{p^2 -1} \right)  \left(\frac{3}{\pi}\right)^{\frac13} (p-1)^{-\frac13}$$
which proves (\ref{eq_thm_intro_gaps}).
The proof of the upper bound in (\ref{eq_thm_intro_bandes}) is similar. We estimate both $c_{2j}-\tilde{c}_{2j-1}$ and $c_{2j+1}-\tilde{c}_{2j}$ for $j\geq 1$ by using (\ref{eq_approx_xi2j}), (\ref{eq_approx_txi2j1}), (\ref{eq_approx_txi2j}) and (\ref{eq_approx_xi2j1}) to obtain that $\xi_{2j}-\tilde{\xi}_{2j-1} \in \left[ \frac{\pi}{3}-\frac{7}{3\pi}\frac{2j+\frac13}{2j(2j+\frac23)}, \frac{\pi}{3}+\frac{7}{3\pi}\frac{2j+\frac13}{2j(2j+\frac23)} \right]$ and $\xi_{2j+1}-\tilde{\xi}_{2j} \in \left[ \frac{\pi}{3}-\frac{5}{9\pi}\frac{2j+1+\frac13}{(2j+1)(2j+1+\frac23)}, \frac{\pi}{3}+\frac{5}{9\pi}\frac{2j+1+\frac13}{(2j+1)(2j+1+\frac23)} \right]$. We also have that, for every $ p\in \{2,\ldots, k_0\}$, $\tilde{\xi}_{p}^{-\frac13} \leq   \left( \frac{(p+1)\pi}{2}\right)^{-\frac13}$. Since $\frac{7}{3\pi} > \frac{5}{9\pi}$, we get the upper bound of (\ref{eq_thm_intro_bandes}) by using Lemma \ref{lem_a_b}.
\vskip 2mm

\noindent It remains to prove the lower bound in (\ref{eq_thm_intro_bandes}). We have to find a lower bound of $\tilde{c}_p-c_p$ for every $p\in \{2,\ldots, k_0\}$. Using (\ref{eq_approx_txi2j}) and (\ref{eq_approx_xi2j}) we get for every $j\geq 1$, 
$$\tilde{\xi}_{2j} - \xi_{2j} \in \left[ \frac{\pi}{6}-\frac{229}{432\pi},  \frac{\pi}{6}+\frac{229}{432\pi} \right] \subset \left[ \frac{\pi}{9},  \frac{2\pi}{9}\right].$$
We have
$$\frac{\tilde{\xi}_{2j} - \xi_{2j}}{ \xi_{2j}} \leq \frac{ \frac{2\pi}{9} }{ \frac{5\pi}{12}+\pi - \frac{7}{16\pi} } \leq  \frac{ \frac{2\pi}{9} }{ \frac{17\pi}{12} } \leq \frac{8}{51} \leq \frac16 .$$ 
Moreover,  since $\frac56 + \frac{7}{8\pi^2} <1$,
$$(\xi_{2j})^{-\frac13} \geq \left( \frac{5\pi}{12}+\frac{2j\pi}{2} + \frac{7}{16\pi} \right)^{-\frac13} \geq \left( \frac{\pi}{2} \right)^{-\frac13} (2j +1 )^{-\frac13}.$$
Thus, we take $\eta = \frac16$ and use the lower bound in (\ref{eq_lem_a_b}) to get 
$$I\left(\left( \tfrac76 \right)^{\frac23} \right) \frac{2^{\frac13}\pi^{\frac23}  }{9} (2j+1)^{-\frac13} \leq \tilde{c}_{2j}-c_{2j}.$$
For $p=2j+1$, since $\frac{\pi}{9} < \frac{\pi}{6}-\frac{97}{264\pi}$ and $\frac56 + \frac{5}{33\pi^2} <1$, and taking $\eta=\frac12$, we get a larger lower bound which is
$$I\left(\left( \tfrac32 \right)^{\frac23} \right) \frac{2^{\frac13}\pi^{\frac23}  }{9} (2j+2)^{-\frac13} \leq \tilde{c}_{2j+1}-c_{2j+1}. $$
It allows to conclude that the lower bound valid for every $p\in \{2,\ldots, k_0\}$ is the one obtained for $p$ even. This proves the lower bound in (\ref{eq_thm_intro_bandes}).
\end{demo}

\begin{demo}
 (of Corollary \ref{cor_spectral_density}). Recall that the integer $k_0$ defined in (\ref{def_k_0}) is equal to $\left[ \frac{4}{3\pi} \mathsf{c}^{\frac32} \right]$ or to  $\left[ \frac{4}{3\pi} \mathsf{c}^{\frac32} \right]-1$ where, for any $x\in \R$, $[x]$ denotes the integer part of $x$. 
Using (\ref{eq_thm_intro_bandes}), one has:
\begin{equation}\label{eq_demo_spec_dens_1}
\forall \mathsf{c} > c_0,\ 0< \frac{1}{\mathsf{c}} \sum_{p=2}^{k_0} \delta_p \leq   \frac{1}{\mathsf{c}} \sum_{p=2}^{k_0} \left( \frac{\pi}{3}+\frac{7}{3\pi}\frac{p+\frac13}{p(p+\frac23)} \right)  \left( \frac{3}{\pi} \right)^{\frac{1}{3}}\frac{1}{p^{\frac13}}.
\end{equation}
Since $x\mapsto \frac{1}{x^{\frac13}}$ and  $x\mapsto \frac{1}{x^{\frac13}} \frac{x+\frac13}{x(x+\frac23)}$ are decreasing functions on $[1,+\infty)$, by comparison between sums and integrals, 
$$\frac32 \left( (k_0 +1 )^{\frac23} -2^{\frac23} \right) \leq \sum_{p=2}^{k_0} \frac{1}{p^{\frac13}} \leq \frac32 \left( k_0^{\frac23} -1 \right)  $$
and 
$$\int_{2}^{k_0 +1} \frac{1}{x^{\frac13}}\frac{x+\frac13}{x(x+\frac23)} \mathrm{d} x \leq  \sum_{p=2}^{k_0} \frac{p+\frac13}{p(p+\frac23)}\frac{1}{p^{\frac13}} \leq \int_{1}^{k_0 } \frac{1}{x^{\frac13}}\frac{x+\frac13}{x(x+\frac23)} \mathrm{d} x.$$
Since 
\begin{align*}
\int \frac{1}{x^{\frac13}} & \frac{x+\frac13}{x(x+\frac23)} \mathrm{d} x  = -\mfrac32 \frac{1}{x^{\frac13}} - \left( \mfrac{3}{16} \right)^{\frac13} \ln\left( x^{\frac13}  +  \left( \mfrac{2}{3} \right)^{\frac13}\right) + \\
& \left( \mfrac{3}{2^7} \right)^{\frac13} \ln\left( x^{\frac23}  -  \left( \mfrac{2}{3} x \right)^{\frac13} + \left( \mfrac{2}{3} \right)^{\frac23}  \right) + \mfrac{3^{\frac56}}{2^{\frac43}} \mathrm{arctan}\left( \mfrac{2^{\frac23}}{ 3^{\frac16}} (x^{\frac13}-\mfrac{1}{3^{\frac12} } )\right),
\end{align*}
and $k_0=\left[ \frac{4}{3\pi} \mathsf{c}^{\frac32} \right]$ or $\left[ \frac{4}{3\pi} \mathsf{c}^{\frac32} \right]-1$, one gets that
$$ \frac{1}{\mathsf{c}} \sum_{p=2}^{k_0} \frac{1}{p^{\frac13}} \xrightarrow[\mathsf{c} \to +\infty]{}  \frac32 \left( \frac{4}{3\pi} \right)^{\frac23} $$
and 
$$ \frac{1}{\mathsf{c}} \sum_{p=2}^{k_0}\frac{p+\frac13}{p(p+\frac23)}\frac{1}{p^{\frac13}} \xrightarrow[\mathsf{c} \to +\infty]{} 0. $$
Thus, 
\begin{equation}\label{eq_demo_spec_dens_2}
  \frac{1}{\mathsf{c}} \sum_{p=2}^{k_0}  \left( \frac{\pi}{3}+\frac{7}{3\pi}\frac{p+\frac13}{p(p+\frac23)} \right)  \left( \frac{3}{\pi} \right)^{\frac{1}{3}}\frac{1}{p^{\frac13}} \xrightarrow[\mathsf{c} \to +\infty]{} \left( \frac23 \right)^{\frac13}.  
\end{equation}

\noindent The function $\mathsf{c} \mapsto  \frac{1}{\mathsf{c}} \sum_{p=0}^{k_0} \delta_p $ is 
increasing and: 
$$\forall \mathsf{c}>0,\ \frac{1}{\mathsf{c}} \sum_{p=0}^{k_0} \delta_p  \leq \frac{1}{\mathsf{c}} (0-\mathbf{E}_{\mathrm{min}}^0) \leq 1,$$
since the spectrum of $H$ is included in $[-V_0,+\infty]$ by Proposition \ref{prop_spec}. 
Thus $\mathsf{c} \mapsto D(\mathsf{c})$ admits a limit in $\R$ at $+\infty$. Then, (\ref{eq_demo_spec_dens_1}) and (\ref{eq_demo_spec_dens_2}) imply (\ref{eq_cor_spectral_density}).
\end{demo}

\section{Spectral bands and spectral gaps in the semiclassical regime}\label{sec_bands_gaps}

\vskip 5mm
\noindent Proposition \ref{prop_Ekmin_Ekmax} allows us to identify the spectral band edges among the solutions of (\ref{eq_vp_Ha_0_a}), (\ref{eq_vp_Ha_0_b}), (\ref{eq_vp_Ha_d_a}) and (\ref{eq_vp_Ha_d_b}). Using proofs similar to those of the estimates and asymptotics of $\mathbf{E}_{\mathrm{min}}^0$, $\mathbf{E}_{\mathrm{max}}^0$ and $\mathbf{E}_{\mathrm{min}}^1$, one gets the following estimates for the spectral band edges in the semiclassical regime.

\begin{prop}\label{prop_asym_band_edges}
Let $j\geq 0$.
\begin{enumerate}
\item Let $\tau >0$ be arbitrary small. For every $\mathsf{h} \in (0, (\tilde{a}_{j+1}+\tau)^{-\frac32}]$,
\begin{equation}
 \mathbf{E}_{\mathrm{min}}^{2j}  =  -\mathsf{h}^{-\frac23}+\tilde{a}_{j+1} -\alpha \sqrt{3} \frac{(u'(-\tilde{a}_{j+1}))^2 }{\tilde{a}_{j+1} }\ex^{-\frac43 \mathsf{h}^{-1} + 2\tilde{a}_{j+1}\mathsf{h}^{-\frac13}} \left( 1 +  \mathcal{O} \left(\mathsf{h}^{\frac13} \right) \right) \label{prop_asym_band_edges_min2j}
\end{equation}
and for every $\mathsf{h} \in (0, c_{2j}^{-\frac32})$,
\begin{equation}
\mathbf{E}_{\mathrm{max}}^{2j}  =  -\mathsf{h}^{-\frac23}+\tilde{a}_{j+1}  +\alpha \sqrt{3} \frac{(u'(-\tilde{a}_{j+1}))^2  }{\tilde{a}_{j+1} }\ex^{-\frac43 \mathsf{h}^{-1} + 2\tilde{a}_{j+1}\mathsf{h}^{-\frac13}} \left( 1 +  \mathcal{O} \left(\mathsf{h}^{\frac13} \right) \right). \label{prop_asym_band_edges_max2j}
\end{equation}

\item Let $\tau >0$ be arbitrary small. For every $\mathsf{h} \in (0, (a_{j+1}+\tau)^{-\frac32}]$,
\begin{equation}
 \mathbf{E}_{\mathrm{min}}^{2j+1}  =  -\mathsf{h}^{-\frac23}+a_{j+1} -\alpha \sqrt{3} (u(-a_{j+1}))^2 \ex^{-\frac43 \mathsf{h}^{-1} + 2a_{j+1}\mathsf{h}^{-\frac13}} \left( 1 +  \mathcal{O} \left(\mathsf{h}^{\frac13} \right) \right) \label{prop_asym_band_edges_min2j1}
\end{equation}
and for every $\mathsf{h} \in (0, c_{2j+1}^{-\frac32})$,
\begin{equation}
 \mathbf{E}_{\mathrm{max}}^{2j+1}  =  -\mathsf{h}^{-\frac23}+a_{j+1}+\alpha \sqrt{3} (u(-a_{j+1}))^2 \ex^{-\frac43 \mathsf{h}^{-1} + 2a_{j+1}\mathsf{h}^{-\frac13}} \left( 1 +  \mathcal{O} \left(\mathsf{h}^{\frac13} \right) \right).
\label{prop_asym_band_edges_max2j1}
\end{equation}
\item Let $h\in (0,c_1^{-\frac32})$ be arbitrary small. For every $p\in \{0,\ldots,  \left[\frac{4}{3\pi}\frac{1}{\mathsf{h}} \right]  \}$, (\ref{prop_asym_band_edges_min2j}) and (\ref{prop_asym_band_edges_max2j}) hold true if $p=2j$ is even and  (\ref{prop_asym_band_edges_min2j1}) and (\ref{prop_asym_band_edges_max2j1}) hold true if $p=2j+1$ is odd. 
\end{enumerate}
\end{prop}

\vskip 2mm
\noindent Note that this amount to the behaviour of $\mathbf{E}_{\mathrm{min}}^p + \mathbf{V}(0)$ and of $\mathbf{E}_{\mathrm{max}}^p + \mathbf{V}(0)$. We notice that one has the correct scaling near the minimum value of the potential.

\vskip 3mm

\begin{demo}
For $j=0$ we already obtained the estimate of $\mathbf{E}_{\mathrm{min}}^{0}$.  For every $j\geq 1$, $\mathbf{E}_{\mathrm{min}}^{2j}$ is the unique solution of (\ref{eq_vp_Ha_d_a}) with $-\mathsf{c} -\mathbf{E}_{\mathrm{min}}^{2j} \in [-\tilde{a}_{j+1},-\tilde{c}_{2j-1})$. Using Lemma \ref{lem_minor_uv_alpha}, the function $x\mapsto -\frac{x}{(u'(x))^2}$ is greater than $\alpha$ on the interval $[-\tilde{a}_{j+1},-\tilde{c}_{2j-1})$. Thus, the scheme of the proof of the estimate of $\mathbf{E}_{\mathrm{min}}^{0}$ can be followed and leads to (\ref{prop_asym_band_edges_min2j}).
\vskip 3mm

\noindent For every $j\geq 0$,  $\mathbf{E}_{\mathrm{max}}^{2j}$ is the unique solution of (\ref{eq_vp_Ha_0_b})  with $-\mathsf{c} -\mathbf{E}_{\mathrm{max}}^{2j} \in [-c_{2j},-\tilde{a}_{j+1})$. Since the function $x\mapsto -\frac{x}{(u'(x))^2}$ is greater than $\alpha$ on the interval $[-c_{2j},-\tilde{a}_{j+1}]$, the scheme of the proof of the estimate of  $\mathbf{E}_{\mathrm{max}}^{0}$ can be followed and leads to (\ref{prop_asym_band_edges_max2j}).
\vskip 3mm

\noindent For every $j\geq 0$, $\mathbf{E}_{\mathrm{min}}^{2j+1}$ is the unique solution of (\ref{eq_vp_Ha_0_a}) with $-\mathsf{c} -\mathbf{E}_{\mathrm{min}}^{2j+1} \in [-a_{j+1}, -\tilde{c}_{2j})$. Since the function $\frac{1}{u^2}$ is greater than $\frac{\alpha}{2}$ on this interval, the scheme of the proof of the estimate of $\mathbf{E}_{\mathrm{min}}^{1}$ can be followed and leads to (\ref{prop_asym_band_edges_min2j1}).
\vskip 3mm

\noindent  For every $j\geq 0$, $\mathbf{E}_{\mathrm{max}}^{2j+1}$ is the unique solution of (\ref{eq_vp_Ha_d_b}) with $-\mathsf{c} -\mathbf{E}_{\mathrm{max}}^{2j+1} \in [-c_{2j+1},-a_{j+1})$. Since the function $\frac{1}{u^2}$ is greater than than $\frac{\alpha}{2}$ on the interval $ [-c_{2j+1},-a_{j+1}]$, combining the proofs of the estimates of  $\mathbf{E}_{\mathrm{max}}^{0}$  and $\mathbf{E}_{\mathrm{min}}^{1}$ leads to (\ref{prop_asym_band_edges_max2j}).
\vskip 3mm

\noindent The last statement is a direct consequence of the counting of the number of rescaled spectral bands in the range of $\mathbf{V}$ done in Theorem \ref{thm_intro_pbands} and the fact that both sequences $(c_p^{-\frac32})_{p\geq 0}$ and $(\mathfrak{a}_p^{-\frac32})_{p\geq 0}$ are strictly decreasing and interlaced.

\end{demo}
\vskip 5mm

\noindent These estimates lead to the estimates of the widths of the $p$-th spectral band and the $p$-th gap in the semiclassical regime.

\begin{prop}\label{prop_asymp_width}
Let $j\geq 0$. For every $\mathsf{h} \in (0, c_{2j}^{-\frac32})$,
\begin{equation}
\delta_{2j} =  2 \alpha \sqrt{3} \mfrac{(u'(-\tilde{a}_{j+1}))^2  }{\tilde{a}_{j+1} }  \ex^{-\frac43 \mathsf{h}^{-1} + 2\tilde{a}_{j+1}\mathsf{h}^{-\frac13}} \left( 1 +  \mathcal{O} \left(\mathsf{h}^{\frac13} \right) \right), \label{eq_prop_width_2j}
\end{equation}
and for every $\mathsf{h} \in (0, c_{2j+1}^{-\frac32})$,
\begin{equation}
\delta_{2j+1} =  2 \alpha \sqrt{3} (u(-a_{j+1}))^2  \ex^{-\frac43 \mathsf{h}^{-1} + 2a_{j+1}\mathsf{h}^{-\frac13}} \left( 1 +  \mathcal{O} \left(\mathsf{h}^{\frac13} \right) \right). \label{eq_prop_width_2j1} 
\end{equation}
\end{prop}

\begin{demo}
The estimate (\ref{eq_prop_width_2j}) is obtained from (\ref{prop_asym_band_edges_min2j}), (\ref{prop_asym_band_edges_max2j}) and the estimate (\ref{eq_prop_width_2j1}) is a direct consequence of (\ref{prop_asym_band_edges_max2j1}) and (\ref{prop_asym_band_edges_min2j1}).
\end{demo}

\noindent Proposition \ref{prop_asym_band_edges} and Proposition \ref{prop_asymp_width} together imply Theorem \ref{thm_intro_width}. 

\begin{demo} (of Theorem \ref{thm_intro_width}). The first statement of Theorem \ref{thm_intro_width} is about the convergence of the rescaled spectral bands to the zeroes of the Airy function $Ai$ and its derivative, and is a consequence of the first two terms of (\ref{prop_asym_band_edges_min2j}), (\ref{prop_asym_band_edges_max2j}), (\ref{prop_asym_band_edges_min2j1}) and (\ref{prop_asym_band_edges_max2j1}). More precisely, going back to the intermediate computations done in the end of the proof of the second point of Theorem \ref{thm_bottom}, for every $j\geq 0$ and every $\mathsf{h} \in (0, c_{2j}^{-\frac32})$,
$$\mathsf{h}^{-\frac23} + \mathbf{E}_{\mathrm{min}}^{2j} = \tilde{a}_{j+1} + \mathcal{O} \left( \ex^{-\frac43  \frac{1}{\mathsf{h}} (1-\tilde{a}_{j+1}\mathsf{h}^{\frac23})^{\frac32}} \right)$$
and 
$$\mathsf{h}^{-\frac23} + \mathbf{E}_{\mathrm{max}}^{2j} = \tilde{a}_{j+1} + \mathcal{O} \left( \ex^{-\frac43 \frac{1}{\mathsf{h}} (1-\tilde{a}_{j+1}\mathsf{h}^{\frac23})^{\frac32}} \right)$$ 
and both $\mathsf{h}^{-\frac23} + \mathbf{E}_{\mathrm{min}}^{2j} $ and $\mathsf{h}^{-\frac23} + \mathbf{E}_{\mathrm{max}}^{2j}$ tends to $\tilde{a}_{j+1}$ when $\mathsf{h}$ tends to $0$. 

\noindent Using (\ref{prop_asym_band_edges_min2j1}) and (\ref{prop_asym_band_edges_max2j1}) we prove similarly that both $\mathsf{h}^{-\frac23} + \mathbf{E}_{\mathrm{min}}^{2j+1} $ and $\mathsf{h}^{-\frac23} + \mathbf{E}_{\mathrm{max}}^{2j+1}$ tends to $a_{j+1}$ when $\mathsf{h}$ tends to $0$. 

\noindent The estimates for $\delta_{2j}$ and $\delta_{2j+1}$ in Theorem \ref{thm_intro_width} are exactly the same as those in Proposition \ref{prop_asymp_width}. 

\noindent The last statement is a direct consequence of the counting of the number of rescaled spectral bands in the range of $\mathbf{V}$ done in Theorem \ref{thm_intro_pbands} and the fact that the sequence $(c_p^{-\frac32})_{p\geq 0}$ is strictly decreasing. 
\end{demo}

\begin{prop}\label{prop_asymp_gap}
Let $j\geq 0$. Let $\tau >0$ be arbitrary small. For every $\mathsf{h}\in (0,(a_{j+1}+\tau)^{-\frac32}]$,
 \begin{equation}
   \gamma_{2j} =  a_{j+1}-\tilde{a}_{j+1} - \alpha \sqrt{3} (u(-a_{j+1}))^2 \ex^{-\frac43 \mathsf{h}^{-1} + 2a_{j+1}\mathsf{h}^{-\frac13}} \left( 1 +  \mathcal{O} \left(\mathsf{h}^{\frac13} \right) \right), \label{eq_prop_gap_2j}
\end{equation}
and for every $\mathsf{h}\in (0,(\tilde{a}_{j+2}+\tau)^{-\frac32} ]$,
 \begin{equation}
\gamma_{2j+1} =  \tilde{a}_{j+2}-a_{j+1} - \alpha \sqrt{3}   \mfrac{(u'(-\tilde{a}_{j+2}))^2  }{ \tilde{a}_{j+2} } \ex^{-\frac43 \mathsf{h}^{-1} + 2\tilde{a}_{j+2}\mathsf{h}^{-\frac13}} \left( 1 +  \mathcal{O} \left(\mathsf{h}^{\frac13} \right) \right).\label{eq_prop_gap_2j1}
\end{equation} 
\end{prop}
 
\begin{demo}
The estimate (\ref{eq_prop_gap_2j}) is a consequence of (\ref{prop_asym_band_edges_min2j1}), (\ref{prop_asym_band_edges_max2j}) and the fact that $a_{j+1} > \tilde{a}_{j+1}$.

\noindent  The estimate (\ref{eq_prop_gap_2j1}) is a consequence of (\ref{prop_asym_band_edges_min2j}), (\ref{prop_asym_band_edges_max2j1}) and the fact that $\tilde{a}_{j+2} > a_{j+1}$.
\end{demo}

\noindent This proposition implies directly Theorem \ref{thm_intro_gap} of the introduction. The last statement of Theorem \ref{thm_intro_gap} is a direct consequence of the counting of the number of rescaled spectral gaps in the range of $\mathbf{V}$ done in Theorem \ref{thm_intro_pbands} and the fact that the sequence $(\mathfrak{a}_p^{-\frac32})_{p\geq 0}$ is strictly decreasing.

\section{Conclusion}\label{sec_conclusion}

Let us summarize some of the results obtained in this article.

\begin{enumerate}

\item We have been able to get estimates of the widths of the spectral bands and the spectral gaps in the semiclassical regime for a periodic potential which is not analytic and not even differentiable at its maxima or minima. It is an example of non-regular periodic potential of physical interest. This was done thanks to the accurate asymptotic expansions of the classical Airy functions and their derivatives and thanks to a bootstrap analysis argument developed in the proof of Theorem \ref{thm_bottom}. 

\item The spectral bands are exponentially thin, but the exponential decay is not as good as the classical exponential decay for a regular potential. The widening of the spectral bands is similar to a Gevrey $3$ effect of the singularity of the potential at its minima.

\item We defined a notion of semiclassical regime. 
In this regime the results are stated for a fixed value of the semiclassical parameter $\mathsf{h}$ in an interval $(0,\mathsf{h}_0)$. We give explicit values of $\mathsf{h}_0$ for which the behaviour of the considered spectral quantity is similar to its behaviour in the semiclassical limit for every $\mathsf{h}\in (0,\mathsf{h}_0)$. In our results, these explicit values are the zeroes of the Airy function $Ai$ and its derivative and the zeroes of the derivative of the odd canonical solution of the Airy equation.

\item We count the number of spectral bands which lie in the range of the potential $V$. This number depends only on the value of the counting parameter $\mathsf{c}$ compared to the values of the zeroes of the canonical solutions of the Airy equation and their derivatives. In Proposition \ref{thm_zero_edges} we give a dynamical picture of the successives entrances in the range of $V$ of the spectral bands and gaps when $\mathsf{c}$ grows and takes the successives values $c_p$ and $\tilde{c}_p$ for any $p\geq 0$.

\item We give explicit bounds for the widths of the spectral bands and gaps which lie in the range of the potential. The bounds for the $p$-th spectral band or gap depend only on $p$. 

\item We prove an upper bound on the integrated spectral density in the range of the potential.

\end{enumerate}

\noindent We address some open questions which naturally arose in our research on the periodic Airy-Schr\"odinger operator. 

\begin{enumerate}
 \item  A first question is to generalize our results for a potential which is no longer our explicit potential. We consider a function $V$ which is analytic, such that $V(0)=0$ and $V'(0)\neq 0$ and then we consider the $2L_0$-periodic function $W$ defined on $[-L_0,L_0]$ by:
$$\forall x \in [-L_0,L_0],\ W(x)=V(|x|).$$  
Using perturbation theory techniques like those developed in \cite{BBL} or \cite{BG}, one would like to obtain similar results as Theorem \ref{thm_intro_transition}, \ref{thm_bottom}, \ref{thm_intro_width} or  \ref{thm_intro_gap}.
\vskip 3mm

\item We hope to use our results on the spectrum of $H$ to study other periodic operators at least in dimension $2$. First results were presented in \cite{BL2}. In this case, we were to tackle the case where our operator decompose into a tensor product of two periodic Airy-Schr\"odinger operators with possibly two different characteristic lengths $L_0$ and $L_1$. In this case the spectrum of the two-dimensional operator is the superposition of the band spectra of the two one-dimensional operators. It would certainly lead to difficulties linked to the compared arithmetic natures of $L_0$ and $L_1$, like those presented in \cite{FK}.
\vskip 3mm

\item  Another interesting generalization would be to obtain asymptotic results similar to Theorem \ref{thm_intro_width} or Theorem \ref{thm_intro_gap} for our non-analytic potential, in all the regimes considered in \cite{marz}. 
\vskip 3mm

\item  Another question is the one of the meromorphic continuation of the resolvent of the periodic Airy-Schr\"odinger operator to the spectral bands, using techniques like those in \cite{GHH,KM,ramond,Z2}. It would lead to the question of the existence and the description of the resonances for the periodic Airy-Schr\"odinger operator. The Gevrey $3$ behaviour that is observed in the estimates may give rise to other estimates.
\end{enumerate}

\vskip 5mm

\appendix

\section{Proofs of Propositions \ref{prop_asym_uv}, \ref{prop_loc_ck}, \ref{prop_variations_uv} and of Lemma \ref{lem_asymp_ck} }\label{sec_proofs}

\noindent We start with the proof of Proposition \ref{prop_asym_uv}.
\vskip 3mm

\begin{demo}(of Proposition \ref{prop_asym_uv}). Let $x> 0$ and  $\xi = \frac{2}{3} x^{\frac{3}{2}}$. The functions $Ai$, $Bi$ and their derivatives are related to the Bessel functions through the relations (see \cite[10.4.15 and after]{AS})
$$Ai(-x)=\frac13 \sqrt{x} \left( J_{\frac13}(\xi)+J_{-\frac13}(\xi) \right),\quad  Bi(-x)= \sqrt{ \frac{x}{3} } \left( J_{-\frac13}(\xi)-J_{\frac13}(\xi) \right)$$
and
$$Ai'(-x)=-\frac13 x \left( J_{-\frac23}(\xi)-J_{\frac23}(\xi) \right),\quad  Bi'(-x)=  \frac{x}{\sqrt{3}}  \left( J_{-\frac23}(\xi)+J_{\frac23}(\xi) \right). $$
Thus, we have the following expressions for the Airy functions and their derivatives on the negative half-line:

\begin{align}
Ai(-x) & =\pi^{-\frac12} x^{-\frac14} \left( \cos\left(\xi - \tfrac{\pi}{4}\right)P\left(\tfrac{1}{3},\xi \right) - \sin\left(\xi - \tfrac{\pi}{4}\right) Q\left(\tfrac{1}{3},\xi \right) \right), \label{eq_Ai_asym} \\
 Ai'(-x) & =\pi^{-\frac12} x^{\frac14} \left( \sin\left(\xi - \tfrac{\pi}{4}\right)P\left(\tfrac{2}{3},\xi \right) + \cos \left(\xi - \tfrac{\pi}{4}\right)Q\left(\tfrac{2}{3},\xi\right) \right),\label{eq_Aip_asym}\\
 Bi(-x) & =-\pi^{-\frac12} x^{-\frac14} \left( \sin\left(\xi - \tfrac{\pi}{4}\right)P\left(\tfrac{1}{3},\xi\right) + \cos\left(\xi - \tfrac{\pi}{4}\right)Q\left(\tfrac{1}{3},\xi\right) \right),\label{eq_Bi_asym}\\
 Bi'(-x)& =\pi^{-\frac12} x^{\frac14} \left( \cos\left(\xi - \tfrac{\pi}{4}\right)P\left(\tfrac{2}{3},\xi\right) - \sin\left(\xi - \tfrac{\pi}{4}\right)Q\left(\tfrac{2}{3},\xi\right) \right).\label{eq_Bip_asym}
\end{align}
\vskip 3mm

\noindent Before getting similar expressions for the canonical solutions $u$ and $v$, let us start by rewriting $u$ and $v$, observing that $\frac{Bi'(0)}{Ai'(0)} = -\sqrt{3} =- \tan(\frac{\pi}{3})$:
\begin{equation} \label{eq_u_trigo}
\begin{split}
\forall x \in \R,\ u(x) & =   \pi  (Bi'(0) Ai(x)-Ai'(0) Bi(x)) \\
&  =   -2 \pi Ai'(0) \left( \cos\left(\tfrac{\pi}{3}\right) Bi(x) + \sin\left(\tfrac{\pi}{3}\right) Ai(x)\right).
\end{split}
\end{equation}
Similarly,
\begin{equation}\label{eq_v_trigo}
 \forall x \in \R,\ v(x)  =  2 \pi Ai(0) \left( \cos\left(\tfrac{\pi}{3}\right) Bi(x) - \sin\left(\tfrac{\pi}{3}\right) Ai(x)\right).
\end{equation}
Combining (\ref{eq_u_trigo}), (\ref{eq_Ai_asym}) and (\ref{eq_Bi_asym}) one gets, for every $x>0$,
{\small \begin{align*}
u(-x) & =  -2 \pi Ai'(0) \left( \cos\left(\tfrac{\pi}{3}\right) Bi(-x) + \sin\left(\tfrac{\pi}{3}\right) Ai(-x)\right)  \\
 & =  -2 \pi^{\frac12} x^{-\frac14} Ai'(0) \left( \left(-\cos\left(\tfrac{\pi}{3}\right)\sin\left(\xi-\tfrac{\pi}{4}\right)  + \sin\left(\tfrac{\pi}{3}\right)\cos\left(\xi-\tfrac{\pi}{4}\right)   \right)P\left(\tfrac{1}{3},\xi\right) + \right.  \\
&  \left. \qquad \left(-\cos\left(\tfrac{\pi}{3}\right)\cos\left(\xi-\tfrac{\pi}{4}\right)  - \sin\left(\tfrac{\pi}{3}\right)\sin\left(\xi-\tfrac{\pi}{4}\right)   \right)Q\left(\tfrac{1}{3},\xi\right)   \right)  \\
 & =   - 2 \pi^{\frac12} x^{-\frac14} Ai'(0) \left( \sin\left(-\xi + \tfrac{\pi}{4} +\tfrac{\pi}{3} \right)P\left(\tfrac{1}{3},\xi\right) - \cos\left(-\xi + \tfrac{\pi}{4} +\tfrac{\pi}{3} \right)Q\left(\tfrac{1}{3},\xi\right)  \right) \\
 & =     2 \pi^{\frac12} x^{-\frac14} Ai'(0) \left(  \sin\left(\xi - \tfrac{7\pi}{12}\right)P\left(\tfrac{1}{3},\xi\right) + \cos\left(\xi - \tfrac{7\pi}{12}\right)Q\left(\tfrac{1}{3},\xi\right)  \right). 
\end{align*}}
\noindent By derivating (\ref{eq_u_trigo}) and doing similar computations as in (\ref{eq_u_asym}) one gets, for every $x>0$,
$$u'(-x) = - 2 \pi^{\frac12} x^{\frac14} Ai'(0) \left(  \cos\left(\xi - \tfrac{7\pi}{12}\right)P\left(\tfrac{2}{3},\xi\right)  - \sin\left(\xi - \tfrac{7\pi}{12}\right)Q\left(\tfrac{2}{3},\xi\right)   \right).$$
The expressions for $v$ and $v'$ are obtained the same way.
\end{demo}

\vskip 3mm

\noindent Before proving simultaneously Proposition \ref{prop_loc_ck} and Proposition \ref{prop_variations_uv}, we need the following technical lemma. 
\vskip 3mm

\begin{lem}\label{lem_PQ}
For every $\xi>\frac{1}{\sqrt{26}}$, $P(\frac{1}{3},\xi)>0$ and we have
\begin{equation}\label{eq_ineq_QP13}
\forall \xi>\frac{1}{\sqrt{13}},\ \left| \frac{Q(\frac{1}{3},\xi)}{ P(\frac{1}{3},\xi)}\right| <  \frac{5}{36\xi}. 
\end{equation}

\noindent For every $\xi>\frac{1}{\sqrt{22}}$, $P(\frac{2}{3},\xi)>0$ and we have
\begin{equation}\label{eq_ineq_QP23}
\forall \xi>\frac{1}{\sqrt{11}},\ \left| \frac{Q(\frac{2}{3},\xi)}{ P(\frac{2}{3},\xi)}\right| < \frac{7}{12\xi}.
\end{equation}
\end{lem}

\begin{demo}
Using \cite[9.2.9 and 9.2.10]{AS}, we have
\begin{equation}\label{eq_maj_P13}
\forall \xi >0,\ \left| P\left(\tfrac{1}{3},\xi\right) -1 \right| \leq \frac{5\times 77}{81 \times 128 \times \xi^2} < \frac{1}{26\xi^2}
\end{equation}
and
\begin{equation}\label{eq_maj_Q13}
 \forall \xi >0,\ \left| Q\left(\tfrac{1}{3},\xi\right) +\frac{5}{72\xi} \right| \leq \frac{5\times 77\times 221}{6\times 9^3 \times 8^3 \times \xi^3} < \frac{1}{26\xi^3}.
\end{equation}
In particular, we deduce from (\ref{eq_maj_P13}) that for every $\xi>\frac{1}{\sqrt{26}}$, $P(\frac{1}{3},\xi)>0$. Then, (\ref{eq_maj_P13}) and (\ref{eq_maj_Q13}) imply 
$$\forall \xi>\frac{1}{\sqrt{13}},\ \left| \frac{Q(\frac{1}{3},\xi)}{ P(\frac{1}{3},\xi)}\right| < \frac{\frac{5}{72\xi}-\frac{1}{26\xi^3}}{1-\frac{1}{26\xi^2}} < \frac{5}{36\xi},$$
which proves (\ref{eq_ineq_QP13}). Indeed, for $\xi>\frac{1}{\sqrt{13}}$, $ \frac{1}{1-\frac{1}{26\xi^2}} <2$.

\noindent Using again \cite[9.2.9 and 9.2.10]{AS}, we also have
\begin{equation}\label{eq_maj_P23}
\forall \xi >0,\ \left| P\left(\tfrac{2}{3},\xi\right) -1 \right| \leq \frac{7\times 65}{81 \times 128 \times \xi^2} < \frac{1}{22\xi^2} 
\end{equation}
and
\begin{equation}\label{eq_maj_Q23}
 \forall \xi >0,\ \left| Q\left(\tfrac{2}{3},\xi\right) -\frac{7}{72\xi} \right| \leq \frac{7\times 65\times 209}{6\times 9^3 \times 8^3 \times \xi^3} < \frac{1}{22\xi^3}.
\end{equation}
In particular, for every $\xi>\frac{1}{\sqrt{22}}$, $P(\frac{2}{3},\xi)>0$. Then, (\ref{eq_maj_P23}) and (\ref{eq_maj_Q23}) imply 
$$\forall \xi>\frac{1}{\sqrt{11}},\ \left| \frac{Q(\frac{2}{3},\xi)}{ P(\frac{2}{3},\xi)}\right| < \frac{\frac{7}{72\xi}+\frac{1}{22\xi^3}}{1-\frac{1}{22\xi^2}} < \frac{7}{36\xi} + \frac{1}{11\xi^3} < \frac{7}{12\xi},$$
which proves (\ref{eq_ineq_QP23}). Indeed, for $ \xi>\frac{1}{\sqrt{11}}$, $\frac{1}{ 1-\frac{1}{22\xi^2}} <2$ and $ \frac{7}{36\xi} + \frac{1}{11\xi^3} < \frac{7}{12\xi}$. 
\end{demo}
%
\vskip 3mm

\begin{demo}(of Propositions \ref{prop_loc_ck} and \ref{prop_variations_uv}). Before starting the proof, since $-(\frac{\pi}{2})^{\frac23} > -\tilde{a}_1 > -c_0$, we stress that $v'$ does not vanish in the interval $(-(\frac{\pi}{2})^{\frac23},+\infty)$ which justifies the starting point for the numbering of the $-c_{2j}$. Similarly, the only root of $v$ in the interval $(-(\frac{5\pi}{4})^{\frac23},+\infty)$ is $0$, which justifies the numbering of the  $-c_{2j+1}$, $u$ does not vanish in $(-(\frac{3\pi}{4})^{\frac23},+\infty)$ and the only root of $u'$ in the interval $(-(\frac{3\pi}{2})^{\frac23},+\infty)$ is $0$, which justifies the numbering of respectively the $-\tilde{c}_{2j}$ and the $-\tilde{c}_{2j+1}$.
\vskip 2mm

\noindent We prove only the assertion on $u$ and $-\tilde{c}_{2j}$, the others are proved in a completely similar way. We use the following method: thanks to the bounds of $\frac{Q}{P}$ for $\xi > \frac{1}{\sqrt{11} }$, we show that $u$ changes its sign at the two boundary values of the considered interval while $u'$ is of constant sign in the interval.

\noindent Let $j\geq 0$. Using (\ref{eq_u_asym}) for $x=\left(\frac32 (j\pi + \frac{2\pi}{3}) \right)^{\frac23}$ and thus $\xi=j\pi + \frac{2\pi}{3}$, one gets:
\begin{align*}
u\left(-\left(\tfrac32 (j\pi + \tfrac{2\pi}{3}) \right)^{\frac23}\right) & =   2 \pi^{\frac12} \left(\tfrac32 (j\pi + \tfrac{2\pi}{3}) \right)^{-\frac16} Ai'(0)(-1)^j \sin(\tfrac{\pi}{12})\times \\
&  P(\tfrac{1}{3},j\pi + \tfrac{2\pi}{3}) \left( 1 + \mathrm{cotan}(\tfrac{\pi}{12}) \frac{Q(\tfrac{1}{3},j\pi + \tfrac{2\pi}{3})}{P(\tfrac{1}{3},j\pi + \tfrac{2\pi}{3})} \right).
   \end{align*}
But, $\mathrm{cotan}(\frac{\pi}{12})=2+\sqrt{3}$ and since using (\ref{eq_ineq_QP13}),
$$\left| \mathrm{cotan}(\tfrac{\pi}{12}) \frac{Q(\frac{1}{3},j\pi + \frac{2\pi}{3})}{P(\frac{1}{3},j\pi + \frac{2\pi}{3})} \right| < \frac{5(2+\sqrt{3}) }{36} < \frac23,$$
we get that 
$$\left( 1 + \mathrm{cotan}(\tfrac{\pi}{12}) \frac{Q(\frac{1}{3},j\pi + \frac{2\pi}{3})}{P(\frac{1}{3},j\pi + \frac{2\pi}{3})} \right) >0.$$
Since $\sin(\frac{\pi}{12})>0$, $P(\frac{1}{3},j\pi + \frac{2\pi}{3})>0$ and  $Ai'(0)<0$,
$$(-1)^j u\left(-\left(\tfrac32 (2j\pi + \tfrac{9\pi}{12}) \right)^{\frac23}\right) <0.$$
Then, using (\ref{eq_u_asym}) for $x=\left(\frac32 (j\pi + \frac{\pi}{3}) \right)^{\frac23}$ and thus $\xi=j\pi + \frac{\pi}{3}$, one gets:
\begin{align*}
u\left(-\left(\tfrac32 (j\pi + \tfrac{\pi}{3}) \right)^{\frac23}\right) & = 2 \pi^{\frac12} \left(\tfrac32 (j\pi + \tfrac{\pi}{3}) \right)^{-\frac16} Ai'(0) (-1)^j \sin(-\tfrac{\pi}{12}) \times \\
 & P(\tfrac{1}{3},j\pi + \tfrac{\pi}{3})  \left( 1 + \mathrm{cotan}(-\tfrac{\pi}{12}) \frac{Q(\frac{1}{3},j\pi + \frac{\pi}{3})}{P(\frac{1}{3},j\pi + \frac{\pi}{3})} \right).
\end{align*}
But, using (\ref{eq_ineq_QP13}),
$$\left| \mathrm{cotan}(-\tfrac{\pi}{12}) \frac{Q(\frac{1}{3},j\pi + \frac{\pi}{3})}{P(\frac{1}{3},j\pi + \frac{\pi}{3})} \right| < \frac{5(2+\sqrt{3}) }{36} < \frac23,$$
and we get that 
$$\left( 1 + \mathrm{cotan}(-\tfrac{\pi}{12}) \frac{Q(\frac{1}{3},j\pi + \frac{\pi}{3})}{P(\frac{1}{3},j\pi + \frac{\pi}{3})} \right) <0.$$
Since $P(\frac{1}{3},j\pi + \frac{\pi}{3})>0$,
$$(-1)^j u\left(-\left(\tfrac32 (j\pi + \tfrac{\pi}{3}) \right)^{\frac23}\right) > 0.$$
If $x\in [\left(\frac32 (j\pi + \frac{\pi}{3}) \right)^{\frac23},\left(\frac32 (j\pi + \frac{2\pi}{3}) \right)^{\frac23} ]$ then $\xi-\frac{7\pi}{12} \in [j\pi -\frac{\pi}{12}, j\pi+\frac{\pi}{12}]$ and one has 
$$\frac{\sqrt{2} }{2} \leq (-1)^j\cos\left(\xi-\tfrac{7\pi}{12}  \right) \leq \frac{\sqrt{2+\sqrt{3} } }{2} \ \mbox{ and }\ \frac{1}{2+\sqrt{3} } \leq \tan\left(\xi-\tfrac{7\pi}{12}  \right) \leq 1.$$
Moreover, using (\ref{eq_ineq_QP23}) and $\xi \geq \frac{\pi}{2}$,
$$\left| \frac{Q(\frac{2}{3},\xi) }{P(\frac{2}{3},\xi) } \right| < \frac{7}{12\xi} \leq \frac{14}{12\pi } < \frac12.$$
Then, using (\ref{eq_up_asym}), for every $x \in \left[\left(\frac32 (j\pi + \frac{\pi}{3}) \right)^{\frac23},\left(\frac32 (j\pi + \frac{2\pi}{3}) \right)^{\frac23} \right]$,
{\small $$(-1)^j u'(-x)=2 \pi^{\frac12} x^{\frac14} Ai'(0)   \cos(\xi - \tfrac{7\pi}{12}) P\left(\tfrac{2}{3},\xi\right) \left(1 - \tan(\xi - \tfrac{7\pi}{12})\frac{ Q(\frac{2}{3},\xi)}{ P(\frac{2}{3},\xi)}\right) > 0.$$}
We deduce that $u$ is continuous, strictly increasing for $j$ even (respectively decreasing for $j$ odd) from a negative value to a positive one (respectively from a positive value to a negative one) and thus has a unique zero in the interval $\left[-\left(\frac32 (j\pi + \frac{2\pi}{3}) \right)^{\frac23},-\left(\frac32 (j\pi + \frac{\pi}{3}) \right)^{\frac23} \right]$, for every $j\geq 0$. 

\noindent It remains to verify that $u$ does not vanish on the interval 
$$\left(-\left(\tfrac32 ((j+1)\pi + \tfrac{\pi}{3}) \right)^{\frac23},-\left(\tfrac32 (j\pi + \tfrac{2\pi}{3}) \right)^{\frac23}\right).$$
If  $x\in (\left(\frac32 (j\pi + \frac{2\pi}{3}) \right)^{\frac23}, \left(\frac32 ((j+1)\pi + \frac{\pi}{3}) \right)^{\frac23} )$ then $\xi-\frac{7\pi}{12} \in [j\pi +\frac{\pi}{12}, j\pi +\frac{3\pi}{4}]$ and one has 
$$ \frac{1}{2\sqrt{2+\sqrt{3} } } \leq (-1)^j \sin\left(\xi-\tfrac{7\pi}{12}  \right) \leq \frac{\sqrt{2}}{2}  \ \mbox{ and }\ 1 \leq \mathrm{cotan} \left(\xi-\tfrac{7\pi}{12}  \right) \leq 2+ \sqrt{3}.$$
Moreover, using (\ref{eq_ineq_QP23}) and $\xi \geq \frac{2\pi}{3}$,
$$\left|  \mathrm{cotan}\left(\xi-\tfrac{7\pi}{12}  \right)\frac{Q\left(\tfrac{1}{3},\xi\right) }{P\left(\tfrac{1}{3},\xi\right) } \right| < (2+\sqrt{3})\cdot \frac{7}{36\xi} < \frac23$$
and using (\ref{eq_up_asym}), for every $x \in \left(\left(\frac32 (j\pi + \frac{2\pi}{3}) \right)^{\frac23}, \left(\frac32 ((j+1)\pi + \frac{\pi}{3}) \right)^{\frac23} \right)$,
\begin{equation}
 \begin{split}
 &(-1)^j u(-x)= \nonumber \\
& 2 \pi^{\frac12} x^{-\frac14} Ai'(0) (-1)^j \sin(\xi - \tfrac{7\pi}{12}) P(\tfrac{1}{3},\xi) \left( 1+ \mathrm{cotan}(\xi - \tfrac{7\pi}{12})\frac{ Q(\frac{1}{3},\xi)}{ P(\frac{1}{3},\xi)}\right) <0. 
 \end{split}
\end{equation}
\noindent As $-\tilde{c}_0 \in  [-\pi^{\frac23}, -\left(  \frac{3\pi}{4}\right)^{\frac23} ]$, we deduce (\ref{eq_interval_tc2j}) by counting the constants $-\tilde{c}_{2j}$ and the intervals in which $u$ vanishes.
\end{demo}
\vskip 3mm

\noindent We finish this Appendix with the proof of Lemma \ref{lem_asymp_ck}.
\vskip 3mm

\begin{demo}(of Lemma \ref{lem_asymp_ck}).
\noindent Let $j\geq 0$. Applying (\ref{eq_u_asym}) with $x=\tilde{c}_{2j}$,
{\small $$  2 \pi^{\frac12} \tilde{c}_{2j}^{-\frac14} Ai'(0) \left(  \sin(\tilde{\xi}_{2j} - \tfrac{7\pi}{12})P(\tfrac{1}{3},\tilde{\xi}_{2j}) + \cos(\tilde{\xi}_{2j} - \tfrac{7\pi}{12})Q(\tfrac{1}{3},\tilde{\xi}_{2j}) \right) = u(-\tilde{c}_{2j}) =0.$$}
Thus, using Lemma \ref{lem_PQ} and $\tilde{\xi}_{2j}>\frac{1}{\sqrt{13}}$ (thanks to $\tilde{a}_0 >\left(\frac{3}{2\sqrt{13}}\right)^{\frac23}$),  we have  
\begin{equation}\label{eq_tan_PQ_1}
P(\tfrac{1}{3},\tilde{\xi}_{2j})>0 \quad \mbox{ and }\quad \tan\left(\tilde{\xi}_{2j} - \tfrac{7\pi}{12}\right) = \frac{ \sin(\tilde{\xi}_{2j} - \frac{7\pi}{12})}{ \cos(\tilde{\xi}_{2j} - \frac{7\pi}{12})} = -\frac{Q(\frac{1}{3},\tilde{\xi}_{2j})}{P(\frac{1}{3},\tilde{\xi}_{2j})}. 
\end{equation}
\noindent With (\ref{eq_ineq_QP13}) and (\ref{eq_tan_PQ_1}) we get:
\begin{equation}\label{eq_tan_PQ_3}
\left| \tan\left(\tilde{\xi}_{2j} - \tfrac{7\pi}{12}\right) \right| < \frac{5}{36\cdot \tilde{\xi}_{2j}}. 
\end{equation}
From (\ref{eq_interval_tc2j}) and (\ref{eq_tan_PQ_3}), it leads to 
\begin{equation}\label{eq_ineq_txi2j}
\left| \tilde{\xi}_{2j}-\frac{7\pi}{12} -j \pi  \right| <  \mbox{ arctan}\left(  \frac{5}{36\cdot \tilde{\xi}_{2j}} \right) \leq  \frac{5}{36\cdot \tilde{\xi}_{2j}} \leq \frac{5}{36( j\pi + \frac{\pi}{2})},
\end{equation}
which proves the assertion on the interval of localization of $\tilde{\xi}_{2j}$. 

\vskip 3mm
\noindent Since $\tilde{a}_0 >\left(\frac{3}{2\sqrt{11}}\right)^{\frac23}$, $\xi_{2j}>\frac{1}{\sqrt{11}}$ and, using Lemma \ref{lem_PQ}, we have $P(\frac{2}{3},\xi_{2j})>0$.
\vskip 3mm

\noindent Applying (\ref{eq_vp_asym}) with $x=c_{2j}$,
{\small $$ -2 \pi^{\frac12} c_{2j}^{\frac14} Ai(0) \left( - \cos(\xi_{2j} + \tfrac{\pi}{12})P(\tfrac{2}{3},\xi_{2j}) + \sin(\xi_{2j} + \tfrac{\pi}{12})Q(\tfrac{2}{3},\xi_{2j})  \right) = v'(-c_{2j}) = 0.$$}
Thus,
\begin{equation}\label{eq_tan_PQ_2}
 \frac{Q(\frac{2}{3},\xi_{2j})}{ P(\frac{2}{3},\xi_{2j})} = \mbox{cotan}\left(\xi_{2j}+\tfrac{\pi}{12} \right) = - \tan\left(\xi_{2j} + \tfrac{\pi}{12}-\tfrac{\pi}{2} \right) =- \tan\left(\xi_{2j} - \tfrac{5\pi}{12} \right). 
\end{equation}
Using (\ref{eq_ineq_QP23}) and (\ref{eq_tan_PQ_2}) we get:
\begin{equation}\label{eq_tan_PQ_4}
\left| \tan\left(\xi_{2j} - \tfrac{5\pi}{12}\right) \right| < \frac{7}{12\cdot \xi_{2j}},
\end{equation}
and the rest of the proof of the interval of localization of $\xi_{2j}$ is similar to what we have done for $\tilde{\xi}_{2j}$, thanks to  (\ref{eq_interval_c2j}). The intervals of localization of $\xi_{2j+1}$ and $\tilde{\xi}_{2j+1}$ are obtained in a similar way.
\vskip 3mm

\noindent In order to prove (\ref{eq_approx_c2j}), we need a more precise estimate on $\xi_{2j}$. Using (\ref{eq_tan_PQ_2}) and  \cite[9.2.9 and 9.2.10]{AS},
{\small \begin{align}
 \tan\left( \xi_{2j}-\tfrac{5\pi}{12} \right) & =  \frac{\frac{7}{72\xi_{2j}}-\frac{7\times 65\times 209}{6\times 8^3 \times 9^3 \xi_{2j}^3 } + \mathcal{O}\left(\frac{1}{\xi_{2j}^5 } \right) }{1-\frac{7\times 65}{2\times 8^2 \times 9^2 \xi_{2j}^2 } + \mathcal{O}\left(\frac{1}{\xi_{2j}^4 } \right) } \nonumber \\
  & =  \left( \mfrac{7}{72\xi_{2j}}-\tfrac{7\times 65\times 209}{6\times 8^3 \times 9^3 \xi_{2j}^3 } + \mathcal{O}\left(\mfrac{1}{\xi_{2j}^5 } \right) \right) \left(1+\tfrac{7\times 65}{2\times 8^2 \times 9^2 \xi_{2j}^2 } + \mathcal{O}\left(\mfrac{1}{\xi_{2j}^4 } \right) \right) \nonumber \\
& = \mfrac{7}{72\xi_{2j}}  + \mathcal{O}\left(\mfrac{1}{\xi_{2j}^3 } \right) \nonumber.
\end{align}}
Thus, thanks to $\xi_{2j}-\frac{5\pi}{12}-j\pi \in [-\frac{\pi}{12},\frac{\pi}{12} ]$,
$$\xi_{2j}-\mfrac{5\pi}{12}-j\pi =  \mbox{arctan}\left(\mfrac{7}{72\xi_{2j}} +  \mathcal{O}\left(\mfrac{1}{\xi_{2j}^3 } \right)  \right)  =  \mfrac{7}{72\xi_{2j}} \left( 1+\mathcal{O}\left(\mfrac{1}{\xi_{2j}^2 }  \right) \right). $$
Since $\xi_{2j}\in \left[ j\pi + \frac{\pi}{3}, j\pi + \frac{\pi}{2} \right]$, there exists a constant $C>0$ such that:
$$ \frac{7}{72(j\pi + \frac{\pi}{2})} \left( 1 - \frac{C}{j^2}\right)   \leq \xi_{2j}-\frac{5\pi}{12}-j\pi \leq \frac{7}{72(j\pi + \frac{\pi}{3})} \left( 1 + \frac{C}{j^2}\right)$$
which proves that
$$\xi_{2j} =\frac{5\pi}{12} +j \pi + \frac{7}{72(j\pi + \frac{\pi}{2})} + \mathcal{O}\left( \frac{1}{j^3} \right).$$
Since $c_{2j}=\left(\frac32 \xi_{2j} \right)^{\frac23 }$, 
{\small \begin{align*}
c_{2j} & =  \left( \frac{3j\pi}{2} + \frac{5\pi}{8} +  \frac{7}{48(j\pi + \frac{\pi}{2})}  + \mathcal{O}\left( \frac{1}{j^3}\right) \right)^{\frac23 }  \\
 & =  \left( \frac{3j\pi}{2} \right)^{\frac23 } \cdot \left( 1+  \frac{5}{12j} +  \frac{7}{72j\pi(j\pi + \frac{\pi}{2})} + \mathcal{O}\left( \frac{1}{j^4}\right) \right)^{\frac23 }  \\
 & =  \left( \frac{3j\pi}{2} \right)^{\frac23 } \cdot \left( 1+  \frac{5}{18j} +  \frac{7}{48j\pi(j\pi + \frac{\pi}{2})} - \frac19 \left( \frac{5}{12j} \right)^2 + \mathcal{O}\left( \frac{1}{j^4}\right) \right) \\
 & =   \left( \frac{3j\pi}{2} \right)^{\frac23 } \cdot \left( 1+  \frac{5}{18j} + \mathcal{O}\left( \frac{1}{j^2}\right) \right), 
\end{align*}}
which proves (\ref{eq_approx_c2j}).

\noindent Similarly, using $\tilde{\xi}_{2j}\in \left( j\pi + \tfrac{\pi}{2}, j\pi + \tfrac{2\pi}{3} \right)$, $\xi_{2j+1}\in \left( j\pi + \tfrac{5\pi}{6}, j\pi + \pi \right)$ and  $\tilde{\xi}_{2j+1}\in \left( j\pi + \pi, j\pi + \tfrac{7\pi}{6} \right),$
one proves (\ref{eq_approx_tc2j}), (\ref{eq_approx_c2j1}) and (\ref{eq_approx_tc2j1}).
\end{demo}

\section{Variations of $u$ and $v$}\label{sec_app_variations}

{\footnotesize\begin{table}[H]
\begin{tikzpicture}
\tkzTabInit[lgt = 1.6, espcl = 1, deltacl = 0.5]{$x$ /1, $(-1)^j u'$ /1, $(-1)^j u$ /1.5, $(-1)^j v'$ /1, $(-1)^j v$ /1.5 } {$-c_{2j+2}$, $-\tilde{a}_{j+2}$ , $-\tilde{c}_{2j+1}$, $-c_{2j+1}$, $-a_{j+1}$ ,  $-\tilde{c}_{2j}$, $-c_{2j}$, $-\tilde{a}_{j+1}$ ,  $-\tilde{c}_{2j-1}$, $-c_{2j-1}$ }
\tkzTabLine{,  , - , , z , , ,  ,  ,  , + , ,   ,  , , ,z  ,   -, } 
\tkzTabVar{+/ $<0$ , R/ ,  -/$<0$ , R/ , R/, R/, R/, R/,  +/$>0$, -/ $>0$ }
\tkzTabIma{3}{9}{6}{$0$} 
\tkzTabLine{z, ,  ,  ,  , , - , , , ,  ,  , z , , , + , , } 
\tkzTabVar{+/$>0$, R/ , R/, R/  , R/ , R/, -/$<0$,  R/, R/ , +/ $0$ }
\tkzTabIma{1}{7}{4}{$0$} 
\end{tikzpicture}
\caption{Variations of $(-1)^j u$ and $(-1)^j v$ on the interval $[-c_{2j+2}, -c_{2j-1}]$ for $j\geq 1$.}\label{tab_var_uv}

\end{table}}

\begin{table}[H]
\begin{center}
$$\begin{array}{|c|c|c|c|c|c|c|c|c|c|c|c|}
\hline c_0 & c_1 & c_2 & c_3 & c_4 & c_5 & c_6 & c_7 & c_8 & c_9 & c_{10}  \\
\hline 1.515 & 2.66 & 3.53 & 4.34 & 5.06 & 5.74 & 6.37 & 6.98 & 7.56 & 8.13 & 8.67  \\
\hline
\end{array}$$

\caption{Approximate values of $c_p$ for $p=0$ to $p=10$.}\label{tab_values_cct}
\end{center}
\end{table}



\begin{figure}[H]
\centering
\resizebox{0.9\textwidth}{!}{%
  \includegraphics{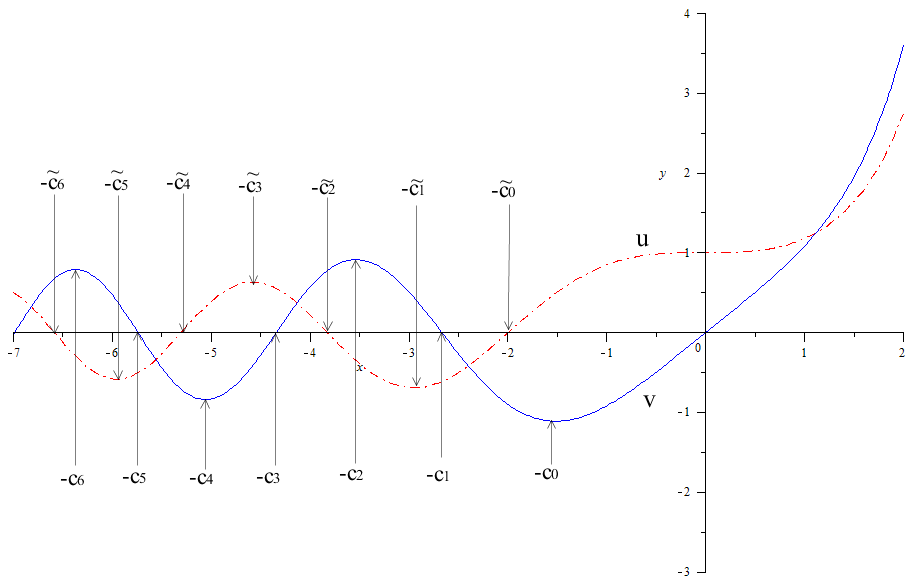}
}
\caption{Graphs of $u$, $v$ and positions of the $-c_p$ and $-\tilde{c}_p$.}
\label{fig_uvcct}       
\end{figure}

\section{Sturm-Picone's lemmas}\label{sec_sturm_picone}

In this Appendix we prove a Sturm's formula and a version of the Sturm-Picone's formula adapted to the setting of the proof of Lemma \ref{lem_Eomax2}. 

\begin{lem}\label{lem_S}
Let $a<b$ be two real numbers and let $g_1,g_2 \in C^0([a,b])$. Let $z$ be a solution of $-z'' +g_1 z=0$ and let $y$ be a solution of $-y''+g_2 y=0$.
Then:
\begin{equation}\label{eq_S_formula}
\left( yz'-zy' \right)' = (g_1 -g_2)yz.
\end{equation} 
\end{lem}

\begin{demo}
We have:
$$y(-z''+g_1 z) - z(-y'' + g_1 y)  =  -yz'' + zy'' =  - (yz'-zy')' $$
We also have:
$$(-y''+g_1y)-(-y'' + g_2 y) = (g_1 -g_2 )y.$$
Then,
$$y(-z''+g_1 z)-z\left(-z'' + g_2z + (g_1-g_2)y \right) =  - (yz'-zy')'.$$
Since  $-z'' +g_1 z=0$ and $-y''+g_2 y=0$, we finally have:
$$- (yz'-zy')' = -(g_1-g_2)yz,$$
which proves (\ref{eq_S_formula}). 
\end{demo}

\begin{lem}\label{lem_SP}
Let $a<b$ be two real numbers and let $q_1,q_2,g \in C^0([a,b])\cap C^1((a,b))$, $q_1>q_2>0$. Let $z$ be a solution of $-(q_1 z')' +gz=0$ and let $y$ be a solution of $-(q_2 y')'+gy=0$ with $y>0$ on $(a,b]$.
Then:
\begin{equation}\label{eq_SP_formula}
 \left( \frac{z}{y} (q_1 y z' - q_2 y' z) \right)' = (q_1-q_2)(z')^2 + q_2 \left( z' - \frac{y' z}{y} \right)^2.
\end{equation}
Moreover, if there exists $\eta >0$ such that $q_1-q_2 >\eta$, then there exists $A >0$ such that
\begin{equation}\label{eq_SP_formula2}
\left[ \frac{z}{y} (q_1 y z' - q_2 y' z) \right]_a^b \geq A \int_a^b z^2(x) \dd x.
\end{equation}
\end{lem}

\begin{demo} We have, on the interval $(a,b)$,
\begin{align*}
 \left( \frac{z}{y} (q_1 y z' - q_2 y' z) \right)' & =   \left( q_1 z z' - q_2 \frac{y'}{y} z^2 \right)'  \\
 & =  (q_1 z')' z + q_1 (z')^2 - (q_2 y')' \frac{z^2}{y} - q_2 y' \left( \frac{z^2}{y}\right)'  \\
& =  gz^2 + q_1(z')^2 - g y \frac{z^2}{y} - q_2 y' \left( \frac{z^2}{y}\right)'  \\
& =  q_1 (z')^2 - q_2 y' \left( \frac{z^2}{y}\right)'  \\
& =  (q_1 -q_2)(z')^2 + q_2 \left( (z')^2 -  y' \left( \frac{z^2}{y}\right)'\right)\\
& =   (q_1 -q_2)(z')^2 + q_2 \left( (z')^2 - 2 y' \frac{ z z'}{y} + \frac{(y')^2}{y^2} z^2 \right)  \\
& =  (q_1 -q_2)(z')^2 + q_2 \left( z' - \frac{y' z}{y} \right)^2 . 
\end{align*}
This proves (\ref{eq_SP_formula}). Then, integrating (\ref{eq_SP_formula}) between $a$ and $b$ and using Poincar\'e inequality in the last inequality, there exists $A >0$ (depending on $\eta, a$ and $b$) such that
{\small \begin{align*}
 \left[ \frac{z}{y} (q_1 y z' - q_2 y' z) \right]_a^b & =  \int_a^b (q_1(x)-q_2(x))(z'(x))^2 + q_2(x) \left( z'(x) - \frac{y'(x) z(x)}{y(x)} \right)^2 \dd x  \\
 & \geq  \int_a^b (q_1(x)-q_2(x)) (z'(x))^2 \dd x \\
& \geq  \eta \int_a^b (z'(x))^2 \dd x  \geq A \int_a^b (z(x))^2 \dd x. 
\end{align*}}
This proves (\ref{eq_SP_formula2}).
\end{demo}

\section{The monotonicity arguments}\label{sec_zk}

\noindent We have defined the functions $f_x$, $g_x$ and the functions $z_k$ for $k\geq 0$ in Section \ref{sec_2analytic}.

\begin{lem}\label{lem_zeroes_appendix}
Let $k\geq 0$. Then, for every $x\geq 0$,
$$z_k(x)\geq 0\quad \mbox{ and }\quad z_k(x)=x-\psi_k(x).$$
Therefore, $z_k$ is continuous on $[0,+\infty)$. Moreover, for every $j\geq 0$ and every $x\geq 0$, 
\begin{equation}
\begin{split}
  0< x+\tilde{a}_1 < z_0(x) & \leq x+c_0 <  \cdots  < x+ \tilde{a}_{j+1} < z_{2j}(x) \leq x+c_{2j}  \\  
&  <  x+a_{j+1} < z_{2j+1}(x) \leq x+c_{2j+1} \ldots \label{eq_ineg_zk}
\end{split}
\end{equation}

\end{lem}

\begin{demo}
In this proof it will be easier to use the expressions in terms of Airy functions for $f_x$ and $g_x$ since we will use classical properties of the $Ai$ and $Bi$ functions and in particular the fact that $Ai'$ is strictly negative on the positive real half-line, which is not the case for $u'$.

\noindent For $x\geq 0$, $Ai(x)>0$, $Bi(x)>0$, $Bi'$ is strictly positive on $[0,+\infty)$ and $Ai'$ is strictly negative on $[0,+\infty)$. If $z\leq 0$, $x-z \geq 0$ and  $f_x(z)>0$. Thus, $z_0(x) >0$. Then, $0$ is a zero of $g_x$ and since $\frac{Bi}{Ai}$ is strictly increasing on $[0,+\infty)$, and for $z<0$, $x-z > x$ and $g_x(z)>0$. So, $0$ is the first zero of $g_x$. In particular, for every $k\geq 0$, $z_k(x) \geq 0$.
\vskip 2mm

\noindent Let $j\geq 0$. We remark that, by definition of $\psi_{2j}$, we have $f_x(x-\psi_{2j}(x))=0$, and by definition of $\psi_{2j+1}$, we have $g_x(x-\psi_{2j+1}(x))=0$. Moreover, for $x-z \notin \{ -\tilde{c}_{2j+1} \}_{j\geq 0}$, by unicity of $\psi_{2j}(x)$ in $[-c_{2j},-\tilde{a}_{j+1})$, $x-\psi_{2j}(x)$ is the unique zero of $f_x$ in $(x+\tilde{a}_{j+1},x+c_{2j}]$. Since we have $f_x(x+\tilde{a}_{j+1})=Bi'(-\tilde{a}_{j+1}) Ai(x) \neq 0$, the set of the zeroes of $f_x$ is exactly $\{ x-\psi_{2j}(x),\ j\geq 0 \}$. Thus, for every $j\geq 0$, $z_{2j}(x)=x-\psi_{2j}(x)$.
\vskip 2mm

\noindent For $x-z \notin \{ -\tilde{c}_{2j} \}_{j\geq 0}$, by unicity of $\psi_{2j+1}(x)$ in $[-c_{2j+1},-a_{j+1})$, $x-\psi_{2j+1}(x)$ is the unique zero of $g_x$ in $(x+a_{j+1},x+c_{2j+1}]$. Since we have $g_x(x+a_{j+1})=Bi(-a_{j+1}) Ai(x) \neq 0$, the set of the zeroes of $g_x$ is exactly $\{ 0 \} \cup \{ x-\psi_{2j+1}(x),\ j\geq 0 \}$. Thus, for every $j\geq 0$, $z_{2j+1}(x)=x-\psi_{2j+1}(x)$.
\vskip 2mm

\noindent By Lemma \ref{lem_psi_j}, we deduce that $z_k$ is continuous on $[0,+\infty)$. Since for every $j\geq 0$, $\psi_{2j}(x)\in [-c_{2j},-\tilde{a}_{j+1})$ and $\psi_{2j+1}(x) \in [-c_{2j+1},-a_{j+1})$, we deduce (\ref{eq_ineg_zk}). 
\end{demo}

\noindent We now prove further properties of the functions $f_x$ and $g_x$ and in particular their signs and their variations.

\vskip 2mm

\begin{prop}\label{prop_fx_edo}
 For every $x\geq 0$, the functions $f_x$ and $g_x$ from $\R$ to $\R$ have the following properties:
\begin{enumerate}
\item $\forall z\in \R,\ g_x'(z)=-f_x(z)$ and $f_x'(z)=-(x-z)g_x(z)$.
\item $f_x$ satisfies the ordinary differential equation on $\R\setminus \{ x\}$:
\begin{equation}\label{eq_lem_Eomax2_1}
 \left( \frac{f_x'}{x-z} \right)' = f_x,
\end{equation}
and $g_x$ satisfies the Airy equation: $g_x''=(x-z)g_x$.
\end{enumerate}
\end{prop}

\begin{demo}
\noindent \textbf{(1)} We compute the derivative of $f_x$, using the fact that $u$ and $v$ satisfies the Airy equation:
\begin{equation}\label{eq_prop_fx_1}
\forall z\in \R,\ f_x'(z)= (x-z)(u(x-z)v(x)-v(x-z)u(x)). 
\end{equation}
Thus, $f_x'(z)=-(x-z)g_x(z)$, for every $z\in \R$. A direct computation leads to  $g_x'(z)=-f_x(z)$, for every $z\in \R$.

\noindent \textbf{(2)} We assume that $z\neq x$, we divide (\ref{eq_prop_fx_1}) by $x-z$ and by derivation: 
$$\forall z\in \R\setminus \{x\},\ \left( \frac{f_x'}{x-z} \right)' = -u'(x-z)v(x) + v'(x-z) u(x) = f_x(z).$$
The function $g_x$ satisfies the Airy equation since it is a linear combination of solutions of the Airy equation.
\end{demo}

\begin{prop}\label{prop_fx}
For every $x\geq 0$, the functions $f_x$ and $g_x$ from $\R$ to $\R$ have the following properties:
\begin{enumerate}
\item The function $f_x'$ vanishes exactly on $0$, $x$, and $z_{2j+1}(x)$ for every $j\geq 0$. It is strictly negative on $(-\infty,0)$, strictly positive on $(0,x)$, strictly negative on $(x,z_1(x))$ and, for every $j\geq 1$, $(-1)^{j+1} f_x'$ is strictly positive on 
$(z_{2j-1}(x),z_{2j+1}(x))$.

\item The function $f_x$ is strictly positive on $(-\infty,z_0(x))$ and, for every $j\geq 1$,  $(-1)^{j+1} f_x$  is strictly positive on $(z_{2j-2}(x),z_{2j}(x))$.
 
\item The function $g_x'$ vanishes exactly on $z_{2j}(x)$ for every $j\geq 0$. It is strictly negative on $(-\infty,z_0(x))$ and, for every $j\geq 1$,  $(-1)^{j+1} g_x'$  is strictly positive on $(z_{2j-2}(x),z_{2j}(x))$.

\item The function $g_x$ is strictly positive on $(-\infty,0)$, strictly negative on $(0,z_1(x))$ and, for every $j\geq 1$, $(-1)^{j+1} g_x$ is strictly positive on $(z_{2j-1},z_{2j+1})$.
\end{enumerate}
\end{prop}

\vskip 3mm

\begin{demo}
\noindent We will again use the expressions in terms of Airy functions for $f_x$ and $g_x$. 

\noindent \textbf{(1)} From (\ref{eq_prop_fx_1}), it is clear that $f_x'(0)=f_x'(x)=0$. We have already proven in Lemma \ref{lem_zeroes} that for $z\in (-\infty,0)$, $f_x'(z)<0$. Then, for $z\in(0,x)$, $x-z > 0$, $x-z < x$, $\frac{Bi}{Ai}(x-z) < \frac{Bi}{Ai}(x-z)$ and $f_x'(z)>0$. From $f_x'(z)=-(x-z)g_x(z)$ and Lemma \ref{lem_zeroes}, we know that the remaining zeroes of $f_x'$ are exactly the $z_{2j+1}(x)$ for $j\geq 0$.  We also have $f_x'(x+a_1)=a_1 Bi(-a_1)Ai(x)$ with $a_1>0$, $Bi(-a_1)<0$ and $Ai(x)>0$, thus  $f_x'(x+a_1) <0$. Since $f_x'$ is of constant sign in $(x,z_1(x))$, one deduce that $f_x'$ is strictly negative on $(x,z_1(x))$. To finish the proof of point $(1)$, it is sufficient to remark that $f_x'$ is of constant sign on every interval $(z_{2j-1}(x),z_{2j+1}(x))$ for $j\geq 1$. But, $x+a_{2j} \in (z_{4j-3}(x),z_{4j-1}(x))$ and $f_x'(x+a_{2j})=a_{2j} Bi(-a_{2j})Ai(x) >0$, since $Bi(-a_{2j})>0$. Thus, $f_x'$ is strictly positive on $(z_{4j-3}(x),z_{4j-1}(x))$. Similarly,  $x+a_{2j+1} \in (z_{4j-1}(x),z_{4j+1}(x))$ and $f_x'(x+a_{2j+1})=a_{2j+1} Bi(-a_{2j+1})Ai(x) <0$,  since $Bi(-a_{2j+1})<0$. Thus, $f_x'$ is strictly negative on $(z_{4j-1}(x),z_{4j+1}(x))$.
\vskip 2mm

\noindent \textbf{(2)} We have already proven in Lemma \ref{lem_zeroes} that for $z\in (-\infty,0)$, $f_x(z)>0$. We also have  $f_x(0)=\frac{1}{\pi}>0$ since it is the value of the Wronskian of $Ai$ and $Bi$ and thus, for every $z\in (-\infty,x)$, $f_x(z) \geq \frac{1}{\pi}$. Since $z_0(x)$ is the first zero of $f_x$,  this function is strictly positive on $(-\infty,z_0(x))$. We remark that $f_x$ is of constant sign on every interval $(z_{2j-2}(x),z_{2j}(x))$ for $j\geq 1$. But, $x+\tilde{a}_{2j+1} \in (z_{4j-2}(x),z_{4j}(x))$ and $f_x(x+\tilde{a}_{2j+1})=Bi'(-\tilde{a}_{2j+1})Ai(x) >0$, since $Bi'(-\tilde{a}_{2j+1})>0$. Thus, $f_x$ is strictly positive on $(z_{4j-2}(x),z_{4j}(x))$. Similarly,  $x+\tilde{a}_{2j+2} \in (z_{4j}(x),z_{4j+2}(x))$ and $f_x(x+\tilde{a}_{2j+2})=Bi'(-\tilde{a}_{2j+2})Ai(x) <0$,  since $Bi'(-\tilde{a}_{2j+2})<0$. Thus, $f_x$ is strictly negative on $(z_{4j}(x),z_{4j+2}(x))$.
\vskip 2mm

\noindent \textbf{(3)} It is deduced directly from point $(1)$ of Proposition \ref{prop_fx_edo} and point $(2)$. 
\vskip 2mm

\noindent \textbf{(4)} It comes from point $(1)$ of Proposition \ref{prop_fx_edo}, point $(1)$ and the fact that for $z\geq z_1(x)$, $z>x$ and $x-z<0$.
\end{demo}

\vskip 5mm

\noindent We have now all the ingredients needed to prove that $z_k$ is a strictly increasing function.

\begin{lem}\label{lem_Eomax2}
For every $k\geq 0$, the function $z_k$ is strictly increasing from $[0,+\infty)$ to $[c_k,+\infty)$.  
\end{lem}

\begin{demo} We will separate the proof in two cases, depending on the parity of $k$.
\vskip 2mm

\noindent \textbf{Case 1: $k=2j$} for $j\geq 0$. Let $0< x_1 < x_2$.  We want to prove that $z_{2j}(x_1)\leq z_{2j}(x_2)$. Assume that $ z_{2j}(x_2) < z_{2j}(x_1)$. Let $\delta >0$ be such that $z_{2j}(x_2) +\delta < z_{2j}(x_1)$. We use (\ref{eq_ineg_zk}) to get
$$z_{2j-1}(x_1) < x_1+\tilde{a}_{j+1} < x_2+\tilde{a}_{j+1} < z_{2j}(x_2) < z_{2j}(x_2)+\delta <\qquad \qquad$$
$$\qquad \qquad \qquad\qquad  z_{2j}(x_1) < x_1 + a_{j+1} < x_2 +a_{j+1} < z_{2j+1}(x_2).$$
In particular, $x_1 - (z_{2j}(x_2)+\delta) <0$, $x_2 - (z_{2j}(x_2)+\delta) <0$ and $$(z_{2j}(x_2)+\delta) \in (z_{2j-1}(x_1), z_{2j}(x_1)) \cap (z_{2j}(x_2), z_{2j+1}(x_2)).$$ Thus, using Proposition \ref{prop_fx},
\begin{equation}\label{eq_demo_Eomax31}
(-1)^j f_{x_1}(z_{2j}(x_2)+\delta) >0,\quad (-1)^j f_{x_1}'(z_{2j}(x_2)+\delta) <0,
\end{equation}
\begin{equation}\label{eq_demo_Eomax32}
\quad \mbox{ and }\qquad (-1)^j f_{x_2}(z_{2j}(x_2)+\delta) <0,  \quad (-1)^j f_{x_2}'(z_{2j}(x_2)+\delta) <0.
\end{equation}
There exists $\eta >0$ such that, for every $z\in [z_{2j}(x_2)+\delta , z_{2j}(x_1)]$, $\frac{1}{x_1 -z} - \frac{1}{x_2 -z} \geq \eta$. Moreover, since $(z_{2j}(x_2) +\delta , z_{2j}(x_1)) \subset (z_{2j-1}(x_1) , z_{2j}(x_1))$, for every $z\in (z_{2j}(x_2) +\delta , z_{2j}(x_1)) $, $(-1)^{j+1}f_{x_1}(z)>0$. Similarly, since we have the inclusion $(z_{2j}(x_2) +\delta , z_{2j}(x_1)) \subset (z_{2j}(x_2),z_{2j+1}(x_2)) $, for every $z\in (z_{2j}(x_2) +\delta , z_{2j}(x_1)) $, $(-1)^{j+1} f_{x_2}(z)<0$. Then, applying Lemma \ref{lem_SP},
\begin{equation}\label{eq_demo_Eomax3_SP}
\int_{z_{2j}(x_2)+\delta}^{z_{2j}(x_1)} \left( \frac{f_{x_1}(z)}{f_{x_2}(z)} \left( \frac{f_{x_1}'(z)f_{x_2}(z)}{x_1 -z}- \frac{f_{x_1}(z)f_{x_2}'(z)}{x_2 -z} \right) \right)' \dd z  >0.
\end{equation}
Since $f_{x_1}(z_{2j}(x_1))=0$, the integral in the left side of equality (\ref{eq_demo_Eomax3_SP}) is equal to
\begin{equation}\label{eq_demo_Eomax4}
-  \frac{ f_{x_1}(z_{2j}(x_2)+\delta) f_{x_1}'(z_{2j}(x_2)+\delta)}{x_1 -z_{2j}(x_2) -\delta} + \frac{ f_{x_1}^2(z_{2j}(x_2)+\delta)}{x_2 -z_{2j}(x_2) -\delta }\cdot  \frac{f_{x_2}'(z_{2j}(x_2)+\delta)}{f_{x_2}(z_{2j}(x_2)+\delta)} <0
\end{equation}
by the use of (\ref{eq_demo_Eomax31}) and (\ref{eq_demo_Eomax32}). But (\ref{eq_demo_Eomax4}) contradicts (\ref{eq_demo_Eomax3_SP}) and thus we must have $ z_{2j}(x_1) \leq  z_{2j}(x_2)$. The function $z_{2j}$ is an increasing function from $[0,+\infty)$ to $[c_k,+\infty)$.
\vskip 2mm

\noindent It remains to prove that $z_{2j}$ is strictly increasing. If $z_{2j}$ is not strictly increasing, since it is increasing and continuous, there exists an interval in $[0,+\infty)$ on which $z_{2j}$ is constant. But, $z_{2j}$ is also analytic on $[0,+\infty)$ since one can prove that actually the functions $\psi_{2j}$ are analytic. Thus, if it is constant on an interval, it should be constant everywhere which is not the case, so $z_{2j}$ is actually strictly increasing. 
\vskip 3mm

\noindent \textbf{Case 2: $k=2j+1$} for $j\geq 0$. Let $x_1 < x_2$. We will show by induction on $j\geq 0$ that $z_{2j+1}(x_1) < z_{2j+1}(x_2)$.

\noindent For $j=0$, we directly apply the classical interlacing zeroes theorem of Sturm with potentials $q(z)= -(x_2 -z) < p(z)=-(x_1 -z)$, since $g_{x_1}$ satisfies $-g_{x_1}'' + pg_{x_1}=0$ and  $g_{x_2}$ satisfies $-g_{x_2}'' + qg_{x_2}=0$. Applying this theorem between $0$ which is a common zero to $g_{x_1}$ and $g_{x_2}$ and $z_1(x_2)$ which is the first strictly positive zero of $g_{x_2}$ one gets that $g_{x_1}$ admits a zero in the interval $(0,z_1(x_2))$. Since $z_1(x_1)$ is the smallest strictly positive zero of $g_{x_1}$, we necessarily have $z_1(x_1)\in (0,z_1(x_2))$ and $z_1(x_1) < z_1(x_2)$. Thus, $z_1$ is strictly increasing. 
\vskip 2mm

\noindent Let $j\geq 1$. We assume by induction that $z_{2j-1}(x_1) < z_{2j-1}(x_2)$ and we want to prove that $z_{2j+1}(x_1)< z_{2j+1}(x_2)$. We assume the contrary: $z_{2j+1}(x_2) \leq z_{2j+1}(x_1)$. Then we have
\begin{equation}\label{eq_demo_Eomax6}
z_{2j-1}(x_1) < z_{2j-1}(x_2) < z_{2j}(x_2) < z_{2j+1}(x_2) \leq z_{2j+1}(x_1). 
\end{equation}
We apply Lemma \ref{lem_S} to $g_{x_1}$ and $g_{x_2}$ between $z_{2j-1}(x_2)$ and $z_{2j+1}(x_2)$ to get
{\small \begin{equation}\label{eq_demo_Eomax7_S}
 \int_{z_{2j-1}(x_2)}^{z_{2j+1}(x_2)} \left( g_{x_2}(z)g_{x_1}'(z)-g_{x_1}g_{x_2}'(z) \right)' \dd z = \int_{z_{2j-1}(x_2)}^{z_{2j+1}(x_2)} (x_1 - x_2) g_{x_1}(z)g_{x_2}(z) \dd z. 
\end{equation}}
But, using (\ref{eq_demo_Eomax6}), we have $(z_{2j-1}(x_2),z_{2j+1}(x_2)) \subset (z_{2j-1}(x_1),z_{2j+1}(x_1))$. Using Proposition \ref{prop_fx}, 
$$\forall z \in (z_{2j-1}(x_2),z_{2j+1}(x_2)),\ (-1)^j g_{x_1}(z) < 0\ \mbox{ and }\ (-1)^j g_{x_2}(z) <0.$$ 
Since $x_1-x_2 <0$,
\begin{equation}\label{eq_demo_Eomax8}
  \int_{z_{2j-1}(x_2)}^{z_{2j+1}(x_2)} (x_1 - x_2) g_{x_1}(z)g_{x_2}(z) \dd z < 0.
\end{equation}
We have, 
\begin{align*}
\int_{z_{2j-1}(x_2)}^{z_{2j+1}(x_2)} &\left( g_{x_2}(z)  g_{x_1}'(z)-g_{x_1}g_{x_2}'(z) \right)' \dd z = \\
&  -g_{x_1}(z_{2j+1}(x_2)) g_{x_2}'(z_{2j+1}(x_2))+ g_{x_1}(z_{2j-1}(x_2)) g_{x_2}'(z_{2j-1}(x_2)). 
\end{align*}

\noindent But, using again Proposition \ref{prop_fx},
$$\forall z \in [z_{2j-1}(x_2),z_{2j}(x_2)),\ (-1)^j g_{x_2}'(z) < 0$$
and
$$\forall z \in (z_{2j}(x_2),z_{2j+1}(x_2)],\ (-1)^j g_{x_2}'(z) > 0.$$
In particular, 
$$(-1)^j g_{x_1}(z_{2j+1}(x_2)) < 0,\quad  (-1)^j g_{x_2}'(z_{2j+1}(x_2)) >0,$$
and
$$(-1)^j g_{x_1}(z_{2j-1}(x_2)) <0,\quad (-1)^j g_{x_2}'(z_{2j-1}(x_2))<0.$$
Thus,
$$\int_{z_{2j-1}(x_2)}^{z_{2j+1}(x_2)} \left( g_{x_2}(z)g_{x_1}'(z)-g_{x_1}g_{x_2}'(z) \right)' \dd z >0$$
which contradicts (\ref{eq_demo_Eomax8}). So we have $z_{2j+1}(x_1)<z_{2j+1}(x_2)$ and $z_{2j+1}$ is strictly increasing.
\vskip 2mm

\noindent We have thus proven by induction that for every $j\geq 0$, $z_{2j+1}$ is strictly increasing from $[0,+\infty)$ to $[c_{2j+1},+\infty)$. This finishes the proof of Lemma \ref{lem_Eomax2}.
\end{demo}

\end{document}